\documentclass[10pt,oneside,english]{amsart}
\usepackage[T1]{fontenc}
\usepackage[latin9]{inputenc}
\usepackage{geometry}
\usepackage{amsthm}
\usepackage{amssymb}
\usepackage{esint}

\makeatletter

\providecommand{\tabularnewline}{\\}

\theoremstyle{plain}
\newtheorem{thm}{\protect\theoremname}[section]
  \theoremstyle{definition}
  \newtheorem{defn}[thm]{\protect\definitionname}
  \theoremstyle{plain}
  \newtheorem{prop}[thm]{\protect\propositionname}
  \theoremstyle{plain}
  \newtheorem{cor}[thm]{\protect\corollaryname}
  \theoremstyle{plain}
  \newtheorem{lem}[thm]{\protect\lemmaname}
  \theoremstyle{definition}
  \newtheorem{problem}[thm]{\protect\problemname}
  \theoremstyle{remark}
  \newtheorem{claim}[thm]{\protect\claimname}

\DeclareMathOperator{\supp}{supp}

\DeclareMathOperator{\esssup}{esssup}
\DeclareMathOperator{\essinf}{essinf}

\DeclareMathOperator{\cnt}{cent}

\DeclareMathOperator{\id}{id}
\DeclareMathOperator{\dsc}{dsc}
\DeclareMathOperator{\diam}{diam}
\DeclareMathOperator{\lip}{lip}

\DeclareMathOperator{\ldim}{\underline{\dim}}
\DeclareMathOperator{\bdim}{bdim}
\newcommand{\sqr}{{}^{{}^{}_\square}}

\usepackage{epsfig}
\usepackage{color}

\newcommand{\hide}[1]{}

\makeatother

\usepackage{babel}
  \providecommand{\claimname}{Claim}
  \providecommand{\corollaryname}{Corollary}
  \providecommand{\definitionname}{Definition}
  \providecommand{\lemmaname}{Lemma}
  \providecommand{\problemname}{Problem}
  \providecommand{\propositionname}{Proposition}
\providecommand{\theoremname}{Theorem}

\begin{document}

\title{Dynamics on fractals and fractal distributions}

\author{Michael Hochman}

\thanks{Partially supported by NSF grant 0901534 and ERC grant 306494.}


\begin{abstract}
We study fractal measures on Euclidean space through the dynamics
of ``zooming in'' on typical points. The resulting family of measures
(the ``scenery''), can be interpreted as an orbit in an appropriate
dynamical system which often equidistributes for some invariant distribution.
The first part of the paper develops basic properties of these limiting
distributions and the relations between them and other models of dynamics
on fractals, specifically to Z\"{a}hle distributions and Furstenberg's
CP-processes. In the second part of the paper we study the geometric
properties of measures arising in these contexts, specifically their
behavior under projection and conditioning on subspaces.
\end{abstract}
\maketitle
\tableofcontents{}

\section{\label{sec:Introduction}Introduction}

This work is a systematic study of a class of measures, called uniformly
scaling measures (USMs), and associated distributions on measures,
called fractal distributions (FDs), which capture the notion of self-similarity
(or ``fractality'') of a measure on Euclidean space in terms of
the dynamics of rescaling and translation. USMs were introduced abstractly
by Gavish \cite{Gavish09} though examples were studied earlier by
many authors, e.g. Patzschke and Z\"{a}hle \cite{PatzschkeZahle90},
Bandt \cite{Bandt1992,Bandt1997}, and Graf \cite{Graf95}. Fractal
distributions, which we define here, generalize Z\"{a}hle's scale-invariant
distributions \cite{Zahle88} and are very closely related to Furstenberg's
CP-processes \cite{Furstenberg70,Furstenberg08}. They may also be
viewed as ergodic-theoretic analogues of the scenery flow for sets
which was studied, among others, by Bedford and Fisher \cite{BedfordFisher96,BedfordFisher97,BedfordFisherUrbanski02}. 

The models we discuss here are sufficiently general so as to unify
the treatment of many examples of interest in fractal geometry and
dynamics, but at the same time are sufficiently structured that their
geometric behavior is far better than general measures.

In particular, this work was motivated by some recent results by Furstenberg
\cite{Furstenberg08}, Peres and Shmerkin \cite{PeresShmerkin08},
and Hochman and Shmerkin \cite{HochmanShmerkin09}, on the geometry
of measures which arise from certain combinatorial constructions or
as invariant measures for certain dynamics. The behavior of such measures
under projections and conditioning on subspaces was shown to be more
regular than that of general measures. One of our motivations for
the present work was to place these results in a more general framework
and to clarify some of the objects involved. We also obtain new results
for general measures by relating them to structured measures which
arise from them using a limiting procedure.

In the remainder of this introduction we present the main definitions
and results. More detailed discussion, examples and proofs are provided
in the subsequent sections.

\subsection{\label{sub:Standing-notation}Standing notation}

Throughout the paper $d$ will be a fixed integer dimension. Equip
$\mathbb{R}^{d}$ with the norm $\left\Vert x\right\Vert =\sup|x_{i}|$
and the induced metric. $B_{r}(x)$ is the closed ball of radius $r$
around $x$; in particular, $B_{1}(0)=[-1,1]^{d}$. We abbreviate
\[
B_{r}=B_{r}(0)
\]
Let $\lambda$ denote Lebesgue measure on $\mathbb{R}^{d}$ and $\delta_{a}$
the point mass at $a$. 

Let $\mathcal{M}=\mathcal{M}(\mathbb{R}^{d})$ denote the space of
Radon measures on $\mathbb{R}^{d}$ with the weak topology, i.e. the
weakest topology which makes the maps $\mu\mapsto\int fd\mu$ continuous
for every compactly supported continuous function $f:\mathbb{R}^{d}\rightarrow\mathbb{R}$.
For a measurable space $X$ (e.g. a topological space with the Borel
structure) let $\mathcal{P}(X)$ denote the space of probability measures
on $X$. When $X$ is a compact metric space we give $\mathcal{P}(X)$
the weak topology, which is also compact and metrizable. 

We reserve the term \emph{distribution} for members of $\mathcal{P}(\mathcal{M})$
and for probability measures on similarly ``large'' spaces, denoting
them by $P,Q,R$, and use the term \emph{measure }exclusively for
members of $\mathcal{M}$, which are denoted $\mu,\nu,\eta$ etc.
We generally use brackets to denote operations which produce distributions
from measures, e.g. the operations $\left\langle \mu\right\rangle _{U}$,
$\left\langle \mu\right\rangle _{x,T}$, $\left\langle \mu,x\right\rangle _{N}$
defined below. 

Write $\mu\sim P$ to indicate that $\mu$ is chosen randomly according
to the distribution $P$, and similarly $x\sim\mu$. We also write
$\mu\sim\nu$ (or $P\sim Q$) to indicate equivalence of measures
i.e. that $\mu,\nu$ (or $P,Q$) have the same null sets. Which of
these is intended will be clear from the context. 

If $\mu\in\mathcal{M}$ and $\mu(A)>0$ then $\mu|_{A}\in\mathcal{M}$
is the restricted measure on $A$, i.e.
\[
\mu|_{A}(B)=\mu(A\cap B)
\]
and, assuming $0<\mu(A)<\infty$, let $\mu_{A}\in\mathcal{M}$ be
normalized  version of $\mu|_{A}$, i.e.
\[
\mu_{A}(B)=\frac{1}{\mu(A)}\mu(A\cap B)
\]

Finally, define two (partial) normalization operations $*,\sqr:\mathcal{M}\rightarrow\mathcal{M}$
by 
\begin{eqnarray*}
\mu^{*} & = & \frac{1}{\mu(B_{1})}\mu\\
\mu^{\sqr} & = & (\frac{1}{\mu(B_{1})}\mu)|_{B_{1}}
\end{eqnarray*}
so that $\mu^{\sqr}=\mu_{B_{1}}$. The operations are defined on $\{\mu\in\mathcal{M}\,:\,\mu(B_{1})>0\}$,
which is a measurable subset of $\mathcal{M}$. We apply these operations
pointwise to sets and distributions, so if $\mathcal{A}\subseteq\mathcal{M}$
then $\mathcal{A}^{*}=\{\mu^{*}\,:\,\mu\in\mathcal{A}\}$, etc. In
particular,
\begin{eqnarray*}
\mathcal{M}^{*} & = & \{\mu\in\mathcal{M}\,:\,\mu(B_{1})=1\}\\
\mathcal{M}^{\sqr} & = & \{\mu\in\mathcal{M}\,:\,\mu\mbox{ is a probability measure on }B_{1}\}\\
 & \cong & \mathcal{P}(B_{1})
\end{eqnarray*}
 In the same way, $P^{*}$ is the push-forward of the distribution
$P$ through $\mu\mapsto\mu^{*}$, and similarly $P^{\sqr},$etc.

For $x\in\mathbb{R}^{d}$ let $T_{x}:\mathbb{R}^{d}\rightarrow\mathbb{R}^{d}$
denote the translation taking $x$ to the origin:
\[
T_{x}(y)=y-x
\]
and given a ball $B=B_{r}(x)$ we write $T_{B}$ for the orientation-preserving
homothety mapping $B$ onto $B_{1}$, i.e.
\[
T_{B}(y)=\frac{1}{r}(y-x)
\]
Finally, we define scaling operators on $\mathbb{R}^{d}$ by
\[
S_{t}(y)=e^{t}y
\]
Note the exponential time scale, which makes $S=(S_{t})_{t\in\mathbb{R}}$
into an additive $\mathbb{R}$-action on $\mathbb{R}^{d}$ by linear
transformations (i.e. $S_{s+t}=S_{s}S_{t}$).

For $f:\mathbb{R}^{d}\rightarrow\mathbb{R}^{d}$ and $\mu\in\mathcal{M}$
write $f\mu$ for the push-forward of $\mu$ through $f$, that is,
$f\mu(A)=\mu(f^{-1}A)$. Thus the operators $T_{x},S_{t}$ induce
translation and scaling operations on $\mathcal{M}$, given by
\begin{eqnarray*}
T_{x}\mu(A) & = & \mu(A+x)\\
S_{t}\mu(A) & = & \mu(e^{-t}A)
\end{eqnarray*}
We also append $*,\sqr$ to operations on measures to indicate post-composition
with the corresponding normalization operator:
\begin{eqnarray*}
S_{t}^{*}\mu & = & (S_{t}\mu)^{*}\\
S_{t}^{\sqr}\mu & = & (S_{t}\mu)^{\sqr}
\end{eqnarray*}
and similarly $T_{x}^{\sqr}$ etc. Then $S^{*}=(S_{t}^{*})_{t\in\mathbb{R}}$
is again an additive $\mathbb{R}$-action on the measurable set $\{\mu\in\mathcal{M}^{*}\,:\,0\in\supp\mu\}$,
and $S^{\sqr}=(S_{t}^{\sqr})_{t\in\mathbb{R}^{+}}$ is an additive
$\mathbb{R}^{+}$-action on $\{\mu\in\mathcal{M}^{\sqr}\,:\,0\in\supp\mu\}$.
Note that the restriction of measures to $B_{1}$ is not an invertible
operation, so the maps $S_{t}^{\sqr}$ are not invertible for $t\neq0$.
Note also that these actions are discontinuous at measures $\mu$
with $\mu(\partial B_{1})>0$.

See Section \ref{sub:Summary-of-notation} below for a summary of
the notation.

\subsection{\label{sub:Fractal-distributions}Fractal distributions}

As usual in dynamical systems theory, one may study a dynamical system
either globally, via the invariant sets or measures of the system,
or in terms of the behavior of individual orbits (for background in
ergodic theory see Section \ref{sec:Ergodic-theory}, or \cite{Walters82,Glasner2003}
for more comprehensive introductions). While individual orbits are
often of most interest, it is usually easier to obtain global results
or results for typical orbits.

We shall similarly have two perspectives in our study of fractals:
a global one, which is concerned with distributions (on measures)
which possess certain invariance properties, and an individual one,
which deals with individual measures whose ``orbits'' display some
regularity. 

Our ``global'' objects are distributions on measures which exhibit
certain invariance under change of scale and translation, mirroring
the idea of ``self-similarity'' of a measure. In this section we
describe these objects axiomatically.

Recall that a probability distribution $P$ on $\mathcal{M}^{*}$
is $S^{*}$ invariant if $P((S_{t}^{*})^{-1}A)=P(A)$ for all measurable
sets $A\subseteq\mathcal{M}^{*}$ and all $t\in\mathbb{R}$. If $P$
is an $S^{*}$-invariant distribution then $P^{\sqr}$ (the push-forward
of $P$ through $\mu\mapsto\mu^{\sqr}$) is an $S^{\sqr}$-invariant
distributions called the \emph{restricted version }of $P$, and $P$
is the \emph{extended version} of $P^{\sqr}$. This is a 1-1 correspondence
between $S^{*}$- and $S^{\sqr}$-invariant distributions (see Lemma
 \ref{pro:restriction-of-distributions}).

For a subset $U\subseteq\mathbb{R}^{d}$ and $\mu\in\mathcal{M}(\mathbb{R}^{d})$
define the $U$-\emph{diffusion} of $\mu$ by
\[
\left\langle \mu\right\rangle _{U}=\int_{U}\delta_{T_{x}\mu}\, d\mu(x)
\]
i.e. $\left\langle \mu\right\rangle _{U}$ is the distribution of
the random measure $\nu$ obtained by choosing $x\in U$ according
to $\mu$ and setting $\nu=T_{x}\mu$. As usual we may apply normalization,
and usually do: $\left\langle \mu\right\rangle _{U}^{*}=\int_{U}\delta_{T_{x}^{*}\mu}\, d\mu(x)$
and $\left\langle \mu\right\rangle _{U}^{\sqr}=\int_{U}\delta_{T_{x}^{\sqr}\mu}\, d\mu(x)$.
\begin{defn}
\label{def:quasi-palm}A distribution $P$ on $\mathcal{M}^{*}$ is
\emph{quasi-Palm }if for every bounded, open neighborhood $U$ of
the origin,
\[
P\sim\int\left\langle \mu\right\rangle _{U}^{*}dP(\mu)
\]

\end{defn}
Thus $P$ is quasi-Palm when the following condition holds: $P(\mathcal{A})=0$
if and only if the random measure $T_{x}^{*}\mu$, obtained by selecting
$\mu\sim P$ followed by $x\sim\mu$, is almost surely not in $\mathcal{A}$.
As the name indicates, this definition generalizes Palm distributions,
which are distributions $P$ such that $\int\mu(B_{1})\, dP(\mu)<\infty$
and $P=\int\left\langle \mu_{U}\right\rangle \, dP(\mu)$ for all
symmetric neighborhoods $U$ of the origin (we have replaced equivalence
with equality, and have not normalized the diffusion). See \cite[Chapter 11]{Kallenberg83}.
\begin{defn}
\label{def:fractal-distributions}A \emph{fractal distribution}%
\footnote{The terminology derives from Furstenberg's notion of an ergodic fractal
measure \cite{Furstenberg08}, which essentially is a generic measure
for a CP-distribution (defined below). Later we shall see that a typical
measure for a FD is, up to a translation, an ergodic fractal measure
in this sense.%
} (FD) is a probability distribution on $\mathcal{M}^{*}$ which is
$S^{*}$-invariant and quasi-Palm. An \emph{ergodic fractal distribution}
(EFD) is a FD which is ergodic with respect to the $S^{*}$-action.
\end{defn}
When $P$ is a FD we shall often refer to its restricted version $P^{\sqr}$
as a FD as well, even though strictly speaking it is not. This is
justified by the 1-1 correspondence between restricted and extended
versions of $S^{*}$-invariant distributions. See section \ref{sub:remarks-on-normalization}
for some remarks on the definition.

This notion is closely related to the $\alpha$-scale-invariant distributions
of Z\"{a}hle \cite{Zahle88}. These are Palm distributions which
are invariant under a family of scaling operations depending on a
parameter $\alpha$. FDs are strictly more general objects and apply
in many cases where Z\"{a}hle distributions are inappropriate. The
relationship between them can be described precisely: up to normalization,
ergodic Z\"{a}hle distributions are those EFDs which are supported
on measures which have positive and finite second order densities
a.e.. See Section \ref{sub:Palm-and-Zahle-distributions}.

As we shall see later there is no shortage of FDs, but constructing
a non-trivial one directly requires a little work. At this point we
give two trivial examples. One may, first, take $P=\delta_{\lambda^{*}}$,
i.e. the distribution which is concentrated on (normalized) Lebesgue
measure $\lambda^{*}$. Clearly $\lambda^{*}$ is a fixed point for
both translation $T_{x}^{*}$ and for the scaling operators $S_{t}^{*}$,
so $P$ is an EFD. Second, one can consider the distribution $P=\delta_{\delta_{0}}$
. Since $\delta_{0}$ is $S^{*}$-invariant, and $\left\langle \delta_{0}\right\rangle _{U}^{*}=\delta_{\delta_{0}}$
for every neighborhood $U$ of the origin, $P$ is an EFD.

As an example of an $S^{*}$-invariant distribution which is not a
FD, consider the measure $\eta=(\lambda|_{[-\infty,0)})^{*}$ and
let $P=\delta_{\eta}$. Since $\eta$ is an $S^{*}$-fixed point $P$
is $S^{*}$-invariant, but for any neighborhood $U$ of the origin,
$\left\langle \eta\right\rangle _{U}$ is supported on measures which
a.s. give positive mass to $(0,\infty)$, whereas $\eta((0,\infty))=0$,
so $\left\langle \eta\right\rangle _{U}^{*}\not\sim\delta_{\eta}$.
Therefore $\delta_{\eta}$ is not a FD. 

The ergodic components of a FD with respect to $S^{*}$ are of course
$S^{*}$-invariant, but not a-priori quasi-Palm. This is the content
of:
\begin{thm}
\label{thm:EFD-ergodic-components}The ergodic components of a FD
are EFDs.
\end{thm}
We discuss this and other decompositions in Section \ref{sub:Ergodic-and-spatial-decomposition}.

\subsection{\label{sub:Uniformly-scaling-measures}Uniformly scaling measures}

Next, we examine individual measures which display dynamical regularity
upon ``zooming in'' to typical points. This idea has a long history;
one may view the density theorems of Lebesgue and Besicovitch as early
manifestations of it, and similarly the work of D. Preiss on tangent
measures. More recently the dynamical perspective has been taken up
by many authors \cite{PatzschkeZahle90,Bandt1992,Graf95,BedfordFisher96,Bandt1997,BedfordFisher97,KriegMorters98,Morters98,MortersPreiss98,BedfordFisherUrbanski02,Patzschke04,Fisher04,Gavish09,HochmanShmerkin09}.

Given a measure $\mu\in\mathcal{M}$ and $x\in\supp\mu$ one may translate
$x$ to the origin, forming $T_{x}\mu$, and consider the orbit of
this measure under $S^{*}$, which is called the \emph{scenery} of
$\mu$ at $x$. We are interested in measures for which for typical
$x$ the scenery equidistributes in the space of measures for some
distribution, i.e. the uniform measure on the orbit up to time $T$
converges, as $T\rightarrow\infty$, to some distribution. Unfortunately,
because $\mathcal{M}$ is not locally compact, the space $\mathcal{P}(\mathcal{M})$
does not carry a good topology.%
\footnote{Some authors have used vague convergence, but this causes various
complications, such as degenerate limit points and ``inseparability''
of the topology (though it is not really a topology). We prefer to
avoid these.%
} We therefore work in $\mathcal{P}(\mathcal{M}^{\sqr})$, which is
compact and metrizable in the weak topology. 
\begin{defn}
\label{def:sceneries-and-scenery-distributions}Let $\mu\in\mathcal{M}$
and $x\in\supp\mu$. The parametrized family $(\mu_{x,t}^{\sqr})_{t>0}$
given by 
\[
\mu_{x,t}=S_{t}(T_{x}\mu)
\]
is called the \emph{scenery} of $\mu$ at $x$. The \emph{scenery
distribution up to time }$T$ is obtained by placing length measure
on initial $T$-segments of the scenery: 
\[
\left\langle \mu\right\rangle _{x,T}=\frac{1}{T}\int_{0}^{T}\delta_{\mu_{x,t}}\, dt
\]
We almost always apply the $\sqr$ operation: $\mu_{x,t}^{\sqr}=S_{t}^{\sqr}T_{x}\mu$
and $\left\langle \mu\right\rangle _{x,T}^{\sqr}=\frac{1}{T}\int_{0}^{T}\delta_{\mu_{x,t}^{\sqr}}\, dt$.
These are defined for all $t$ and $T$ as along as $0\in\supp\mu$.
\end{defn}
Similar definitions have appeared previously in the literature. Our
method of normalization appears in Bandt \cite{Bandt1992,Bandt1997},
Graf \cite{Graf95}, and Gavish \cite{Gavish09}. Other authors \cite{Zahle88,MortersPreiss98,Patzschke04,PatzschkeZahle90}
have studied the normalization in which $\mu_{x,t}^{\alpha}=e^{\alpha t}S_{t}(T_{x}\mu)$,
which is appropriate when $\mu$ has second-order density $\alpha$
at $x$, but does not apply to many common cases, e.g. when second-order
densities do not exist. See section \ref{sub:Palm-and-Zahle-distributions}.

There is in general no reason the scenery distributions should converge
as $T\rightarrow\infty$. In the case where they do we introduce the
following terminology (Gavish \cite{Gavish09}):%
\footnote{Gavish uses the term ``measures with uniformly scaling sceneries'',
but we prefer the shorter name.%
}
\begin{defn}
\label{def:limit-distributions-and-scaling-measures}Let $\mu\in\mathcal{M}$
and $x\in\supp\mu$. 
\begin{enumerate}
\item If $\left\langle \mu\right\rangle _{x,T}^{\sqr}\rightarrow P$ as
$T\rightarrow\infty$ we say that $\mu$ \emph{generates} $P$ at
$x$.
\item If $\mu$ generates a distribution $P_{x}$ at $\mu$-a.e. $x$, we
say that $\mu$ is a \emph{scaling measure }(SM).
\item If $\mu$ generates the same distribution $P$ at a.e. point then
$\mu$ is a \emph{uniformly scaling measure }(USM) and $\mu$ \emph{generates}
$P$.
\end{enumerate}
\end{defn}
There is a close relation between FDs and SMs, similar to the relation
between an invariant measure on a dynamical system and a generic point
for the measure. In one direction we have:
\begin{thm}
\label{thm:FD-to-USM}If $P$ is a FD then $P$-a.e. measure is a
USM generating the ergodic component of $P$ to which $\mu$ belongs. 
\end{thm}
For EFDs this is an easy consequence of the definitions and the ergodic
theorem, which in turn implies the non-ergodic case using Theorem
\ref{thm:EFD-ergodic-components}. For a detailed proof see Section
\ref{sub:CP-to-USM}. The converse is also true but less trivial:
\begin{thm}
\label{thm:limit-distributions-are-FDs}If $\mu\in\mathcal{M}$ then
for $\mu$-a.e. $x$, every accumulation point of $\left\langle \mu\right\rangle _{x,T}^{\sqr}$
is a restricted FD. In particular, if $\mu$ generates $P$ then $P$
is a restricted FD.
\end{thm}
See section \ref{sec:Equivalence-of-the}. While it is fairly obvious%
\footnote{Actually this is somewhat tedious to prove directly, due to the discontinuity
of the action of $S^{\sqr}$.%
} that an accumulation point of the scenery distributions is $S^{\sqr}$-invariant,
it is remarkable that the distribution should have the additional
spatial invariance of a FD. This fact is analogous to a similar phenomenon
discovered by P. M\"{o}rters and D. Preiss \cite{MortersPreiss98},
namely that the limiting distributions of the sceneries $(\mu_{x,t}^{\alpha})_{t\geq0}$
are scale-invariant Palm distributions.

Note also that in general the distributions generated by scenery distributions
of a measure need not be $S^{*}$-ergodic. See Section \ref{sub:USM-generating-non-ergodic-FD}.

One important property of the FDs and the limiting distributions of
sceneries is that they behave nicely under linear maps, and, since
they are infinitesimal notions, under local diffeomorphisms. 
\begin{prop}
\label{prop:linear-images-of-FDs}If $L\in GL(\mathbb{R}^{d})$ is
a linear map then the induced map $L^{*}:\mathcal{M}^{*}\rightarrow\mathcal{M}^{*}$
maps EFDs to EFDs. 
\end{prop}

\begin{prop}
\label{pro:smooth-images-of-usms}Let $\mu\in\mathcal{M}$ be supported
on an open set $U\subseteq\mathbb{R}^{d}$ and let $f:U\rightarrow V\subseteq\mathbb{R}^{d}$
be a diffeomorphism. Then for $x\in U$, if $\left\langle \mu\right\rangle _{x,T_{i}}^{\sqr}\rightarrow P$
for some sequence of times $T_{i}\rightarrow\infty$ and $P$ is an
FD, then $\left\langle f\mu\right\rangle _{fx,T_{i}}^{\sqr}\rightarrow(Df(x))^{\sqr}P$.
In particular if $\mu$ generates $P_{x}$ at $x$ then $f\mu$ generates
$Df(x)^{*}P_{x}$ at $f(x)$.
\end{prop}
See Section \ref{sub:Linear-images-of-FDs}. For a closely related
result see \cite{AspbergEkstromPerssonSchmeling2013}.

\subsection{\label{sub:CP-distributions}CP-distributions}

Another closely related family of ``self-similar'' distributions
are Furstenberg's CP-processes \cite{Furstenberg08}, or, in our formalism,
CP-distributions. Fix an integer $b\geq2$ and let $\mathcal{D}_{b}=\mathcal{D}_{b}(\mathbb{R})$
denote the partition of $\mathbb{R}$ into intervals right-open intervals
of length $2/b$, with the origin in the center of one of them. Formally,
$\mathcal{D}_{b}=\{[\frac{k}{b},\frac{k+2}{b})\,:\, k=1\bmod2\}$.
This partition is chosen so that it divides $[-1,1)$ into $b$ congruent
intervals.%
\footnote{It would be more natural to work with the partition of $\mathbb{R}$
into into $b$-adic intervals, i.e. intervals of the form $[\frac{k}{b},\frac{k+1}{b})$,
which divides $[0,1]$ into $b$ subintervals. The definitions of
CP-distributions and processes in \cite{Furstenberg08,HochmanShmerkin09}
follow this path. This would make sense if we used $[0,1]$ as our
``basic'' interval instead of $B_{1}=[-1,1]$, but much of our notation
is adapted to $B_{1}$, e.g. the maps $T_{B}$ and normalization operators,
and because of this it is more efficient to adopt the present non-standard
definition.%
} For $d>1$ define $\mathcal{D}_{b}=\mathcal{D}_{b}(\mathbb{R}^{d})$
to be the partition of $\mathbb{R}^{d}$ into cells of the form $I_{1}\times\ldots\times I_{d}$,
$I_{i}\in\mathcal{D}_{b}(\mathbb{R})$, so $B_{1}$ is partitioned
into $b^{d}$ congruent cells. We denote by $\mathcal{D}_{b}(x)$
the unique element of $\mathcal{D}_{b}$ containing $x$. Note that
$\mathcal{D}_{b^{n+1}}^{d}$ refines $\mathcal{D}_{b^{n}}^{d}$. 
\begin{defn}
\label{def:magnification-operator}The base-$b$ magnification operator
$M_{b}:\mathcal{M}\times B_{1}\rightarrow\mathcal{M}\times B_{1}$
is
\[
M_{b}(\mu,x)=(T_{\mathcal{D}_{b}(x)}\mu,T_{\mathcal{D}_{b}(x)}x)
\]
We denote by $M_{b}^{*}$, $M_{b}^{\sqr}$ the associated operators
in which, after $M_{b}$ is applied, $*,\sqr$ are applied, respectively,
to the measure component of the output. These are defined on pairs
$(\mu,x)$ with $\mu(\mathcal{D}_{b}(x))>0$.
\end{defn}
Both $M_{b}$ and its domain are measurable, though $M_{b}$ may be
discontinuous at $(\mu,x)$ when $\mu(\partial\mathcal{D}_{b}(x))>0$.
Starting from $(\mu,x)$, the orbit $(M_{b}^{\sqr})^{n}(\mu,x)=M_{b^{n}}^{\sqr}(\mu,x)$,
$n=0,1,2,\ldots$, assuming it is defined, may be viewed as the scenery
obtained by zooming in to $x$ along $\mathcal{D}_{b}$-cells. We
call this the $b$-\emph{scenery} at $x$. Note that this is a discrete
time sequence, and the point $x$ is generally not in the ``center
of the frame'' as it is for sceneries.
\begin{defn}
\label{def:adapted-distributions}A distribution $Q$ on $\mathcal{M}^{\sqr}\times B_{1}$
is \emph{adapted} if, conditioned on the first component being $\nu$,
the second component is distributed according to $\nu$. Equivalently,
for every $f\in C(\mathcal{M}^{\sqr}\times B_{1})$,
\[
\int f(\nu,x)\, dQ(\nu,x)=\int\left(\int f(\nu,x)\, d\nu(x)\right)dQ(\nu)
\]

\end{defn}
Here we have written $dQ(\nu)$ instead of $dQ(\nu,y)$ since $y$
is not used. We shall often make this implicit identification between
$Q$ and its marginal on to $\mathcal{M}^{\sqr}$.
\begin{defn}
\label{def:CP-distributions}A restricted base-$b$ \emph{CP-distribution}%
\footnote{The terminology follows \cite{Furstenberg08}. ``CP'' stands for
conditional probability, hinting at the role of adaptedness in the
definition.%
} is an $M_{b}^{\sqr}$-invariant and adapted probability distribution
on $\mathcal{M}^{\sqr}\times B_{1}$.
\end{defn}
Given $b$ it is easy to see that the measure component of a base-$b$
CP-distribution determines the original distribution. If $b$ is not
given this may not hold: distributions supported on Lebesgue measure
or atomic measures are CP-distributions for several bases. It is not
clear if there are any other examples. One may suspect that the answer
is negative based on the analogy with of Furstenberg's $\times2,\times3$
problem.

To every restricted CP-distribution $P$ there is an associated extended
version on $\mathcal{M}^{*}\times B_{1}$ which is invariant under
$M_{b}^{*}$ and projects to $Q$ under $(\mu,x)\mapsto(\mu^{\sqr},x)$.
This distribution is not unique but the construction is canonical;
see Section \ref{sub:Restricted-and-extended-versions}.

The relation of CP-distributions and USMs was previously examined
by Gavish \cite{Gavish09}, who showed that CP-distributions typically
give rise to USMs. We strengthen this to show that they give rise
to FDs. 
\begin{defn}
\label{def:centering}For $\mu\in\mathcal{M}$ and $x\in\supp\mu$,
let 
\[
\cnt_{0}(\mu,x)=T_{x}^{*}\mu
\]
Thus $\cnt_{0}$ is a partially defined map $\mathcal{M}\times\mathbb{R}^{d}\rightarrow\mathcal{M}^{*}$.
The \emph{discrete centering }of an extended base-$b$ CP-distribution
$Q$ is the distribution $\cnt_{0}Q$. The \emph{continuous centering
}$\cnt Q$ is the distribution
\[
\cnt Q=\frac{1}{\log b}\int_{0}^{\log b}S_{t}^{*}\cnt_{0}Q\, dt
\]
\end{defn}
\begin{thm}
\label{thm:CP-to-EFD}The continuous centering $P=\cnt Q$ of a CP-distribution
$Q$ is a FD, and the measure-preserving system $(P,S^{*})$ is a
factor of the suspension of $(Q,M_{b})$ by the function with constant
height $\log b$. 
\end{thm}
For the proof see Section \ref{sub:CP-to-USM}. There is also a converse:
\begin{thm}
\label{thm:FD-to-CP}If $P$ is a FD and $b\geq2$ then $P$ is the
continuous centering of some base-$b$ CP-distribution, and each ergodic
component of $P$ with respect to the map $S_{\log b}^{*}$ is the
discrete centering of some base-$b$ CP-distribution.
\end{thm}
Suspensions are discussed in Section \ref{sub:Suspenssions}. For
the proof see Section \ref{sec:Equivalence-of-the}. 

The CP-distribution in the theorem above is highly non-unique. We
record one useful manifestation of this:
\begin{prop}
\label{pro:CP-distribution-in-arbitrary-coords}The CP-distributions
in Theorem \ref{thm:FD-to-CP} may be chosen so that its second component
is distributed according to Lebesgue measure on $B_{1}$.
\end{prop}
The centering operation establishes (by definition) a measure-preserving
map under which typical measure $\nu$ from the CP-distribution are
mapped to a translated, scaled and normalized version of $\nu$. These
operations preserve most geometric properties of $\nu$. This is significant
because CP-distributions, while being powerful analytic tools, are
tied inflexibly to a particular coordinate system and base. For example,
the results of Furstenberg on dimension conservation \cite{Furstenberg08}
are for the projection of a CP-distribution onto one of the coordinate
planes. One would like to be able to change the coordinates and base.
This is made possible by the two theorems above, which allow one to
pass from a CP-distribution to a FD, the latter object being defined
in a coordinate-free way.%
\footnote{In fact there is a mild dependence on coordinates in FDs since $B_{1}$
depends on the coordinates and is used to define normalization, but
this is insignificant. See Section \ref{sub:remarks-on-normalization}.%
} Thus one can pass from a CP-distribution to an FD and back to a CP-distribution
in a new coordinate system, and the net effect is only that the measures
have been translated and normalized. In particular, a geometric property
of measures which is unchanged by translations and scaling of measures
will hold a.s. for the old CP-distribution if and only if it holds
almost surely for the new one. We summarize this as follows:
\begin{cor}
\label{cor:CP-change-of-coords}For any extended FD or CP-distribution
$Q$ and for any choice of base $b$ and coordinate system, there
is an extended CP-distribution $Q'$ in this base and coordinate system
such that, for any Borel set $\mathcal{A}\subseteq\mathcal{M}$ which
is invariant under translations, scaling and normalization, we have
$Q(\mathcal{A})=Q'(\mathcal{A})$.
\end{cor}
See Section \ref{sub:Changing-coordinates}.

Finally, we note that in certain situations it is useful to generalize
the notion of CP-distributions to allow for more complicated partitions
than $\mathcal{D}_{b}$. Some examples are provided in Section \ref{sec:Examples}.
See also \cite{HochmanShmerkin09}.

\subsection{\label{sub:Dimension}Dimension }

We now turn to the geometric properties of USMs and EFDs. We recall
below some background on dimension. Falconer's books \cite{Falconer90,Falconer97}
are good introductions to the topic.

Denote the Hausdorff dimension of a set $A$ by $\dim A$. For a measure
$\mu\in\mathcal{M}$ and $x\in\supp\mu$ the \emph{upper and lower
local dimension }of $\mu$ and $x$ are 
\begin{eqnarray*}
\overline{D}_{\mu}(x) & = & \limsup_{r\searrow0}\frac{\log\mu(B_{r}(x))}{\log r}\\
\underline{D}_{\mu}(x) & = & \liminf_{r\searrow0}\frac{\log\mu(B_{r}(x))}{\log r}
\end{eqnarray*}
If these quantities are equal then their value is the \emph{exact
(or pointwise) dimension }of $\mu$ at $x$, denoted $D_{\mu}(x)$.
If $D_{\mu}(x)=\alpha$ at $\mu$-a.e. $x$, we say that $\mu$ is
\emph{exact dimensional }and has \emph{exact dimension $\alpha$},
and we write $\dim\mu=\alpha$. The exact dimension need not exist,
but one always can define the lower (Hausdorff) dimension of $\mu$,
\[
\ldim\mu=\essinf_{x\sim\mu}\underline{D}_{\mu}(x)
\]
Equivalently \cite[Theorem 1.2]{FanLauRao02},
\begin{equation}
\ldim\mu=\inf\{\dim A\,:\,\mu(A)>0\}.\label{eq:ldim-by-dim-of-supporting-sets}
\end{equation}

\begin{lem}
\label{lem:dimension-of-FDs}If $P$ is an EFD then $P$-a.e. $\mu$
is exact dimensional and the dimension is $P$-a.s. constant. Writing
$\dim P$ for the a.s. value of the dimension, for every $0<r<1$
we have
\[
\dim P=\int\frac{\log\mu(B_{r}(0))}{\log r}\, dP(\mu)
\]

\end{lem}
The fact that $P$-a.s. the dimension is constant follows from ergodicity
of $P$ by noticing that the map $\mu\mapsto\dim\mu$ is $S^{*}$-invariant.%
\footnote{One may check that the map is measurable. We shall generally omit
these routine verifications.%
} Existence of exact dimension along with the formula for it is proved
in Section \ref{sub:Dimension-of-FDs}.

In general, we define the dimension of an FD $P$ by
\begin{eqnarray*}
\dim P & = & \int\dim\nu\, dP(\nu)
\end{eqnarray*}
If $P$ is a FD with ergodic decomposition $P=\int P_{\nu}\, dP(\nu)$,
this is the same as $\dim P=\int\dim P_{\nu}\, dP(\nu)$. Note that
when $P$ is not ergodic it is not true that $P$-a.e. $\nu$ has
dimension $\dim P$; rather, $P$-a.e. $\nu$ has dimension $\dim P_{\nu}$
where $P_{\nu}$ is the ergodic component of $\nu$. 

For $\mu\in\mathcal{M}$ let $\mathcal{V}_{x}$ denote the accumulation
points in $\mathcal{P}(\mathcal{M}^{\sqr})$ of $\left\langle \mu\right\rangle _{x,T}^{\sqr}$
as $T\rightarrow\infty$. We have the following bounds for the local
dimension of $\mu$:
\begin{prop}
\label{pro:dimension-via-accumulation-pts}If $\mu\in\mathcal{M}$
then for $\mu$-a.e. $x$ 
\[
\underline{D}_{\mu}(x)\geq\inf_{P\in\mathcal{V}_{x}}\dim P
\]
and
\[
\overline{D}_{\mu}(x)\leq\sup_{P\in\mathcal{V}_{x}}\dim P
\]
In particular, if $\mu$ is a USM generating a FD $P$ then $\mu$
is exact dimensional and $\dim\mu=\dim P$. 
\end{prop}
See Section \ref{sub:Dimension-of-FDs}.

\subsection{\label{sub:Geometric-properties}Projections of measures}

One of the fundamental phenomena in fractal geometry is that if one
projects a set or measure through a typical linear map $\pi:\mathbb{R}^{d}\rightarrow\mathbb{R}^{k}$
(or pushes it through a typical $C^{1}$ map) then the image measure
has dimension which is ``as large as it can be''. This result has
many variants but is often known generally as Marstrand's theorem.
More precisely, let $\Pi_{d,k}$ denote%
\footnote{We shall mostly be interested in the elements of $\Pi_{d,k}$ only
up to change of coordinates in the range. Some authors identify $\Pi_{d,k}$
modulo this relation with the space of orthogonal projections from
$\mathbb{R}^{d}$ to $k$-dimensional subspaces.%
} the space of linear maps $\mathbb{R}^{d}\rightarrow\mathbb{R}^{k}$.
Note that $\Pi_{d,k}$ is a smooth manifold and carries a natural
measure class (for example, $\Pi_{d,k}$ can be identified with the
space of $k\times d$ matrices, with the measure induced from Lebesgue
measure on $\mathbb{R}^{dk}$). Marstrand's theorem for measures is:
\begin{thm}
[Hunt-Kaloshin, \cite{HuntKaloshin97}] \label{thm:Hunt-Kaloshin}
If $\mu$ is an exact-dimensional finite measure on $\mathbb{R}^{d}$,
then for a.e. $\pi\in\Pi_{d,k}$, the image measure $\pi\mu$ is exact
dimensional and
\[
\dim\pi\mu=\min\{k,\dim\mu\}
\]

\end{thm}
This general statement has two shortcomings. First, it is an almost-everywhere
statement, and gives no information about $\pi\mu$ for particular
$\pi\in\Pi_{d,k}$. Second, in general the map $E_{\mu}:\Pi_{d,k}\rightarrow\mathbb{R}$
given by 
\[
E_{\mu}(\pi)=\ldim\pi\mu
\]
is Borel, but can be highly discontinuous. 

In the context of CP-distributions, however, Hochman and Shmerkin
\cite{HochmanShmerkin09} recently proved that some continuity exists.
To state this, we note first that for $\mu\in\mathcal{M}$ the projection
$\pi\mu$ may not be Radon, so we define%
\footnote{Although it would be more precise to write $\dim(\pi,\mu)$ and not
$\dim\pi\mu$, no confusion should arise.%
}
\[
\ldim\pi\mu=\lim_{R\rightarrow\infty}(\ldim\pi(\mu|_{B_{R}}))
\]
(it is a general fact that $\nu\ll\nu'$ implies $\ldim\nu\leq\ldim\nu'$,
therefore $\ldim\pi(\mu|_{B_{R}})$ is non-decreasing and the limit
exists). We make a similar definition for $\dim$, if it exists. Next,
if $P$ is a distribution on $\mathcal{M}$, write 
\begin{eqnarray*}
E_{P}(\pi) & = & \int E_{\mu}(\pi)\, dP(\mu)\\
 & = & \int\ldim\pi\mu\, dP(\mu)
\end{eqnarray*}
If $P$ is $S^{*}$-ergodic then the map $\mu\mapsto\ldim\pi\mu$
is $S^{*}$-invariant and hence $P$-a.s. constant, so the integral
in the definitions of $E_{P}(\pi)$ trivializes and $E_{P}(\pi)$
is just the almost-sure value of $\ldim\pi\mu$. 
\begin{thm}
[Hochman-Shmerkin \cite{HochmanShmerkin09}, Theorem 1.10] \label{thm:semi-continuity-for-CP-processes}
Let $P$ be an ergodic CP-distribution.
\begin{enumerate}
\item $E_{P}(\cdot)$ is lower semi-continuous and equal almost everywhere
to $\min\{k,\dim P\}$.
\item For $P$-a.e. $\mu$ and every regular $\varphi\in C^{1}(\mathbb{R}^{d},\mathbb{R}^{k})$,
\[
\ldim\varphi\mu\geq\essinf_{x\sim\mu}E_{P}(D\varphi(x))
\]
In particular, $\ldim\pi\mu\geq E_{P}(\pi)$ for $\pi\in\Pi_{d,k}$.
\end{enumerate}
\end{thm}
Using  Theorem \ref{thm:FD-to-CP}, this transfers immediately to
EFDs. We can also improve it in several ways. First, Furstenberg \cite{Furstenberg08}
showed that if $Q$ is a (restricted) CP-distribution and $\pi$ is
the projection to a coordinate plane, i.e. $\pi(x)=(x_{i_{1}},\ldots,x_{i_{k}})$
for some indices $1\leq i_{1}<\ldots<i_{k}\leq d$, then $\pi\mu$
is exact dimensional for $Q$-a.e. $\mu$. Given some other $\pi'\in\Pi_{d,k}$
and an FD- or CP-distribution $Q$, we can pass to a CP-distribution
$Q'$ such that $\pi'$ is the coordinate projection in those coordinates.
The result of Furstenberg above implies that $Q'$-a.e. measure projects
to an exact dimensional measure under $\pi'$. Since this property
is invariant under translation, scaling and normalization (and may
be seen to be Borel), by Corollary \ref{cor:CP-change-of-coords}
we have:
\begin{thm}
\label{thm:exact-dimension-of-projections}If $P$ is a CP-distribution
or a FD then for every $\pi\in\Pi_{d,k}$, for $P$-a.e. $\mu$ the
image $\pi\mu$ is exact dimensional.
\end{thm}
The last two theorems are still a.e results, the uncertainty being
about the measure rather than the projection. Instead, one would like
to obtain information about individual measures. For example, if $Q$
is a CP-distribution, is it true for $Q$-a.e. $\mu$ that $\pi\mu$
is exact dimensional, and that $E_{\mu}(\pi)=E_{Q}(\pi)$, for \emph{all}
$\pi\in\Pi_{d,k}$? Our results only give this with the quantifiers
in the weaker, reverse order. Since a typical measure for a FD is
a USM, one approach to these questions is to explore the validity
of the results above for USMs and, more generally, SMs.

In this spirit, one can get lower bounds for projections of scaling
measures. For a $\mu\in\mathcal{M}$ let $\mathcal{V}_{x}$ again
denote the accumulation points of $\left\langle \mu\right\rangle _{x,T}^{\sqr}$
as $T\rightarrow\infty$, and let 
\begin{equation}
E_{\mu,x}(\pi)=\inf_{P\in\mathcal{V}_{x}}E_{P}(\pi)\label{eq:Ex}
\end{equation}

\begin{thm}
\label{thm:projection-of-SMs}Let $\mu\in\mathcal{M}$. Then for regular
$f\in C^{1}(\mathbb{R}^{d},\mathbb{R}^{k})$, 
\[
\ldim f\mu\geq\essinf_{x\sim\mu}E_{\mu,x}(Df(x))
\]
In particular, if $\mu$ is a USM generating an EFD $P$, then for
all $\pi\in\Pi_{d,k}$, 
\[
\ldim\pi\mu\geq E_{P}(\pi)
\]

\end{thm}
See Section \ref{sub:Dimension-of-projected-measures}. However, one
cannot hope for equality or exact dimension:
\begin{prop}
\label{pro:sm-with-bad-projections}Let $\mu$ be a USM generating
an EFD $P$ and $\pi\in\Pi_{d,k}$. Then it is possible that $\pi\mu$
is not exact dimensional and that $\ldim\pi\mu>E_{P}(\pi)$.
\end{prop}
We give such examples in Section \ref{sub:examples}. As a consequence,
since USMs are exact dimensional we have:
\begin{cor}
\label{pro:projection-of-USM-may-not-be-USM}The linear image of a
USM need not be a USM.
\end{cor}
We have not been able to settle the following:
\begin{problem}
\label{problem:good-projections}If $P$ is an EFD on $\mathbb{R}^{d}$,
is it true that for $P$-a.e $\mu$, for every $\pi\in\Pi_{d,k}$
the image $\pi\mu$ is exact dimensional with $\dim\pi\mu=E_{P}(\pi)$?
\end{problem}

\subsection{\label{sub:Conditional-measures}Conditional measures on subspaces}

Another classical question, which is in some sense dual to the problem
of understanding projections, is to understand the conditional measures
of $\mu$ on the fibers of a projection. More precisely, given a Borel
map $f:\mathbb{R}^{d}\rightarrow\mathbb{R}^{k}$, define for $x\in\mathbb{R}^{d}$
the $f$-fiber through $x$ by
\begin{eqnarray*}
[x]_{f} & =\{y\in\mathbb{R}^{d}\,:\,\pi y=\pi x\}= & f^{-1}(f(x))
\end{eqnarray*}
and the corresponding partition into fibers, 
\begin{eqnarray*}
\mathcal{F}_{f} & = & \{f^{-1}(y)\,:\, y\in\mathbb{R}^{k}\}\\
 & = & \{[x]_{f}\,:\, x\in\mathbb{R}^{d}\}
\end{eqnarray*}
Given a measure $\mu\in\mathcal{P}(\mathbb{R}^{d}$), there is a family
of measures $\mu_{[x]_{\pi}}\in\mathcal{P}(\mathbb{R}^{d})$ for $x\in\mathbb{R}^{d}$,
called the \emph{conditional measures of $\mu$ on the fibers $[x]_{\pi}$},
with the following properties:
\begin{enumerate}
\item $x\mapsto\mu_{[x]_{\pi}}$ is measurable with respect to the $\sigma$-algebra
$\{\pi^{-1}(U)\,:\, U\subseteq\mathbb{R}^{k}\mbox{ Borel}\}$ (this
justifies the notation $\mu_{[x]_{\pi}}$, since it means that the
measure depends on $[x]_{\pi}$, not $x$).
\item $\mu_{[x]_{\pi}}$ is supported on $[x]_{\pi}$ for $\mu$-a.e. $x$.
\item For every Borel set $E\subseteq\mathbb{R}^{d}$,
\[
\mu(E)=\int\mu_{[x]_{\pi}}(E)\, d\nu(x)
\]

\item The family of measures is unique in the sense that any other map $x\mapsto\mu'_{x}$
satisfying the three conditions above satisfies $\mu'_{x}=\mu_{[x]_{\pi}}$
at $\nu$-a.e. $x$.
\end{enumerate}
For a discussion of measure-valued functions and integration see Section
\ref{sub:Measure-valued-integration}. More generally, such a family
of conditional measures exists for any measurable map $f:X\rightarrow Y$
between standard Borel spaces and Borel probability measure on $X$.
See \cite[Theorem 6.3]{Kallenberg83}. 

If $\mu\in\mathcal{M}$ then the conditional measures are locally
defined up to a constant in the neighborhood of a.e. point. If $B,B'$
are two balls, then for $\mu$-a.e. $x\in B\cap B'$ then the restriction
to $B\cap B'$ of the conditional measure measures $(\mu_{B})_{[x]_{\pi}}$
and $(\mu_{B'})_{[x]_{\pi}}$ are proportional to each other (see
Lemma \ref{lem:restriction-of-conditionals}). Thus local properties
of conditional measures, such as exact dimensionality, local dimension
etc., are well defined for fiber measures of $\mu\in\mathcal{M}$.
Furthermore one can define the conditional measure $\mu_{[x]_{\pi}}$
up to a constant when $\mu\in\mathcal{M}$ is infinite, by gluing
together the conditional measures of $\mu_{B}$ for a countable number
of balls $B$ covering $\mathbb{R}^{d}$, and using the consistency
property above. When we speak of fiber measures often we do so in
this sense.

There are many classical results about fiber measures. For example,
\begin{thm}
[Matilla \cite{Mattila95}]\label{thm:Matilla-conditional-measures}
Let $\mu\in\mathcal{M}$ be a probability measure that is absolutely
continuous with respect to the $\alpha$-dimensional Hausdorff measure
for some $\alpha$, in particular $\mu$ is exact dimensional and
$\dim\mu=\alpha$. Then for a.e. $\pi\in\Pi_{d,k}$ and $\mu$-a.e.
$x$, the conditional measure $\mu_{[x]_{\pi}}$ is exact dimensional,
with dimension
\begin{equation}
\dim\mu_{[x]_{\pi}}=\max\{0,\dim\mu-k\}\label{eq:typical-fiber-dimension}
\end{equation}
and $\mu_{[x]_{\pi}}$ is absolutely continuous with respect to the
$\dim\mu_{[x]_{\pi}}$-dimensional Hausdorff measure. In particular
for a.e. $\pi\in\Pi_{d,k}$, 
\begin{equation}
\dim\pi\mu+\dim\mu_{[x]_{\pi}}=\dim\mu\qquad\mbox{for }\mu\mbox{-a.e. }x\label{eq:dimension-conservation}
\end{equation}

\end{thm}
The last equation is a ``dimension conservation'' phenomenon: for
an exact-dimensional measure $\mu$ and a.e. projection $\pi\in\Pi_{d,k}$,
the sum of dimensions of the image $\pi\mu$ and a typical fiber $\mu_{[x]_{\pi}}$
is precisely $\dim\mu$. However, in general this only holds for typical
maps $\pi$. It can fail drastically for particular measures and maps.
For example let $\pi(x,y)=x$ and let $\mu$ be the appropriate Hausdorff
measure on the graph $(t,W_{t})$ where $W_{t}$ is a typical Brownian
motion path. Then it is well known that $\dim\mu>1$. But $\dim\pi\mu\leq1$
and $\pi^{-1}(x)$ is a single point for all $x$, hence dimension
conservation fails.

We now turn to the behavior of FDs under conditioning, which is discussed
in more detail in Section \ref{sub:Exact-dimension-of-of-conditionals}.
Let us first make some remarks about the global picture. Although
the conditional measures of $\mu$ with respect to $\pi$ are defined
only for $\mu$-a.e. fiber, if $P$ is a $d$-dimensional FD and $\pi\in\Pi_{d,k}$
then the conditional measure of a $P$-typical $\mu$ on $[x]_{\pi}$
is defined $\mu$-a.e., so by the quasi-Palm property, $P$-a.e. it
is defined on $[0]_{\pi}$. The fiber $[0]_{\pi}$ can be identified
with $\mathbb{R}^{d-k}$, and one may verify without difficulty that
the map $\mu\mapsto(\mu_{\pi^{-1}(0)})^{*}$ intertwines $S^{*}$
and that the image $Q$ of $P$ by this map is a $k$-dimensional
FD. We thus have:
\begin{prop}
\label{pro:exact-dimension-of-conditionals}If $P$ is a $d$-dimensional
EFD and $\pi\in\Pi_{d,k}$, then the push-forward $Q$ of $P$ by
$\mu\mapsto(\mu_{[0]_{\pi}})^{*}$ is an EFD and the given map is
a factor map between $(P,S^{*})$ and $(Q,S^{*})$. In particular,
for $P$-a.e. $\mu$ the conditional measures $\mu_{[0]_{\pi}}$ is
well defined, as are $\mu_{[x]_{\pi}}$ for $\mu$-a.e. $x$, and
these measures are exact dimensional and their dimension is a.s. independent
of $\mu$ and $x$.
\end{prop}
See Section \ref{sub:Exact-dimension-of-of-conditionals}. 

In \cite{Furstenberg08}, Furstenberg established a version of Theorem
\ref{thm:Matilla-conditional-measures} in which the measure is a
typical measures for a CP-processes, but only for projections to coordinate
planes in the coordinate system in which the CP-distribution is defined.
Furstenberg also proved the following:
\begin{thm}
[Furstenberg \cite{Furstenberg08}]\label{thm:Furstenberg-dimension-conservation-1}
If $P$ is an ergodic CP-distribution and $\pi$ is a coordinate projection,
then for $P$-a.e. $\mu$,\textup{
\[
\dim\pi\mu+\dim\mu_{[x]_{\pi}}=\dim\mu\qquad\mbox{for }\mu\mbox{-a.e. }x
\]
}
\end{thm}
Since this property is invariant under $M_{b}^{\sqr}$, it holds for
the extended version of a CP-distribution also.

Now, given an FD $P$ and $\pi\in\Pi_{d,k}$, we can fix a coordinate
system in which $\pi$ is the coordinate projection $(x_{1},\ldots,x_{d})\mapsto(x_{1},\ldots,x_{k})$.
By Theorem \ref{thm:FD-to-CP}, there is a CP-distribution $Q$ in
this coordinate system (and an arbitrary base $b$) such that $P=\cnt Q$.
By the previous theorem, $Q$-a.e. measure satisfies dimension conservation.
Since this property is preserved by the centering operation, we have
\begin{thm}
\label{thm:fibers-for-EFDs-1} Let $P$ be an EFD and $\pi\in\Pi_{d,k}$.
Then for $P$-a.e. $\mu$ and $\mu$-a.e. $x$, the fiber measure
$\mu_{[x]_{\pi}}$ is exact dimensional, and 
\begin{equation}
\dim\pi\mu+\dim\mu_{[x]_{\pi}}=\dim\mu\qquad\mbox{for }\mu\mbox{-a.e. }x\label{eq:FD-dimension-conservation}
\end{equation}

\end{thm}
See Section \ref{sub:Exact-dimension-of-of-conditionals}. 

Turning to USMs, we have a dual version of Theorem \ref{thm:projection-of-SMs},
giving an upper bound on the dimension of fibers:
\begin{thm}
\label{thm:conditional-measures-of-USMs}If $\mu\in\mathcal{M}$ and
$E_{x}$ is defined as in equation \eqref{eq:Ex}, then for every
regular $f\in C^{1}(\mathbb{R}^{d},\mathbb{R}^{k})$,
\[
\ldim\mu_{[x]_{f}}\leq\esssup_{x\sim\mu}(\dim\mu-E_{x}(Df(x)))\qquad\mbox{ for }\mu\mbox{-a.e.}x
\]
In particular, if $\mu$ is a USM generating an EFD $P$ and $\pi\in\Pi_{d,k}$
then 
\[
\ldim\mu_{[x]_{f}}\leq\dim P-E_{P}(\pi))\qquad\mbox{ for }\mu\mbox{-a.e.}x
\]

\end{thm}
In contrast, dimension conservation can fail rather dramatically for
USMs:
\begin{prop}
\label{pro:dimension-conservation-counterexample}There exists a USM
$\mu$ on $\mathbb{R}^{2}$ generating an EFD $P$, such that $\dim(\mu)>1$
but the projection $(x,y)\mapsto x$ is an injection on $X$.
\end{prop}
See Section \ref{sub:non-conservation-example}. We are left the dual
version of Problem \ref{problem:good-projections}:
\begin{problem}
If $P$ is a FD on $\mathbb{R}^{d}$, is it true that for $P$-a.e
$\mu$, for every $\pi\in\Pi_{d,k}$ a.e. fiber $\mu_{\pi^{-1}(y)}$
is exact dimensional with $\dim\mu_{\pi^{-1}(y)}=\dim\mu-E_{P}(\pi)$?
\end{problem}
Finally, the same construction that proves Proposition \ref{pro:dimension-conservation-counterexample}
provides a counterexample to another question of Furstenberg. Following
\cite{Furstenberg08}, for a compact set $X\subseteq\mathbb{R}^{d}$
define a miniset to be a non-empty set of the form $B_{1}\cap T_{E}X$,
where $E$ is a ball, and define a microset to be a limit of minisets,
with respect to the Hausdorff metric on the space of closed subsets
of $B_{1}$. The set $X$ is called \emph{homogeneous} if every microset
is a subset of some miniset. For such sets Furstenberg constructed
an ergodic CP-process of the same dimension as $X$ whose measures
are supported on minisets of $X$, and used this to show that, for
linear projections, a dimension conservation phenomenon holds for
$X$. He also asked \cite{Furstenberg22001} whether this holds for
smooth maps. In Section \ref{sub:non-conservation-example} we give
a negative answer:
\begin{prop}
\label{pro:dim-conservation-of-smooth-maps-1}There exists a homogeneous
set $X\subseteq\mathbb{R}^{2}$ and a $C^{\infty}$ map $f:\mathbb{R}^{2}\rightarrow\mathbb{R}$
which is injective on $X$ (so each fiber is a singleton) and such
that $\dim fX<\dim X$.
\end{prop}

\subsection{\label{sub:FDs-with-additional-invariance}FDs with additional invariance}

For USMs and FDs which enjoy additional geometric invariance properties
one can draw stronger conclusions than the above. We begin with the
notion of a homogeneous measure, which is a modification of a similar
notion of Gavish \cite{Gavish09}.%
\footnote{Gavish's original definition is flawed in a number of ways. Most measures
of interest are not homogeneous in his sense; it is simple, for example,
to show that any measure on a cantor set in the line has atomic micromeasures
and hence if the original measure is non-atomic it cannot be homogeneous
in Gavish's sense. Our definition seems to capture better Gavish's
intention.%
}
\begin{defn}
\label{def:homogeneous-measures} A point $x$ is \emph{homogeneous
}for a measure $\mu\in\mathcal{M}$ if for every accumulation point
$\nu$ of $(\mu_{x,t}^{\sqr})_{t>0}$ there is a ball $B$ with $\mu\ll T_{B}\nu$.
A measure $\mu$ is homogeneous if $\mu$-a.e. point is homogeneous
for $\mu$.
\end{defn}
Examples of homogeneous measures include self-similar measures for
iterated function systems with strong separation whose contractions
are homotheties or, more generally, if the orthogonal part of the
similarities generate a finite group (see Section \ref{sec:Examples}). 
\begin{prop}
\label{pro:homogeneous-implies-usm}If $\mu$ is a homogeneous measure
then it is a USM and generates an EFD $P$ supported on homogeneous
measures.
\end{prop}
With these facts in hand it is not hard to use our results from the
previous sections to deduce:
\begin{thm}
\label{thm:properties-of-homogeneous-measures}If $\mu$ is a homogeneous
measure then $\pi\mapsto\dim\pi\mu$ is lower semi-continuous; for
every $\pi\in\Pi_{d,k}$ the image $\pi\mu$ and a.e. conditional
fiber measure is exact dimensional, and furthermore a.e. fiber measure
is a USM; and $\mu$ satisfies dimension conservation for every $\pi\in\Pi_{d,k}$,
in the sense of \eqref{eq:dimension-conservation}.
\end{thm}
We next consider EFDs with additional geometric invariance. For any
$k$ there is an action of $GL_{d}(\mathbb{R})$ on $\Pi_{d,k}$ given
by $U:\pi\rightarrow\pi\circ U^{-1}$, and a $GL_{k}(\mathbb{R})$-action
given by $V:\pi\rightarrow V\circ\pi$. These actions commute, giving
a $GL_{k}(\mathbb{R})\times GL_{d}(\mathbb{R})$-action.

Let $A\subseteq GL_{d}(\mathbb{R})$ be a group of linear transformations.
$A$ induces an action on measures and hence an action $A^{*}$ on
distributions. An EFD $P$ is non-singular with respect to $A$ if
$a^{*}P\sim P$ for every $a\in A$.

As a consequence of the semi-continuity of $E_{P}(\cdot)$, we obtain
the following:
\begin{prop}
\label{pro:measures-with-additional-invariance}Let $P$ be an EFD
which is non-singular with respect to a group $A\subseteq GL_{d}(\mathbb{R})$.
Then $E_{P}(\cdot)$ is constant on $\overline{A}$-orbits of $\Pi_{d,k}$.
In particular if an orbit $\Lambda\subseteq\Pi_{d,k}$ of $GL_{k}(\mathbb{R})\times\overline{A}$
has non-empty interior then $E_{P}|_{\Lambda}=\min\{k,\dim P\}$. 
\end{prop}
On way in which FDs with additional invariance arise is from USMs
with additional invariance. For example suppose that $\mu$ is a USM
generating $P$, and $f\mu\ll\mu$ for some $C^{1}$ diffeomorphism
$f$. Then for $f\mu$-a.e. $y$, since $\mu$ generates $P$ at $y$
(Proposition \ref{pro:equivalent-measures-have-same-sceneries}) and
$\mu$ generates $P$ at $x=f^{-1}(y)$, by Proposition \ref{pro:smooth-images-of-usms},
so $P=Df(x)^{*}P$ and $P$ is invariant under $Df(x)^{*}$. Using
this idea and the previous proposition, we can recover many of the
main results from \cite{HochmanShmerkin09}. See Theorem \ref{thm:HS-2}
below.

We end this section with another application of this idea which arose
in discussions with P. Shmerkin. Let us first recall that two exact-dimensional
measures $\mu,\nu\in\mathcal{P}(\mathbb{R})$ are said to \emph{dissonate}
if $\dim(\mu*\nu)=\min\{1,\dim\mu+\dim\nu\}$. This is the behavior
expected of ``generic'' pairs of measures but is quite difficult
to verify in specific cases, and does sometimes fail. One might expect,
in fact, that typical measures $\mu$ should dissonate with \emph{every}
measure $\nu$. The only non-trivial%
\footnote{Measures of dimension $0$ or $1$ dissonate with every measure for
trivial reasons.%
} examples known of this behavior are Salem measures, i.e. measures
$\mu$ for which the Fourier transform satisfies $|\widehat{\mu}(t)|=O(t^{-\dim\mu/2-\varepsilon})$
for all $\varepsilon>0$, but this condition is hard to verify. Examples
of such sets primarily come from randomized constructions.

Recall that a measure $\nu$ is Ahlfors-regular if there are constants
$\alpha,c>0$ such that $c^{-1}\cdot r^{\alpha}\leq\nu(B_{r}(x))\leq c\cdot r^{\alpha}$
for every $x\in\supp\nu$. We say that it is a.e.-Ahlfors if this
holds at $\nu$-a.e. $x$. For definitions of self-similar measures
see Section \ref{sub:Self-similar-measures}. The following result
shows that a very natural class of measures dissonates with every
Ahlfors-regular measure.
\begin{thm}
\label{thm:dissonance-with-Ahlfors}Let $\mu\in\mathcal{P}(\mathbb{R})$
be a self-similar measure defined by an IFS satisfying strong separation
and such that the contraction ratios generate a dense multiplicative
subgroup of $\mathbb{R}$. Then $\mu$ dissonates with every a.e.-Ahlfors
measure $\nu$, i.e. $\dim\mu*\nu=\min\{1,\dim\mu+\dim\nu\}$.
\end{thm}
The proof is in Section \ref{sub:EFDs-invariant-under-linear transformations}.
We do not know whether the regularity condition can be dropped; it
is at least possible that measures $\mu$ as in the theorem dissonate
with every measure $\nu$. It does not seem to be known whether such
$\mu$ are Salem measures, but it seems unlikely that they are.

\subsection{\label{sec:Remarks-and-problems}Open problems}

We collect here some open problems connected with this work (some
of which we already mentioned above).
\begin{enumerate}
\item Is it true that if $P$ is a FD then $P$-a.e. $\mu$ has exact dimensional
projections for every $\pi\in\Pi_{d,k}$? Exact dimensional conditionals?
dimension conservation?
\item What are the limits of scenery distributions of self-affine measures
and how do they relate to the structure of the measure?
\item What class of EFDs has absolutely continuous projections? We note
that the semi-continuity result for dimension of projections does
not have an analogue for absolute continuity, as is evident from \cite{Patzschke04}.
However there are various special measures for which absolute continuity
has been verified. It seems possible that some mixing assumption,
perhaps quantitative, on the FD might yield similar results (perhaps
what is required is some form of temporal or spatial mixing, or both).
\item More generally, what dynamical properties of an EFD as a measure preserving
system have implications for the geometry of the fractal? For one
result of this type see \cite{Hochman2010}.
\item Does Theorem \ref{thm:dissonance-with-Ahlfors} hold if we drop the
regularity assumptions on $\nu$?
\end{enumerate}

\subsection{Acknowledgments}

This project evolved in tandem with my work with Pablo Shmerkin on
\cite{HochmanShmerkin09}, and I'd like to thank Pablo for many interesting
discussions and references. I am also indebted to Hillel Furstenberg
and Matan Gavish for sharing with me their ideas about CP-processes
and USMs. Finally, I thank the referee for a careful reading.

\subsection{\label{sub:Summary-of-notation}Summary of notation.}

For the reader's convenience we summarize our main notation in the
table below.

\medskip{}
\begin{tabular}{lll}
\hline 
$d$  &  & Dimension of the ambient Euclidean space.\tabularnewline
$B_{r}(x)$ , $B_{r}$ &  & The closed ball of radius $r$ around $x$ (if not specified, $x=0$).\tabularnewline
$\mathcal{M}\,,\,\mathcal{M}(\mathbb{R}^{d})$ &  & Space of Radon measure on $\mathbb{R}^{d}$.\tabularnewline
$\mathcal{P}(X)$ &  & Space of probability measures on $X$.\tabularnewline
$\mu,\nu,\eta,\theta$ &  & Measures (elements of $\mathcal{M}$).\tabularnewline
$P,Q,R$ &  & Distributions (elements of $\mathcal{P}(\mathcal{M})$).\tabularnewline
$\mathcal{U},\mathcal{V},\mathcal{W}$ &  & Subsets of large spaces, e.g.$\mathcal{M}$ or $\mathcal{P}(\mathcal{M})$.\tabularnewline
$\mu^{*}$  &  & $\mu$ normalized to mass 1 on $[-1,1]^{d}$.\tabularnewline
$\mu^{\sqr}$ &  & $\mu$ restricted and normalized to mass 1 on $[-1,1]^{d}$.\tabularnewline
$\mu_{A}$ &  & Conditional measure of $\mu$ on $A$.\tabularnewline
$U_{x}$ , $U_{x}^{*}$ , $U_{x}^{\sqr}$ &  & $U_{x}(y)=y-x$, and normalized variants.\tabularnewline
$T_{B}$ , $T_{B}^{*}$ , $T_{B}^{\sqr}$ &  & Homothety mapping a cube $B$ onto $[-1,1]^{d}$, normalized variants\tabularnewline
$S_{t}$ , $S_{t}^{*}$ , $S_{t}^{\sqr}$ &  & $S_{t}(x)=e^{t}x$, and normalized variants.\tabularnewline
$\mathcal{D}_{b},\mathcal{D}_{b}(x)$  &  & Partition into $b$-adic cells and the cell containing $x$.\tabularnewline
$M_{b}$ , $M{}_{b}^{*}$ , $M_{b}^{\sqr}$ &  & Base-$b$ magnification operator, normalized variants.\tabularnewline
 $\mu_{x,t}^{\sqr}$, $\mu_{x,t}^{*}$  &  & Scenery of $\mu$ at $x$ at scale $t$ (restricted or not).\tabularnewline
$\left\langle \mu\right\rangle _{x,T}\,,\,\left\langle \mu,x\right\rangle _{T}$ &  & Continuous and $b$-scenery distributions at $x$ and scale $T$.\tabularnewline
$\Pi_{d,k}$  &  & Space of linear maps $\mathbb{R}^{d}\rightarrow\mathbb{R}^{k}$\tabularnewline
$\underline{D}_{\mu}(x)\,,\,\overline{D}_{\mu}(x)$ &  & Lower and upper local dimension of $\mu$ at $x$.\tabularnewline
$\ldim$ , $\dim$  &  & Lower Hausdorff dimension and exact dimension of a measure.\tabularnewline
\hline 
\end{tabular}

\section{\label{sec:Ergodic-theory}Preliminaries}

We collect here some ergodic-theoretic and measure-theoretic background.
For more information see \cite{Walters82,Glasner2003}.

\subsection{\label{sub:Measure-preserving-systems}Measure preserving systems}

A measure preserving system is quadruple $(\Omega,\mathcal{B},P,S)$,
where $(\Omega,\mathcal{B},P)$ is a standard probability space and
$S$ is a semigroup or group acting on $\Omega$ by measure-preserving
transformations: if $s\omega$ denotes the action of $s\in S$ on
$\omega\in\Omega$, then $s:\Omega\rightarrow\Omega$ is measurable
and $P(s^{-1}A)=P(A)$ for every $A\in\mathcal{B}$ and $s\in S$.
We usually drop $\mathcal{B}$ in our notation, and sometimes abbreviate
the system to $(P,S)$.

We will only encounter the cases $S=\mathbb{Z}$ or $\mathbb{Z}^{+}$,
which we call discrete time systems, and $S=\mathbb{R}$ or $\mathbb{R}^{+}$,
which we call continuous time systems, or flows. In the discrete time
case one often writes the system as $(\Omega,\mathcal{B},P,s)$ where
$s$ is a generator of $S$.

\subsection{Ergodicity and ergodic decomposition}

A measure preserving system is \emph{ergodic} if the only invariant
sets are trivial, i.e. if 
\[
\left(\forall s\in S\quad P(s^{-1}A\Delta A)=0\right)\qquad\implies\qquad\left(P(A)=0\mbox{ or }1\right)
\]
The ergodic decomposition theorem asserts that for any measure preserving
system $(\Omega,\mathcal{B},P,S)$ there is a map $\Omega\rightarrow\mathcal{P}(\Omega)$,
denoted $\omega\mapsto P_{\omega}$, which is (i) measurable with
respect to the sub-$\sigma$-algebra $\mathcal{I}\subseteq\mathcal{B}$
of $S$-invariant sets, (ii) $P=\int P_{\omega}dP(\omega)$, (iii)
$P$-a.e. $P_{\omega}$ is invariant and ergodic for $S$ and supported
on the atom of $\mathcal{I}$ containing $\omega$. Both $P_{\omega}$
and the atom containing $\omega$ are called the ergodic components
of $\omega$. The map $\omega\rightarrow P_{\omega}$ is unique up
to changes on a set of $P$-measure zero.

\subsection{\label{sub:Suspenssions}Time-1 maps and suspensions}

Given a continuous time system $(\Omega,\mathcal{B},P,S)$ and $s_{0}\in S$
one can consider the discrete-time systems $(\Omega,\mathcal{B},P,s_{0})$
in which the action is by the semigroup (or group if $s_{0}$ is invertible)
generated by $s_{0}$. This system need not be ergodic even if the
$S$-system is.

Conversely there is a standard construction to go from a discrete-time
system to a continuous-time one. Given a measure-preserving transformation
$s$ of $(\Omega,\mathcal{B},P)$, consider the product system $\Omega\times[0,c]$
with measure $P\times\lambda_{[0,c)}$, where $\lambda_{[0,c)}$ is
normalized Lebesgue measure on $[0,c)$, and for $t\in\mathbb{R}^{+}$
define $s_{t}(x,r)=(s^{[t/c]}x,c\{\frac{1}{c}(r+t)\})$, where $[u]$,
$\{u\}$ denote the integer and fractional parts of $u.$ Then $S=(s_{t})_{t\in\mathbb{R}}$
is a measure-preserving flow called the $c$\emph{-suspension }of
$(\Omega,\mathcal{B},P,s)$. Note that under $s_{c}$ the flow decomposes
into ergodic components of the form $P\times\delta_{u}$, $u\in[0,c)$
and the action of $s_{c}$ on each ergodic component is by applying
$s$ to the first coordinate and fixing the second. 

A more general construction is the \emph{flow under a function} construction.
Given a discrete time system $(\Omega,\mathcal{B},P,s)$ and a positive
measurable function $f:\Omega\rightarrow\mathbb{R}^{+}$. Let $\Omega'\subseteq\Omega\times\mathbb{R}^{+}$
denote the set 
\[
\Omega'=\{(\omega,r)\,:\,0\leq r<f(\omega)\}\subseteq\Omega\times\mathbb{R}
\]
and on it put the probability measure $P'=(P\times\lambda)_{\Omega'}$
(this requires $\int f\, dP<\infty$ in order that $P\times\lambda(\Omega')<\infty$).
Define a flow on this set by flowing vertically ``up'' from $(\omega,r)$
at unit speed until the second coordinate reaches $f(\omega)$, then
jumping to $(s\omega,0)$ and continuing to flow up. Formally, one
considers $\Omega'$ to be a factor space of $\Omega\times\mathbb{R}$
by the equivalence relation $\sim$ generated by $(\omega,r)\sim(s\omega,r-f(\omega))$.
Define a flow on $\Omega\times\mathbb{R}$ by $s_{t}(\omega,r)=(\omega,t+r)$.
This flow factors to the desired $P'$-preserving flow on $\Omega'$.
Note that the $c$-suspension is just a flow under the function $f\equiv c$.

\subsection{\label{sub:Isomorphism-and-factors}Isomorphism, factors, natural
extensions and processes}

If $(\Omega,\mathcal{B}),(\Omega',\mathcal{B}')$ are measurable spaces
and $S$ is a semigroup acting measurably on both then a factor map
$\pi:\Omega\rightarrow\Omega'$ is a measurable map which intertwines
the actions, i.e.
\[
s\circ\pi=\pi\circ s\qquad\mbox{for all }s\in S
\]
When the actions preserve measures $P,P'$ respectively the factor
map is also required to preserve measure, i.e. $P(\pi^{-1}A)=P'(A)$
for all $A\in\mathcal{B}'$. 

If $(\Omega,\mathcal{B},P,S)$ is a system with $S=\mathbb{Z}^{+}$
or $\mathbb{R}^{+}$, the natural extension of the system is the $\mathbb{Z}$
or $\mathbb{R}$ system, respectively, obtained as follows. Let $\overline{S}$
denote the (abstract) group generated by $S$ and take 
\[
\overline{\Omega}=(\omega\in\Omega^{\overline{S}}\,:\,\omega_{t+s}=s\omega_{t}\mbox{ for all }t\in\overline{S}\mbox{ and }s\in S\}
\]
Let $S$ act on $\overline{\Omega}$ by translation: $(s\omega)_{t}=\omega_{s+t}$.
Then $\pi:\omega\mapsto\omega_{0}$ is a factor map from $\overline{\Omega}$
with the translation action and $(\Omega,S)$, and there is a unique
invariant measure $\overline{P}$ on $\overline{\Omega}$ which projects
to $P$. The system $(\overline{\Omega},P,\overline{S})$ is the natural
extension of $(\Omega,P,S)$ and is characterized up to isomorphism
by the property that if $(\Omega',P',\overline{S})$ is another $\overline{S}$-system
and $\pi':(\Omega',S)\rightarrow(\Omega,S)$ is a factor map (with
respect to $S$) then there is a factor map $\varphi:(\Omega',\overline{S})=(\overline{\Omega},\overline{S})$
such that $\pi'=\pi\varphi$. 

It will sometimes be useful to identify a dynamical system with a
process. Given $(\Omega,\mathcal{B},P,S)$ we define the $\Omega$-valued
random variables $(X_{s})_{s\in S}$ on $(\Omega,\mathcal{B},P)$
by $X_{s}(\omega)=s\omega$. The family $(X_{s})_{s\in S}$ is a stationary
process in the sense that the joint distribution of any $k$-tuple
$(X_{s_{i}})_{i=1}^{k}$ is the same as the joint distribution of
$(X_{s_{i}+t})_{i=1}^{k}$ for every $t\in S$.

\subsection{\label{sub:The-ergodic-theorem}The ergodic theorem}

Let $(\Omega,\mathcal{B},P,(S_{t})_{t\in\mathbb{R}^{+}})$ be a measure-preserving
system and let $\mathcal{I}\subseteq\mathcal{B}$ be the $\sigma$-algebra
of $S$-invariant sets. Then for $f\in L^{1}(\Omega,P)$,
\[
\lim_{T\rightarrow\infty}\frac{1}{T}\int_{0}^{T}f\circ S_{t}\, dt=\mathbb{E}(f\,|\,\mathcal{I})
\]
$P$-a.e. and in $L^{1}(P)$. For ergodic systems, $\mathcal{I}$
is the trivial $\sigma$-algebra consisting of sets of measure $0$
and $1$, and the right hand side is then $\int f\, dP$. For discrete
time systems the same result holds with the integral is replaced by
a sum from $1$ to $T$.

\subsection{\label{sub:Generic-points}Generic points in topological systems}

Let $S:X\rightarrow X$ be a continuous transformation of a compact
metric space and let $P$ be a Borel probability measure on $X$.
A point $x\in X$ is \emph{generic} for $P$ if 
\begin{equation}
P_{N}=\frac{1}{N}\sum_{n=0}^{N-1}\delta_{S^{n}x}\rightarrow P\label{eq:orbital averages}
\end{equation}
in the weak-{*} topology; that is, for every $f\in C(X)$,
\begin{equation}
\int f\, dP_{N}=\frac{1}{N}\sum_{n=0}^{N-1}f(S^{n}x)\rightarrow\int fdP\label{eq:generic-point-condition}
\end{equation}

If $P$ is an $S$-invariant distribution then $P$-a.e. $x$ is generic
for $P$. Indeed, in order for $x\in X$ to be generic for $P$ it
is enough that \eqref{eq:generic-point-condition} hold for $f$ in
some fixed dense countable family $\mathcal{F}\subseteq C(X)$. Such
a family exists because $X$ is compact and metric. For $f\in\mathcal{F}$
and for $P$-a.e. $x$, \eqref{eq:generic-point-condition} holds
by the ergodic theorem. 

Conversely, if $x$ is generic for $P$ then $P$ is $S$ invariant,
and, more generally any accumulation point $P$ of the averages \eqref{eq:orbital averages}
is $S$-invariant. Indeed, invariance is equivalent to the equality
$\int f\circ S\, dP=\int f\, dP$ for every $f\in C(X)$. Since $S\delta_{y}=\delta_{Sy}$,
for $f\in C(X)$ we have 
\begin{eqnarray*}
\int f\circ S\, dP-\int f\, dP & = & \int(f\circ S)\, d(\lim_{N\rightarrow\infty}\frac{1}{N}\sum_{n=0}^{N-1}\delta_{S^{n}x})-\int f\, d(\lim_{N\rightarrow\infty}\frac{1}{N}\sum_{n=0}^{N-1}\delta_{S^{n}x})\\
 & = & \lim_{N\rightarrow\infty}\frac{1}{N}\left(\sum_{n=0}^{N}f(S(\delta_{S^{n}x}))-\sum_{n=0}^{N}f(\delta_{S^{n}x})\right)\\
 & = & \lim_{N\rightarrow\infty}\frac{1}{N}\left(\sum_{n=0}^{N}f(\delta_{S^{n+1}x}))-\sum_{n=0}^{N}f(\delta_{S^{n}x})\right)\\
 & = & \lim_{N\rightarrow\infty}\left(\frac{1}{N}f(\delta_{S^{N+1}x})-\frac{1}{N}f(\delta_{x})\right)\\
 & = & 0
\end{eqnarray*}
Note that we used continuity of $S$ to deduce that $f\circ S\in C(X)$,
which was used in passing from the first to the second line.

Generic points for a continuous time action of $\mathbb{R}$ on $X$
are defined similarly.

\subsection{\label{sub:Measure-valued-integration}Measure-valued integration}

Let $(\Omega,\mathcal{B},P)$ be a probability space, $X$ a compact
metric space, and let $\mu:\Omega\rightarrow\mathcal{P}(X)$, $\omega\mapsto\mu_{\omega}$
be a measure-valued function. We say that $\mu$ is measurable if
$\omega\mapsto\mu_{\omega}(A)$ is measurable for every Borel subset
$A\subseteq X$. This is the same as requiring that $\omega\mapsto\int f(x)d\mu_{\omega}(x)$
is measurable for every $f\in C(X)$, or that this is true for every
$f$ in a set $\mathcal{F}\subseteq C(X)$ that is dense in the supremum
norm. The integral $\int\mu_{\omega}dP(\omega)$ is then defined to
be the measure $\mu\in\mathcal{P}(X)$ satisfying
\[
\mu(A)=\int\mu_{\omega}(A)dP(\omega)\qquad\mbox{for Borel subsets }A\subseteq X
\]
or equivalently
\[
\int fd\mu=\int\left(f(x)d\mu_{\omega}(x)\right)dP(\omega)\qquad\mbox{for }f\in C(X)
\]

Many definitions and theorems of real analysis lead to measure-valued
analogs by integrating against test functions $f\in C(X)$. We give
two examples that we will use later.

First for $\mu:\Omega\rightarrow\mathcal{P}(X)$ as before and $\mathcal{C}\subseteq\mathcal{B}$
a sub-$\sigma$-algebra, the measure-valued conditional expectation
$\mathbb{E}(\mu|\mathcal{C})$ is the a.s. unique $\mathcal{C}$-measurable
measure-valued function $\omega\mapsto\nu_{\omega}\in\mathcal{P}(X)$
such that 
\[
\int f\, d\nu_{\omega}=\mathbb{E}(F_{f}|\mathcal{C})(\omega)\qquad P\mbox{ -a.e.}\omega
\]
where the right hand side is the real-valued conditional expectation
of the function $F_{f}(\omega)=\int fd\mu_{\omega}$. To prove that
such $\nu_{\omega}$ exist, fix a countable norm-dense $\mathbb{Q}$-linear
subspace $\mathcal{F}\subseteq C(X)$. Let $\Lambda:\mathcal{F}\rightarrow\mathbb{R}$
be the function $\Lambda_{\omega}(f)=\mathbb{E}(F_{f}|\mathcal{C})(\omega)$.
For $f\in\mathcal{F}$ this function is defined for $P$-a.e. $\omega$,
and hence a.e. for all $f\in\mathcal{F}$, and since conditional expectation
is linear and $\mathcal{F}$ and $\mathbb{Q}$ is countable, the countable
number of relations that ensure $\mathbb{Q}$-linearity and positivity
of $\Lambda_{\omega}$ hold $P$-a.e., we conclude that $\Lambda_{\omega}$
is $\mathbb{Q}$-linear and positive $P$ -a.s. Now, $\Lambda$ extends
to a real-linear positive operator on $C(X)$ because $\mathcal{F}$
is dense in $C(X)$, and the Riesz representation theorem gives us
a measure $\nu_{\omega}$ such that $\Lambda_{\omega}(f)=\int f\, d\omega$
for every $f\in\mathcal{F}$. In particular $\int fd\nu_{\omega}=\mathbb{E}(F_{f}|\mathcal{C})(\omega)$
for $f\in\mathcal{F}$. This relation extends from $\mathcal{F}$
to $\overline{\mathcal{F}}=C(X)$ by density of $\mathcal{F}$ and
continuity of conditional expectation on $C(X)$ (with respect to
the sup norm), as required. 

Second, let us prove the martingale convergence theorem for measure-valued
martingales. Suppose $\mu_{n}:\Omega\rightarrow\mathcal{P}(X)$, $\omega\mapsto\mu_{n,\omega}$,
is a martingale with respect to a filtration $\mathcal{C}_{1}\subseteq\mathcal{C}_{2}\subseteq\ldots\subseteq\mathcal{B}$,
i.e. $\mathbb{E}(\mu_{n}|\mathcal{C}_{n-1})=\mu_{n-1}$ a.s. Then
for every $f\in C(X)$ the functions $F_{f,n}:\Omega\rightarrow\mathbb{R}$
given by $F_{f,n}(\omega)=\int fd\mu_{n,\omega}$ form a bounded real-valued
martingale for the same filtration, since by definition of the measure-valued
conditional expectation, $\mathbb{E}(F_{f,n}|\mathcal{C}_{n-1})(\omega)=\int fd\mu_{n-1,\omega}=F_{f,n-1}(\omega)$.
Hence by the real-valued martingale convergence theorem, $F_{f,n}\rightarrow F_{f}$
a.s. for some $F_{f}:\Omega\rightarrow\mathbb{R}$. The convergence
holds simultaneously for all $f$ in a countable dense $\mathbb{Q}$-linear
subspace of $C(X)$, so we can assume as before that $f\mapsto F_{f}(\omega)$
is a.s. $\mathbb{Q}$-linear. Therefore there is a $\mu:\Omega$$\rightarrow\mathcal{P}(X)$
such that $\int fd\mu=F_{f}(\omega)$ for $P$-a.e. $\omega$, which
just means that $\int fd\mu_{\omega}=\lim\int fd\mu_{n,\omega}$ a.e.
on a dense set of $f$ and so for every $f$, so $\mu_{n}\rightarrow\mu$
a.e..

\subsection{\label{sub:Metrics-on-spaces-of-measures}Metrics on spaces of measures}

For a compact metric space $X$ it is convenient to introduce an explicit
metric $d(\cdot,\cdot)$ on $\mathcal{P}(X)$. Let $\lip(f)$ denote
the the Lipschitz constant (possibly $\infty$) of a function $f:X\rightarrow\mathbb{R}$.
For $\alpha,\beta\in\mathcal{P}(X)$ let 
\[
d(\alpha,\beta)=\sup\left\{ \left|\int fd\alpha-\int fd\beta\right|\,:\,\lip(f)\leq1\right\} 
\]
This metric is compatible with the weak topology (e.g. \cite[Chapter 14]{Mattila95})
and in particular is separable, complete and compact. We will sometimes
write $\mu=\nu+O(1)$ or $\mu=\nu+o(1)$, to mean that $d(\mu,\nu)=O(1)$
or $d(\mu,\nu)=o(1)$, respectively, etc. Note that if $\mu=(1-\varepsilon)\mu_{1}+\varepsilon\mu_{2}$
where $\mu_{1},\mu_{2}\in\mathcal{P}(X)$, then $\mu=\mu_{1}+O(\varepsilon)$.

Note that if $\mathcal{E}$ is a partition of $X$ into sets of diameter
$r>0$, and $\alpha,\beta\in\mathcal{P}(X)$ have the property that
$\alpha(E)=\beta(E)$ for all $E\in\mathcal{E}$, then $\left|\int fd\alpha-\int fd\beta\right|\leq cr$,
for all functions with $\lip(f)\leq1$, where $c=\diam B_{1}$ is
independent of $\alpha,\beta,r$. Consequently $d(\alpha,\beta)\leq cr$.
In particular, for $\nu,\mu\in\mathcal{M}^{\sqr}$,
\begin{equation}
\left(\forall A\in\mathcal{D}_{n}\quad\mu(A)=\nu(A)\right)\qquad\implies\qquad d(\mu,\nu)=O(\frac{1}{n})\label{eq:23}
\end{equation}

We apply this primarily in the cases $X=\mathcal{M}^{\sqr}=\mathcal{P}([-1,1]^{d})$
and $X=\mathcal{P}(\mathcal{M}^{\sqr})$, where in the latter case
the base space $\mathcal{M}^{\sqr}$ is given the metric from the
former case. Then if $(\beta_{n})_{0\leq n<N}$ and $(\gamma_{n})_{0\leq n<N}$
are sequences in $\mathcal{M}^{\sqr}$ and we form the uniform distribution
on them, $P_{\beta}=\frac{1}{N}\sum_{n=0}^{N-1}\delta_{\beta_{n}}$
and $P_{\gamma}=\frac{1}{N}\sum_{n=0}^{N-1}\delta_{\gamma_{n}}$,
then 
\begin{eqnarray}
\widetilde{d}(P_{\beta},P_{\gamma}) & = & \sup_{\lip(F)\leq1}\left\{ \frac{1}{N}\sum_{n=0}^{N-1}\int F\, d\delta_{\beta_{n}}-\frac{1}{N}\sum_{n=0}^{N-1}\int F\, d\delta_{\gamma_{n}}\right\} \nonumber \\
 & \leq & \frac{1}{N}\sum_{n=0}^{N-1}\sup_{\lip(F)\leq1}\left(\int F\, d\delta_{\beta_{n}}-\int F\, d\delta_{\gamma_{n}}\right)\nonumber \\
 & = & \frac{1}{N}\sum_{n=0}^{N-1}\sup_{\lip(F)\leq1}\left(F(\beta_{n})-F(\gamma_{n})\right)\nonumber \\
 & = & \frac{1}{N}\sum_{n=0}^{N-1}d(\beta_{n},\gamma_{n})\label{eq:distance-between averages}
\end{eqnarray}
A similar inequality holds for distributions of the form $\frac{1}{T}\int_{0}^{T}\delta_{\beta_{t}}dt$. 

Fixing the metric 
\[
d((\mu,x),(\nu,y))=\max\{d(\mu,\nu),d(x,y)\}
\]
on $\mathcal{M}^{\sqr}\times B_{1}$, for sequences $(\beta_{n},x_{n}),(\gamma_{n},y_{n})\in\mathcal{M}^{\sqr}\times[-1,1]^{d}$
the analogous inequality holds for the distance between the distributions
$P_{\beta}=\frac{1}{N}\sum_{n=0}^{N-1}\delta_{(\beta_{n},x_{n})}$
and $P_{\gamma}=\frac{1}{N}\sum_{n=0}^{N-1}\delta_{(\gamma_{n},y_{n})}$.
\begin{lem}
\label{lem:modulous-of-continuity-for-measures}Let $f:X\rightarrow Y$
be a map between metric spaces, and $\varepsilon,\delta>0$, such
that if $d(x,x')\leq\delta$ then $d(f(x),f(x'))<\varepsilon$. Then,
if $\alpha,\beta\in\mathcal{P}(X)$ and $d(\alpha,\beta)<\delta^{2}$,
then $d(f\alpha,f\beta)<\varepsilon+\delta\diam Y$.\end{lem}
\begin{proof}
This is most easily proved using the dual form of the metric
\[
d^{*}(\alpha,\beta)=\inf_{\gamma}\int d(x,x')\, d\gamma(x,x')
\]
where the infimum is over measures $\gamma\in\mathcal{P}(X\times X)$
whose projections to the first and second coordinates are $\alpha,\beta$,
respectively. Such $\gamma$ are called couplings of $\alpha,\beta$
and it is well known that $d=d^{*}$. By our assumption there is a
coupling $\gamma$ such that $d^{*}(\alpha,\beta)<\delta^{2}$. Consider
the map $F:X\times X\mapsto Y\times Y$ given by $(x,x')\mapsto(f(x),f(x'))$.
Let $\gamma'=F\gamma$. Then it is immediate to check that $\gamma'$
is a coupling of $f\alpha$ and $f\beta$ and hence 
\[
d(f\alpha,f\beta)=d^{*}(f\alpha,f\beta)\leq\int d(y',y'')\, d\gamma'(y',y'')=\int d(f(x),f(x'))\, d\gamma(x,x')
\]
Now, since by assumption $\int d(x,x')\, d\gamma(x,x')<\delta^{2}$,
by Markov's inequality, 
\[
\gamma\{(x,x')\,:\, d(x,x')>\delta\}<\delta
\]
hence
\begin{eqnarray*}
\int d(f(x),f(x'))\, d\gamma(x,x') & \leq & \varepsilon\cdot\gamma\{(x,x')\,:\, d(x,x')\leq\delta\}+\diam(Y)\cdot\gamma\{(x,x')\,:\, d(x,x')>\delta\}\\
 & < & \varepsilon(1-\delta)+\delta\diam(Y)\\
 & < & \varepsilon+\delta\diam Y
\end{eqnarray*}
as claimed.
\end{proof}

\section{\label{sec:Fractal-Distributions}Basic properties of the models}

This section contains some remarks on the definition of FDs and their
relation to other models, and derivation of some of their basic properties.

\subsection{\label{sub:remarks-on-normalization}Remarks on normalization}

If $\mu\in\mathcal{M}$ then (unless there is an atom at $x$) $\mu_{x,t}\rightarrow0$
as $t\rightarrow\infty$. Therefore some form of normalization is
necessary if we wish to study the dynamics of rescaling. However,
outside of some special cases (e.g. when there are second order densities,
see Section \ref{sub:Palm-and-Zahle-distributions}) there is no natural
way to do this. 

The most mathematically straightforward approach is to bypass the
issue of normalization altogether and work in the projective space
$\mathcal{M}/\mathbb{R}^{+}$, in which one identifies measures which
are constant multiples of each other, or, better yet, the factor space
in which equivalent measures are identified. However, this is somewhat
inconvenient in practice, and has few real benefits. 

Our choice of normalization, $\mu\mapsto\mu^{*}$, amounts to choosing
a section over $\mathcal{M}/\mathbb{R}^{+}$. We could of course choose
to normalize some other neighborhood of the origin $U$ using the
normalization $\mu^{U}=\frac{1}{\mu(U)}\mu$, and associated operations.
Although we stated in the introduction that FDs are defined in a coordinate-free
way, our choice of $B_{1}$ as the set on which measures are normalized
depends on both the coordinate system and the norm. However, if we
use $\mu^{U}$ instead then the extended FDs with respect to $U$
are in 1-1 correspondence with extended FDs as we have defined them.
Indeed since the maps $\mu\rightarrow\mu^{U}$ and $\nu\mapsto\nu^{*}$
are inverses of each other, defined on $\mathcal{M}^{U}=\{\mu^{U}\,:\,\mu\in\mathcal{M}\,:\,\mu(U)>0\}$
and $\mathcal{M}^{*}$, respectively, and they induce a bijections
of the spaces of distributions on these sets which, furthermore, intertwine
the respective scaling operators $S_{t}^{*},S_{t}^{\sqr}$. Also,
the normalization does not change the measure class of the measure
and so modulo equivalence of measures, the choice of $U$ is inconsequential.

\subsection{\label{sub:Restricted-and-extended-versions}Restricted and extended
versions}

The following claim holds generally for $S^{*}$-invariant distributions
(the quasi-Palm property is not needed).
\begin{lem}
\label{pro:restriction-of-distributions}The map $\mathcal{M}\rightarrow\mathcal{M}^{\sqr}$,
$\mu\mapsto\mu^{\sqr}$ intertwines the $\mathbb{R}^{+}$-actions
of $S^{*}$ and $S^{\sqr}$, i.e. $(S_{t}^{*}\mu)^{\sqr}=S_{t}^{\sqr}\mu$
for $t\geq0$, and induces a one-to-one correspondence $P\mapsto P^{\sqr}$
between $S^{*}$-invariant distributions on $\mathcal{M}^{*}$ and
$S^{\sqr}$-invariant distributions on $\mathcal{M}^{\sqr}$.\end{lem}
\begin{proof}
The equality $(S_{t}^{*}\mu)^{\sqr}=S_{t}^{\sqr}\mu$ is immediate,
which gives the first statement.

Suppose that $Q$ is an $S^{\sqr}$-invariant distribution; we construct
a $S^{*}$-invariant distribution $\widetilde{Q}$ with $\widetilde{Q}^{\sqr}=Q$,
showing that restriction is surjective. Begin with stochastic process
$(\mu_{n})_{n=0}^{\infty}$ defined by randomly selecting $\mu_{0}\sim Q$
and setting $\mu_{n}=S_{n}^{\sqr}\mu_{0}$. By invariance of $Q$,
this is a stationary process, and it has a two-sided extension, that
is, a stationary process $(\mu_{n})_{n=-\infty}^{\infty}$ such that
all finite marginals agree with the original process \cite[Lemma 10.2]{Kallenberg83}.
In particular, since the original process satisfied $S_{m}^{\sqr}\mu_{n}=\mu_{n+m}$
for all integers $m,n\geq0$, the same is true for the two-sided process
and by throwing out a zero-probability event we can assume it holds
everywhere. Now, for $t\in\mathbb{R}$ define $\mu_{t}=S_{t-[t]}^{\sqr}\mu_{[t]}$
(this is consistent with the previous definition for $t\in\mathbb{N}$).
Since $Q$ was initially $S_{t}^{\sqr}$-invariant for all $t$ (not
just integer $t$), it is routine to verify that now $(\mu_{t})_{t\in\mathbb{R}}$
is a stationary process and that it satisfies 
\begin{equation}
\mu_{t+s}=S_{s}^{\sqr}\mu_{t}\label{eq:natural-extention-of-restricted-FD}
\end{equation}
a.s. for all $t\in\mathbb{R}$ and $s\geq0$. 

For each $t\geq0$ define
\[
\nu_{t}=S_{t}^{*}\mu_{-t}
\]
Note that \eqref{eq:natural-extention-of-restricted-FD} implies that
for $t\geq s\geq0$, 
\begin{equation}
\nu_{t}|_{B_{e^{s}}(0)}=\nu_{s}\label{eq:compatability-of-extension-of-EFDs}
\end{equation}
Therefore we can define a measure $\mu\in\mathcal{M}$ by the condition
\begin{equation}
\nu(A)=\lim_{t\rightarrow\infty}\nu_{t}(A)\label{eq:construction-of-extended-FD}
\end{equation}
for all Borel sets $A$ (the limit exists since the $\nu_{t}(A)$
increases as $t\rightarrow\infty$). The verification that this is
a measure is standard.

Let $\widetilde{Q}$ be the distribution of $\nu$. Then 
\[
\nu^{\sqr}(A)=\lim_{t\rightarrow\infty}\frac{1}{\nu_{t}(B_{1})}\nu_{t}(A\cap B_{1})=\nu_{0}(A\cap B_{1})=S_{0}^{*}\mu_{0}(A\cap B_{1})=\mu_{0}(A)
\]
so $\widetilde{Q}^{\sqr}=Q$. Similarly, 
\begin{eqnarray*}
S_{s}^{*}\nu & = & c\cdot\nu\circ(S_{s})^{-1}\\
 & = & c\cdot(\lim_{t\rightarrow\infty}\nu_{t})\circ(S_{s})^{-1}\\
 & = & c\cdot\lim_{t\rightarrow\infty}(c'(S_{t}\mu_{-t})\circ(S_{s})^{-1})\\
 & = & cc'\lim_{t\rightarrow\infty}(S_{t}(S_{s}\mu_{-t}))\\
 & = & \lim_{t\rightarrow\infty}S_{t}^{*}(\mu_{-(t-s)})
\end{eqnarray*}
where we have used the definition of $\nu$ and \eqref{eq:compatability-of-extension-of-EFDs},
and $c,c'$ are normalizing constants. Since $\mu_{-(t-s)}$ has the
same distribution as $\mu_{-t}$, this shows that $S_{s}^{*}\nu$
has the same distribution as $\lim_{t\rightarrow\infty}S_{t}^{*}(\mu_{-t})=\nu$,
and so $\widetilde{Q}$ is $S^{*}$ invariant.

Finally, note that if $R$ is any other $S^{*}$-invariant distribution
with $R^{\sqr}=Q$ then $\tau\mapsto\tau|_{B_{e^{n}}(0)}$ for $t\geq0$
has the same distributions as $\nu_{n}$ (since it is related to $R^{\sqr}=Q$
in the same way $\nu_{n}$ is). It follows immediately that $R=\widetilde{Q}$.
\end{proof}
Next, for a restricted CP-distribution $Q$ we construct the extended
version $\widetilde{Q}$. This is an $M_{b}^{*}$-invariant distribution
$\widetilde{Q}$ on $\mathcal{M}^{*}\times B_{1}$ whose push-forward
via $(\mu,x)\mapsto(\mu^{\sqr},x)$ is $Q$, but note that these two
properties of $Q$ do not characterize it as in the case of the extended
version of a FD. An example is given below.

The construction is analogous to the construction of extended FDs.
Let $(\mu_{n},x_{n})_{n\in\mathbb{Z}}$ be a stationary process with
marginal $Q$ and such that 
\[
(\mu_{n+1},x_{n+1})=M_{b}^{\sqr}(\mu_{n},x_{n})
\]
(to see that such a process exists, begin again with the one-sided
process $((M_{b}^{\sqr})^{n}(\mu,x))_{n=0}^{\infty}$, $(\mu,x)\sim Q$,
take the two-sided extension, and note that since the original process
satisfied the equation above, up to a nullset the extended one does
too. For $n\geq0$ let 
\[
\nu_{n}=T_{\mathcal{D}_{b^{n}}(x_{-n})}^{*}\mu_{-n}
\]
and 
\[
E_{n}=T_{\mathcal{D}_{b^{n}}(x_{-n})}B_{1}
\]
then $\nu_{n}$ is supported on $E_{n}$ and  $B_{1}=E_{0}\subseteq E_{1}\subseteq E_{2}\subseteq\ldots$.
For $n\geq m\geq0$ we have 
\[
(M_{b}^{\sqr})^{m}(\mu_{-n},x_{-n})=(\mu_{-n+m},x_{-n+m})
\]
so
\[
\nu_{n}|_{E_{m}}=\nu_{m}
\]
and we may define a measure $\nu$ by
\[
\nu(A)=\lim_{n\rightarrow\infty}\nu_{n}(A)
\]
for all Borel sets $A$. This measure is easily seen to be Radon.
The distribution $\widetilde{Q}$ of $\nu$ is a distribution on $\mathcal{M}^{*}\times B_{1}$
and it is $M_{b}^{*}$-invariant. This is the extended version of
$Q$. One may easily verify that the map $(\mu,x)\mapsto(\mu^{\sqr},x)$
is a factor map between the measure preserving systems $(\mathcal{M}^{*}\times B_{1},\widetilde{Q},M_{b}^{*})$
and $(\mathcal{M}^{\sqr}\times B_{1},Q,M_{b}^{\sqr})$.

Here is an example showing that if $Q$ is a restricted CP-distribution
and $R$ is a distribution on $\mathcal{M}^{*}\times B_{1}$ which
is (i) invariant under the map $M_{b}^{*}$ and (ii) factors onto
$Q$ via $(\mu,x)\mapsto(\mu^{\sqr},x)$, it does not follow that
$R=\widetilde{Q}$. Indeed, consider the measure $\mu$ on $\mathbb{R}^{2}$
consisting of Lebesgue measure on the vertical line $x=1$, and $2^{-n}$
times Lebesgue measure on the vertical line $x=1+2^{-n}$, $n\in\mathbb{N}$.
Let $Q=\int\delta_{(\mu^{\sqr},x)}\, d\mu^{\sqr}(x)$, which is a
restricted CP-distribution (and in fact ergodic). Then, writing $I=\{1\}\times[-1,1]$,
the distribution $R=\int_{I}\delta_{(\mu,x)}\, d\mu(x)$ is adapted
and invariant under $M_{2}^{*}$, and factors onto the CP-distribution
$Q$, but $R$ is not the extended version of $Q$. Rather, the extended
version is $\widetilde{Q}=\int_{I}\delta_{(\nu,x)}\, d\nu(x)$, where
$\nu$ consists of Lebesgue measure on the line $x=1$.

The example above is rather special and also typical of what can go
wrong: 
\begin{lem}
Let $Q$ be a restricted ergodic CP-distribution and suppose that
$\int\theta(\partial B_{1})\, dQ(\theta,x)<1$. Then there is a unique
CP-distribution $\widetilde{Q}$ factoring onto $Q$ via $(\mu,x)\mapsto(\mu^{\sqr},x)$.
\end{lem}
Actually, for ergodic $Q$ we always have $\int\theta(\partial B_{1})\, dQ(\theta,x)=0$
or $1$. We do not use this or the lemma anywhere and omit the proofs.

\subsection{\label{sub:Palm-and-Zahle-distributions}Relation to Z\"{a}hle distributions}

Closely related to FDs are Z\"{a}hle's notion of an $\alpha$-scale-invariant
distribution \cite{Zahle88,MortersPreiss98}. For fixed $\alpha>0$,
Z\"{a}hle's notion of scale invariance uses the additive $\mathbb{R}$-action
$S^{\alpha}=(S_{t}^{\alpha})_{t\in\mathbb{R}}$ on $\mathcal{M}$,
defined by 
\[
S_{t}^{\alpha}\mu=e^{\alpha t}S_{t}\mu
\]
Also, a Palm distribution $P\in\mathcal{P}(\mathcal{M})$ is a distribution
with finite intensity, i.e.
\begin{equation}
\int\mu(K)\, dP(\mu)<\infty\label{eq:finite-intensity}
\end{equation}
for every compact $K\subseteq\mathbb{R}^{d}$, and which satisfies
\[
\int\left\langle \mu\right\rangle _{B}dP(\mu)=P
\]
for every non-empty ball $B$ centered at the origin. The last condition
differs from the definition of the quasi-Palm condition both in that
we have used the non-normalized diffusion operation $\left\langle \mu\right\rangle _{B}$,
and we required equality rather than equivalence of the distributions..

An $\alpha$-Z\"{a}hle distribution is an $S^{\alpha}$-invariant
distribution $P\in\mathcal{M}$ which in addition is a Palm distribution.
We say that $P$ it is degenerate if it gives positive mass to the
zero measure. We write $\mathcal{Z}_{\alpha}$ for the space non-degenerate
ergodic $\alpha$-Z\"{a}hle distributions. Note that by $S^{\alpha}$-invariance,
if $P$ is such a distribution then $P$-a.e. $\mu$ has $0\in\supp\mu$.

Define the $\alpha$-scenery of $\mu$ at $x$ to be
\[
\mu_{x,t}^{\alpha}=S_{t}^{\alpha}T_{x}\mu
\]
and the scenery distribution
\[
\left\langle \mu\right\rangle _{x,T}^{\alpha}=\frac{1}{T}\int_{0}^{T}\delta_{\mu_{x,t}^{\alpha}}\, dt
\]
which is a distribution on $\mathcal{M}$. In this space of distributions
we consider weak convergence, i.e. $P_{n}\rightarrow P$ if $\int f\, dP_{n}\rightarrow\int f\, dP$
for every bounded continuous $f:\mathcal{M}\rightarrow\mathbb{R}$.
Note that the weak topology is not compact, and it is entirely possible
that $\mu_{x,t}^{\alpha}$ becomes unbounded as $t\rightarrow\infty$
or converge to $0$. Thus the $\alpha$-scenery distributions may
not have convergent subsequences, or may converge to a degenerate
distribution. Nevertheless M\"{o}rters and Preiss proved the following
predecessor of Theorem \ref{thm:limit-distributions-are-FDs}:
\begin{thm}
\label{thm:MortersPreiss}[M\"{o}rters and Preiss, \cite{MortersPreiss98}]
If $\mu\in\mathcal{M}$ and $\alpha>0$ then for $\mu$-a.e. $x$,
every weak accumulation point of the $\alpha$-scenery at $x$ is
an $\alpha$-Z\"{a}hle distribution. 
\end{thm}
Let us now make precise the analogy between EFDs and non-degenerate
and ergodic $\alpha$-Z\"{a}hle distributions, i.e. the members of
$\mathcal{Z}_{\alpha}$. We consider the map $\mu\mapsto\mu^{*}$
and, first, note that if $P\in\mathcal{Z}_{\alpha}$ then $P^{*}$
is $S^{*}$-invariant, and $\mu\mapsto\mu^{*}$ is a factor map between
$(P,S^{\alpha})$ and $(P^{*},S^{*})$. Since $P$ is $S^{\alpha}$-ergodic,
$P^{*}$ is $S^{*}$-ergodic. It is also easy to check from the definition
that since $P$ is Palm, $P^{*}$ is quasi-Palm. Hence $P^{*}$ is
an EFD and $\mu\mapsto\pi\mu$ maps $\mathcal{Z}_{\alpha}$ to the
class of EFDs. This map is many-to-one for the following reason. Let
us say that distributions $P,Q$ are equivalent up to a constant,
and write $P\approx Q$, if there is a $c>0$ such that $Q$ is the
push-forward of $P$ through the map $\mu\mapsto c\mu$. Note that
if $P\in\mathcal{Z}_{\alpha}$ and $Q\approx P$ then $Q\in\mathcal{Z}_{\alpha}$,
and also $P^{*}=Q^{*}$. Therefore the fibers of the map map $P\mapsto P^{*}$
from $\mathcal{Z}_{\alpha}$ to the space of $S^{*}$-invariant distributions
are saturated with respect to the $\approx$ relation.

For a $\mu\in\mathcal{M}$ and $x\in\supp\mu$, the $\alpha$\emph{-dimensional
second-order density }of $\mu$ at $x$ is
\[
D_{\alpha}(\mu,x)=\lim_{T\rightarrow\infty}\frac{1}{T}\int_{0}^{T}e^{\alpha t}\mu(B_{t}(x))dt
\]
assuming the limit exists. Denote
\[
\mathcal{A}_{\alpha}=\{\mu\in\mathcal{M}\,:\,\mu\mbox{ has positive, finite }\alpha\mbox{-dimensional second order densities a.e.}\}
\]
Note that $D_{\alpha}(\mu,x)$ is homogeneous in $\mu$, that is,
given $c\in\mathbb{R}$ we have $D_{\alpha}(c\mu,x)=cD_{\alpha}(\mu,x)$
at $\mu$-a.e. $x$. Therefore $\mathcal{A}_{\alpha}$ is also saturated
with respect to $\approx$.
\begin{thm}
\label{thm:EFDs-and-Zahle-distributions}For each $\alpha>0$ the
map $Q\mapsto Q^{*}$ is a bijection between $\mathcal{Z}_{\alpha}/\approx$
and the space of EFDs supported on $\mathcal{A}_{\alpha}$.\end{thm}
\begin{proof}
Let $P\in\mathcal{Z}_{\alpha}$, i.e. it is a non-degenerate, ergodic
$\alpha$-Z\"{a}hle distributions. We have already noted that $P^{*}$
is an EFD, so we must show that $P^{*}(\mathcal{A}_{\alpha})=1$.
Note that $D_{\alpha}(\mu,0)$ is an ergodic average with respect
to $S^{\alpha}$ of the function $f(\mu)=\mu(B_{1})$, and since $P$
is Palm, by \eqref{eq:finite-intensity}, $\int\mu(B_{1})dP(\mu)<\infty$.
Therefore, for $P$-a.e. $\mu$ the $\alpha$-dimensional second-order
density at the origin exists and by ergodicity is $P$-a.e. equal
to $c=\int fdP$, which is positive and finite. Let $\mathcal{A}'\subseteq\mathcal{M}$
denote the set of measures whose $\alpha$-dimensional second-order
density at the origin exists and is equal to $c$. Since $P$ is Palm
and $P$-a.e. $\mu$ is in $\mathcal{A}$, for every $R>0$, $P$-a.e.
$\mu$ and $\mu$-a.e. $x\in B_{r}(0)$, also $T_{x}\mu\in\mathcal{A}'$.
Since this holds for all $R>0$, this implies that for $P$-a.e. $\mu$,
the $\alpha$-dimensional second-order density exists $\mu$-a.e.
and is equal to $c$, so $\mu\in\mathcal{A}'$. Since $\mu^{*}=\frac{1}{\mu(B_{1})}\mu$,
it follows that $\mu^{*}\in\mathcal{A}_{\alpha}$ for $P$-a.e. $\mu$;
equivalently, $P^{*}(\mathcal{A}_{\alpha})=1$, as desired. 

Conversely, let $Q$ be an extended EFD supported on $\mathcal{A}_{\alpha}$.
Let $\mathcal{A}''$ denote the set of measures with positive finite
$\alpha$-dimensional average density at $0$. Since $Q$ is quasi-Palm
property, since $Q(\mathcal{A}_{\alpha})=1$, also $Q(\mathcal{A}'')=1$,
and for $Q$-typical $\nu$, $D_{\alpha}(\nu,0)$ is well-defined.
Observe that 
\begin{eqnarray*}
D_{\alpha}(S_{t}^{*}\nu,0) & = & D_{\alpha}(\frac{1}{\nu(B_{e^{-t}})}S_{t}\nu,0)\\
 & = & D_{\alpha}(\frac{e^{-t}}{\nu(B_{e^{-t}})}S_{t}^{\alpha}\nu,0)\\
 & = & \frac{e^{-t}}{\nu(B_{e^{-t}})}D(S_{t}^{\alpha}\nu,0)\\
 & = & \frac{e^{-t}}{\nu(B_{e^{-t}})}D(\nu,0)
\end{eqnarray*}
where in the last step we used the elementary  fact that $D(\cdot,0)$
is $S_{t}^{\alpha}$-invariant when it exists. For $\nu\in\mathcal{A}_{\alpha}$,
define 
\[
\nu_{*}=\frac{1}{D_{\alpha}(\nu,0)}\nu
\]
By the above, 
\begin{eqnarray*}
(S_{t}^{*}\nu)_{*} & = & \frac{1}{D_{\alpha}(S_{t}^{*}\nu,0)}(S_{t}^{*}\nu)\\
 & = & \frac{e^{t}\nu(B_{e^{-t}})}{D_{\alpha}(\nu,0)}\cdot\frac{1}{\nu(B_{e^{-t}})}S_{t}\nu\\
 & = & \frac{1}{D_{\alpha}(\nu,0)}S_{t}^{\alpha}\nu\\
 & = & S_{t}^{\alpha}(\nu_{*})
\end{eqnarray*}
Thus $Q_{*}$ is an $S^{\alpha}$-invariant distribution, supported
on $\mathcal{A}_{\alpha}$, and, being a factor of $(Q,S^{*})$ via
$\nu\rightarrow\nu_{*}$, it is ergodic. 

Rather than showing directly that $Q_{*}$ is a Palm distribution
of finite intensity, we rely on the work of M\"{o}rters and Preiss.
Choose a $Q$-typical $\nu$ and $\nu$-typical $x$ and consider
the pairs $p_{t}=(\nu_{x,t}^{\alpha},\nu_{x,t}^{\sqr})\in\mathcal{M}(B_{1})\times\mathcal{M}^{\sqr}$
for $t\geq0$. Choose a sequence $T_{i}\rightarrow\infty$ such that
$\frac{1}{T_{i}}\int_{0}^{T_{i}}\delta_{p_{t}}dt\rightarrow R\in\mathcal{P}(\mathcal{M}(B_{1})\times\mathcal{M}^{\sqr})$
as $i\rightarrow\infty$. Since $Q$ is an EFD and $(\mu,x)$ are
typical, the marginal of $R$ on the second component is $Q$. We
may pass to a further subsequence such that $\left\langle \nu\right\rangle _{x,T_{i}}^{\alpha}$
converges weakly to an $\alpha$-Z\"{a}hle distribution $P$ (here
we use Theorem \ref{thm:MortersPreiss}, and also the fact that $D_{\alpha}(\nu,x)<\infty$
and \cite[Proposition 2]{MortersPreiss98}). With minor effort it
follows that $\left\langle \nu\right\rangle _{x,T_{i}}^{\alpha}|_{B_{1}}\rightarrow P^{1}$,
where $P^{1}$ is the push-forward of $P$ through $\theta\rightarrow\theta|_{B_{1}}$,
and this is the marginal of $R$ on the first component. By the previous
discussion, for every $t$ we have 
\[
\nu_{x,t}^{\alpha}|_{B_{1}}=(S_{t}^{\alpha}T_{x}\nu)|_{B_{1}}=D_{\alpha}(\nu,x)\cdot(S_{t}^{*}T_{x}\nu)_{*}|_{B_{1}}=c'\cdot(\nu_{x,t}^{\sqr})_{*}
\]
Now, although the map $\theta\mapsto c'\theta_{*}$ is not continuous,
a judicious use of Egorov's theorem allows us to approximate it on
sets of large $Q$-measure by a continuous one, and conclude that
$R$ is supported on pairs $(c''\cdot\theta_{*},\theta)$. Consequently
$(Q^{\sqr})_{*}=P^{1}$. Passing to extended versions, this means
that $Q_{*}=P$. 

In conclusion, for the EFD $Q$ we have found an $\alpha$-Z\"{a}hle
distribution $P$ such that $P=Q_{*}$ and $Q=P^{*}$. This completes
the proof.
\end{proof}
As a consequence of this characterization, we find that EFDs are a
far broader model than Z\"{a}hle distributions:
\begin{cor}
\label{cor:EFDs-are-more-general-than-Zahle}There exist EFDs which
do not arise from a Z\"{a}hle distribution. \end{cor}
\begin{proof}
Any EFD supported on measures which a.s. do not have positive second
order densities is such an example. To be concrete, we may take the
EFD associated to any self-similar measure arising from similarities
satisfying strong separation, and which is singular with respect to
Hausdorff measure at the appropriate dimension. See Section \ref{sec:Examples}.
The fact that such a measure does not have second order densities
follows from Patzschke and Z\"{a}hle \cite{PatzschkeZahle90}, and
that this is true for the corresponding EFD follows from homogeneity
(Sections \ref{sub:FDs-with-additional-invariance} and \ref{sub:Self-similar-measures}).

For a more concrete example consider the measure $\mu$ on $[0,1]$
which is the distribution of the random number $x$ whose binary digits
are chosen independently to be $0$ with probability $0<p<\frac{1}{2}$
and $1$ with probability $1-p$. The associated EFD is not a Z\"{a}hle
distribution, for the same reason as above.
\end{proof}
Finally, we remark that FDs do not necessarily have finite intensity;
that is, if $P$ is a FD then there may be compact $K\subseteq\mathbb{R}^{d}$
with $\int\nu(K)\, dP(\nu)=\infty$. By definition, this phenomenon
does not occur for Z\"{a}hle distributions.

\subsection{\label{sub:Ergodic-and-spatial-decomposition}Ergodic and spatial
decompositions}

Next, we prove the ergodic decomposition theorem for FDs. We  give
a self-contained geometric proof. A shorter proof using the identification
of FDs with CP-distributions, Theorems \ref{thm:CP-to-EFD} and \ref{thm:FD-to-CP},
and using properties of Markov processes, can be given along the lines
of the remark following Proposition 5.1 of \cite{Furstenberg08}

We require the following classical result, see \cite[Corollary 2.14]{Mattila95}:
\begin{lem}
[Besicovitch density theorem] \label{lem:Besicovitch-lemma}If $\mu$
is a Radon measure on $\mathbb{R}^{d}$ and $\mu(A)>0$ then for $\mu$-a.e.
point in $A$,
\[
\lim_{r\rightarrow0}\frac{\mu(A\cap B_{r}(x))}{\mu(B_{r}(x))}=1
\]

\end{lem}
We now prove the ergodic decomposition theorem for FDs:
\begin{proof}
[Proof of Theorem \ref{thm:EFD-ergodic-components}] Let $P$ be a
FD. Let $\widetilde{\mathcal{I}}$ denote the $\sigma$-algebra of
$S^{*}$-invariant sets in $\mathcal{M}^{*}$ and $\mathcal{I}\subseteq\widetilde{\mathcal{I}}$
 a countably-generated sub-$\sigma$-algebra that is dense in $\widetilde{\mathcal{I}}$
mod $P$ . Let $\mathcal{I}(\mu)$ denote the atom of $\mathcal{I}$
to which $\mu$ belongs, and let $P=\int P_{\mu}\, dP(\mu)$ denote
(an $\mathcal{I}$-measurable version of) the ergodic decomposition
of $P$, so that $\mu\mapsto P_{\mu}$ is $\mathcal{I}$-measurable
and for $P$-a.e. $\mu$ the distribution $P_{\mu}$ is $S^{*}$-invariant
and supported on $\mathcal{I}(\mu)$. 

We must show that $P_{\mu}$ is a.s. quasi Palm. It suffices to show
that for any open neighborhood $U$ of $0\in\mathbb{R}^{d}$ we have
$\int\left\langle \nu\right\rangle _{U}dP_{\mu}(\nu)\sim P_{\mu}$
for $P$-a.e. $\mu$. Indeed, if this holds for all $U$ in a countable
basis of neighborhoods of $0$, then for $P$-a.e. $\mu$ the equivalence
will hold for all neighborhoods $U$ of $0$ simultaneously, because
any $U$ can be written as $U=\bigcup_{i=1}^{\infty}U_{i}$ for sets
$U_{i}$ chosen from a fixed countable basis of the topology, and
the validity for $U_{i}$ implies it for $U$.

Fix an open neighborhood $U$ of the origin. Write 
\[
Q=\int\left\langle \nu\right\rangle _{U}^{*}\, dP(\nu)
\]
For $\mu$ such that $P_{\mu}$ is defined, let 
\[
Q_{\mu}=\int\left\langle \nu\right\rangle _{U}^{*}\, dP_{\mu}(\nu)
\]
so, integrating over $\mu\sim P$ and using $P=\int P_{\mu}\, dP(\mu)$,
we have 
\[
Q=\int Q_{\mu}\, dP(\mu)
\]
Since $P$ is quasi-Palm, we also know that 
\[
Q\sim P=\int P_{\mu}\, dP(\mu)
\]
Combining the last two relations, we have
\[
\int Q_{\mu}\, dP(\mu)\sim\int P_{\mu}\, dP(\mu)
\]
We know that $P_{\mu}$ is a.s. supported on $\mathcal{I}(\mu)$,
hence if we can show that $Q_{\mu}$ is also a.s. supported on $\mathcal{I}(\mu)$,
then we deduce that $Q_{\mu}\sim P_{\mu}$ a.s., as desired.

In order to establish that $Q_{\mu}$ is $P$-a.s. supported on $\mathcal{I}(\mu)$
it suffices to show that for every $\mathcal{V}\in\mathcal{I}$ with
$P(\mathcal{V})>0$, the distribution $\int_{\mathcal{V}}Q_{\mu}\, dP(\mu)=\int_{\mathcal{V}}\left\langle \nu\right\rangle _{U}^{*}\, dP(\nu)$
is supported on $\mathcal{V}$ (here we use that $\mathcal{I}$ is
countably generated). To this end, fix $\mathcal{V}$ as above. Given
$\mu\in\mathcal{M}$ let 
\[
V_{\mu}=\{x\in\mathbb{R}^{d}\,:\, T_{x}^{*}\mu\in\mathcal{V}\}
\]
Note that if $x\in V_{\mu}$ then $\mu_{x,t}^{*}\in\mathcal{V}$ for
all $t$, because $\mathcal{V}$ is $S^{*}$-invariant. For $n\in\mathbb{N}$
let
\[
f_{n}(\mu)=\frac{\mu(\{x\in B_{n}(0)\,:\, T_{x}^{*}\mu\in\mathcal{V}\})}{\mu(B_{n}(0))}
\]
Applying the Besicovitch density lemma to the set $V_{\mu}$, we find
that for $\mu$-a.e. $x$, 
\[
\lim_{t\rightarrow\infty}f_{n}(\mu_{x,t}^{*})\in\{0,1\}\qquad\mbox{for }\mu\mbox{-a.e. }x
\]
By the quasi-Palm property, for $P$-a.e. $\mu$ this holds for $x=0$,
i.e.
\[
\lim_{t\rightarrow\infty}f_{n}(S_{t}^{*}\mu)\in\{0,1\}\qquad\mbox{for }P\mbox{-a.e. }\mu
\]
This in turn implies that 
\[
f_{n}(\mu)\in\{0,1\}\qquad\mbox{for }P\mbox{-a.e. }\mu
\]
Indeed, for $\varepsilon>0$ consider the event $\mathcal{A}_{\varepsilon}=\{f_{n}\geq1-\varepsilon\mbox{ or }f_{n}\leq\varepsilon\}\subseteq\mathcal{M}$.
Since $S_{t}^{*}\mu\in\mathcal{A}_{\varepsilon}$ for all large enough
$t$, ergodicity implies that $\mathcal{A}_{\varepsilon}$ has $P$-measure
$1$. Since $\varepsilon$ was arbitrary, this means $\mu\in\mathcal{A}_{0}$
for $P$-a.e. $\mu$, as claimed.

Next, by the definition of $f_{n}$ and the fact that it takes values
$0,1$ $P$-a.s., for $P$-typical $\mu$ the values $f_{n}(\mu)$
agree for all $n$. Hence for such $\mu$, either (i) $T_{x}\mu\in\mathcal{V}$
for $\mu$-a.e. $x$, or (ii) $T_{x}\mu\notin\mathcal{V}$ for $\mu$-a.e.
$x$. 

Using quasi-Palm property of $P$ one more time, the behavior of a
$P$-typical $\mu$ at $0$ is the same as its behavior at $\mu$-typical
points. Therefore with $P$-probability $1$, if $\mu\in\mathcal{V}$
then $T_{x}\mu\in\mathcal{V}$ for $\mu$-a.e. $x$, and if $\mu\notin\mathcal{V}$
then $T_{x}\mu\notin\mathcal{V}$ for $\mu$-a.e. $x$. 

It now follows that $\int_{\mathcal{V}}\left\langle \nu\right\rangle _{U}^{*}\, dP(\nu)$
is supported on $\mathcal{V}$, as desired.
\end{proof}
It is natural to ask whether an ``ergodic decomposition'' exists
for scaling measures on the spatial level. Two simple results of this
type are the following.
\begin{prop}
\label{pro:invariance of USMs-under-restriction}If $\mu$ is a scaling
measure, and $\mu(A)>0$, then $\mu_{A}$ is a scaling measure and
$\mu$ and $\mu_{A}$ generate the same distribution at $\mu_{A}$-a.e.
$x$.\end{prop}
\begin{proof}
This is immediate from the Besicovitch density theorem (Theorem \ref{lem:Besicovitch-lemma}),
which  implies that for $\mu_{A}$-a.e. $x$, 
\[
S_{t}^{\sqr}T_{x}(\mu_{A})-S_{t}^{\sqr}T_{x}\mu\rightarrow0
\]
and hence the sceneries at $x$ for $\mu$ and $\mu_{A}$ are asymptotic. 
\end{proof}
More generally,
\begin{prop}
\label{pro:equivalent-measures-have-same-sceneries}If $\mu$ is a
scaling measure and $0\neq\nu\ll\mu$ then $\nu$ is a scaling measure
and for $\nu$-a.e. $x$ the distributions generated by $\mu$ and
$\nu$ at $x$ are the same.\end{prop}
\begin{proof}
The proof is similar to the one above. Let $\varphi=\frac{d\nu}{d\mu}$,
let $x$ be a $\nu$-typical point, and write $\varphi_{x,t}=\varphi|_{B_{\exp(-t)}(x)}$.
Then for $\varepsilon>0$, applying the Besicovitch density theorem
to the set $\mathcal{A}_{\varepsilon}=\varphi^{-1}(\varphi(x)-\varepsilon,\varphi(x)+\varepsilon)$,
we find that 
\[
\lim_{t\rightarrow\infty}\frac{\nu(y\in B_{\exp(-t)}(x)\,:\,\varphi(y)\in\mathcal{A}_{\varepsilon})}{\mu(B_{\exp(-t)}(x))}=1
\]
Hence $\lim_{t\rightarrow\infty}\left\Vert \mu_{x,t}-\nu_{x,t}\right\Vert \leq\varepsilon$
(where $\left\Vert \cdot\right\Vert $ is the total variation norm)
and hence, since $\varepsilon$ was arbitrary, $\nu_{x,t}^{\sqr}$
and $\mu_{x,t}^{\sqr}$ are asymptotic. 
\end{proof}
One consequence is that if $\mu$ is a scaling measure, $\mathcal{U}$
is a set of distributions and $A$, $B=\mathbb{R}^{d}\setminus A$
is the partition of $\mathbb{R}^{d}$ according to whether $P_{x}\in\mathcal{U}$
or $P_{x}\notin\mathcal{U}$, respectively, then the conditional measures
on the atoms of the partition behave as described above: $\mu_{A}$-a.e.
point generates a distribution in $\mathcal{U}$ and $\mu_{B}$-a.e.
point generates a distribution in $\mathcal{M}^{*}\setminus\mathcal{U}$.
It is now natural to ask if this phenomenon holds for the partition
of $\mathbb{R}^{d}$ according to the generated distributions, that
is, the partition induced by the map $x\mapsto P_{x}$. We may disintegrate
$\mu$ with respect to this partition, and one might expect that typically
the induced measure on the atom $\{x\,:\, P_{x}=P\}$ will itself
be a USM generating $P$. 

Easy examples show that in general this is false. Consider for instance
Lebesgue measure $\sigma$ on the circle $S^{1}\subseteq\mathbb{R}^{2}$.
Then for $x\in S^{1}$ the measure $\sigma$ generates at $x$ the
EFD $P_{x}=\delta_{\lambda(x)}$, where $\lambda(x)$ is Lebesgue
measure on the line through the origin in the direction of the tangent
to $S^{1}$ at $x$. Thus the partition according to $P_{x}$ is the
partition of $S^{1}$ into antipodal points, and the conditional measure
on each atom generates the same trivial EFD everywhere. 

One may also ask the dual question concerning sums and integrals of
SMs. Consider SMs $\mu,\nu$. Then $\mu+\nu$ is equivalent to the
sum of three mutually singular measures, $\mu+\nu\sim\mu'+\theta+\nu'$,
such that $\mu'\ll\mu$, $\nu'\ll\nu$, and $\theta$ absolutely continuous
with respect to both $\mu$ and $\nu$ (the measure $\theta$ is obtained
as the part of $\nu$ which is non-singular with respect to $\mu$,
and $\nu'=\nu-\theta$; then $\mu'$ is the part of $\mu$ singular
with respect to $\theta$). Using the proposition above, we find that
$\mu+\nu$ is a SM and at a typical point generates the distribution
which is generated at some point by $\mu$ or $\nu$ (or both).

However, passing from finite sums to integrals of SMs we lose this
behavior. Suppose $Q$ is a distribution on SMs. Then $\int\mu\, dQ(\mu)$
need not be an SM. Indeed, any probability measure $\theta\in\mathcal{M}$
can be written as $\theta=\int\delta_{x}\, d\theta(x)$. Each $\delta_{x}$
is a USM but $\theta$ need not be.

\subsection{\label{sub:CP-to-USM}From CP-distributions to FDs to USMs}

In this section we show that centerings of CP-distributions are FDs,
and typical measures for a FD are USMs. 

We begin with the second of these. As discussed in Section \ref{sub:Generic-points},
if $P$ is an ergodic distribution for a transformation of a compact
metric space, then $P$-a.e. point is generic for $P$. USMs are analog
of generic points, and in our setting the analogous statement is:
\begin{thm}
\label{thm:EFDs-to-USMs}If $P$ is an EFD then $P$-a.e. $\mu$ is
a USM and generates $P$.\end{thm}
\begin{proof}
We assume that $P$ is in its extended version. For $f\in C(\mathcal{P}(\mathcal{M}^{\sqr}))$
let $\widetilde{f}:\mathcal{M}^{*}\rightarrow\mathbb{R}$ denote the
map $\mu\mapsto f(\mu^{\sqr})$ (which is no longer continuous, but
is still measurable and bounded). Choose a countable norm-dense family
$\mathcal{F}\subseteq C(\mathcal{P}(\mathcal{M}^{\sqr}))$. Then by
the ergodic theorem, for $f\in\mathcal{F}$ and $P$-a.e. $\mu$ we
have
\begin{equation}
\lim_{T\rightarrow\infty}\frac{1}{T}\int_{0}^{T}f(S_{t}^{\sqr}\mu)\, dt=\int\widetilde{f}\, dP\label{eq:scenery-genericity-for-one-function}
\end{equation}
so, since $\mathcal{F}$ is countable, for $P$-a.e. $\mu$ this holds
simultaneously for all $f\in\mathcal{F}$. By the quasi-Palm property,
for every $R>0$ for $P$-a.e. $\nu$, the limit \eqref{eq:scenery-genericity-for-one-function}
holds for $\mu=T_{x}^{*}\nu$ for $\nu$-a.e. $x\in B_{R}$, so, choosing
a sequence $R\rightarrow\infty$, for $P$-a.e. $\nu$ it holds for
$\nu$-a.e. $x\in\mathbb{R}^{d}$. Finally, this implies that for
$P$-a.e. $\nu$ and $\nu$-a.e. $x$ the desired limit holds for
$\mu=T_{x}^{*}\nu$ and all $f\in C(\mathcal{P}(\mathcal{M}^{\sqr}))$,
because $\mathcal{F}$  is dense in $C(\mathcal{P}(\mathcal{M}^{\sqr}))$
and the set of functions satisfying \eqref{eq:scenery-genericity-for-one-function}
is norm closed. This is the desired result.
\end{proof}
Now we drop the ergodicity assumption and prove:
\begin{proof}
[Proof of Theorem \ref{thm:FD-to-USM}] Let $P$ be a FD and $P=\int P_{\nu}\, dP(\nu)$
\,its ergodic decomposition with respect to $S^{*}$. By Theorem
\ref{thm:EFD-ergodic-components}, $P_{\nu}$ is a FD for $P$-a.e.
$\nu$, and so by the theorem above, for $P$-a.e. $\nu$, $P_{\nu}$-a.e.
$\mu$ is a USM generating $P_{\nu}$. By standard properties of the
ergodic decomposition, $P_{\mu}=P_{\nu}$ for $P_{\nu}$-a.e. $\mu$.
Therefore we have shown that for $P$-a.e. $\nu$ and $P_{\nu}$-a.e.
$\mu$, the measure $\mu$ is a USM generating $P_{\mu}$. Since $P=\int P_{\nu}\, dP(\nu)$,
the same statement holds for $P$-a.e. $\mu$, which is what we wanted
to prove.
\end{proof}
Next let $Q$ be a base-$b$ CP-distribution and recall the definition
of the centering operation, Definition \ref{def:centering}. Since
it is immediate that $(\mu',x')=M_{b}^{*}(\mu,x)$ satisfies
\[
S_{\log b}^{*}(T_{x}\mu)=T_{x'}^{*}\mu'
\]
the map
\[
\cnt_{0}:(\mu,x)\mapsto T_{x}^{*}\mu
\]
is a factor map from the discrete-time measure preserving system $(\mathcal{M}^{*}\times B_{1},Q,M_{b}^{*})$
to the discrete-time  measure preserving system $(\mathcal{M}^{*},\cnt_{0}Q,S_{\log b}^{*})$,
and the continuous time system $(\mathcal{M}^{*},\cnt_{0}Q,S^{*})$
is thus a factor of the suspension of $(\mathcal{M}^{*}\times B_{1},Q,M_{b})$
with roof function of height $\log b$. 
\begin{thm}
\label{thm:centering-is-EFD}The centering of an extended CP-distribution
is a FD.\end{thm}
\begin{proof}
Let $Q$ be an extended CP-distribution with restricted version $P$.
Adopting the notation from the construction of $Q$ in Section \ref{sub:Restricted-and-extended-versions},
let $(\mu_{n},x_{n})_{n\in\mathbb{Z}}$ be the two-sided $\mathcal{M}^{\sqr}\times B_{1}$-valued
stationary process with marginal $P$ and $(\mu_{n+1},x_{n+1})=M_{b}^{\sqr}(\mu_{n},x_{n})$,
and let $\pi$ denote the map $(\mu_{n},x_{n})_{n\in\mathbb{Z}}\mapsto\nu$
constructed there, i.e. 
\begin{equation}
\nu=\lim_{n\rightarrow\infty}T_{\mathcal{D}_{b^{n}}(x_{-n})}^{*}\mu_{-n}\label{eq:CP-extension-compatability}
\end{equation}
so $\nu\sim Q$. 

Let $R=\cnt_{0}Q$. Let $B\subseteq\mathbb{R}^{d}$ be a bounded neighborhood
of $0$ and let 
\begin{equation}
R'=\int\left\langle \nu\right\rangle _{B}^{*}\, dR(\nu)=\int\left(\int_{B}\delta_{T_{x}^{*}\nu}\, d\nu(x)\right)\, dR(\nu)\label{eq:diffusion-of-centering}
\end{equation}
we must show that $R\sim R'$. 

Fix a measurable set $\mathcal{U}\subseteq\mathcal{M}^{*}$. Suppose
that $R(\mathcal{U})=0$. In order to show that $R'(\mathcal{U})=0$
it suffices by \eqref{eq:diffusion-of-centering} to show that 
\[
\int\int1_{\mathcal{U}}(T_{x}^{*}\nu)\, d\nu(x)\, dR(\nu)=0
\]
For this it is enough to show that for a.e. realization $(\mu_{n},x_{n})_{n\in\mathbb{Z}}$
of the process and $\nu$ satisfying \eqref{eq:CP-extension-compatability},
we have 
\[
\int1_{\mathcal{U}}(T_{x}^{*}\nu)\, d\nu(x)=0
\]
By \eqref{eq:CP-extension-compatability}, this will follow once we
establish that for a.e. realization of the process and every $n\in\mathbb{N}$,
\[
\int1_{\mathcal{U}}(T_{x}^{*}\nu)\, dT_{\mathcal{D}_{b^{n}}(x_{-n})}^{*}\mu_{-n}(x)=0
\]
Fixing $n$, by stationarity of the process $(\mu_{n},x_{n})_{n\in\mathbb{Z}}$,
and the fact that the map $(\mu_{n},x_{n})\rightarrow\nu$ intertwines
the shift operation and $S_{\log b}^{*}$, we have 
\begin{eqnarray*}
\int1_{\mathcal{U}}(T_{x}^{*}\nu)\, dT_{\mathcal{D}_{b^{n}}(x_{-n})}^{*}\mu_{-n}(x) & = & \int1_{\mathcal{U}}(S_{n\log b}^{*}T_{x}\nu)\, d\mu_{0}(x)\\
 & = & \int1_{S_{-n\log b}^{*}\mathcal{U}}(T_{x}\nu)\, d\mu_{0}(x)\\
 & = & R(S_{-n\log b}^{*}\mathcal{U})\\
 & = & R(\mathcal{U})\\
 & = & 0
\end{eqnarray*}
because $R$ is $S^{*}$-invariant. We have established the claim
a.s. for each $n$ and therefore a.s. for all $n\in\mathbb{N}$, as
desired.

Conversely, for any measure $\nu$, since $0\in B$ we have (recall
that $B_{1}$ is the unit ball) 
\[
\left\langle \nu\right\rangle _{B_{1}}^{*}\ll\int\left\langle \theta\right\rangle _{B}^{*}\, d\left\langle \nu\right\rangle _{B_{1}}^{*}(\theta)
\]
so since $R=\int\left\langle \nu\right\rangle _{B_{1}}^{*}\, dQ(\nu)$,
\[
R=\int\left\langle \nu\right\rangle _{B_{1}}^{*}\, dQ(\nu)\ll\int\int\left\langle \nu\right\rangle _{B}^{*}\, d\left\langle \nu\right\rangle _{B_{1}}^{*}\, dQ(\nu)=\int\left\langle \nu\right\rangle _{B}^{*}dR(\nu)=R'
\]

Finally, the fact that $\cnt Q=\frac{1}{\log b}\int_{0}^{\log b}S_{t}^{*}\cnt_{0}Q\, dt$
is quasi-Palm follows trivially from the same property for $R=\cnt_{0}Q$.
Therefore $\cnt Q$ is a FD.
\end{proof}

\section{\label{sec:Examples}Examples}

In this section we present a variety of examples of EFDs and USMs.
The proofs usually involve constructing a stationary process similar
to a CP-process and centering. Many of our examples have appeared
before in related contexts. For additional examples of similar constructions
see \cite{Graf95,BedfordFisherUrbanski02,HochmanShmerkin09}.

In Section \ref{sec:Examples} we give further examples produced by
more deliberate construction.

\subsection{\label{sub:CP-processes-from-stationary-processes}Example: CP-processes
arising from stationary processes}

Recall that a stationary process is a sequence of random variables
$(X_{n})_{n\in\mathbb{Z}}$ such that for every $n\in\mathbb{Z}$
and $k\geq0$ the $k$-tuples $(X_{0},\ldots,X_{k})$ and $(X_{n},\ldots,X_{n+k})$
have the same distribution. 

The following example appears in \cite{Furstenberg08}. Let process
$(Y_{n})_{n\in\mathbb{Z}}$ be a stationary process with values in
$\{0,\ldots,b-1\}$ and define a random point $x$ and measure $\mu$
by%
\footnote{If we had defined CP-processes on $[0,1)^{d}$ using $b$-adic partitions,
the definition of $x_{n}$ would simplify to $x=\sum_{i=1}^{\infty}b^{-k}Y_{k}$.%
} 
\begin{eqnarray}
x & = & -1+2\sum_{k=1}^{\infty}b^{-k}Y_{k}\nonumber \\
\mu & = & \mbox{distribution of }x\mbox{ given }(Y_{n})_{n\leq0}\label{eq:prediction-measure}
\end{eqnarray}
Then one may verify that the distribution of $(\mu,x)$ is a CP-distribution. 

With the same notation, another example is 
\begin{eqnarray}
x & = & -1+2\sum_{n=1}^{\infty}b^{-n}Y_{n}\nonumber \\
\mu & = & \delta_{x}\label{eq:deterministic-CP-process}
\end{eqnarray}
Notice that the first construction generally gives rise to a non-trivial
FD, while the second gives rise to the FD supported on the measure
$\delta_{0}\in\mathcal{P}(\mathcal{M})$. 

As a by-product of this we have:
\begin{prop}
\label{pro:EFD-has-many-CP-extensions}A given EFD may be the centering
of many distinct CP-processes.\end{prop}
\begin{proof}
The trivial EFD arises as the centering of the CP-distribution $(\delta_{x},x)$
constructed as above from a process $(Y_{n})$. These distributions
are different for different initial processes because, from the distribution
of $(\delta_{x},x)$ we can recover the process, by setting $Y_{n}$=
$n$-th base-$b$ digit of $\frac{1}{2}(x-1)$. 
\end{proof}
Here is a concrete example, let $b=3$ and consider the process $(Y_{n})$
in which each $Y_{n}$ is chosen independently, equal to $0$ or $2$
with probability $1/2$ each. This corresponds to the product measure
on $\{0,1,2\}^{\mathbb{Z}}$ of the measure $\frac{1}{2}\delta_{0}+\frac{1}{2}\delta_{2}$.
Because the $Y_{n}$ are independent, given $Y_{0},Y_{-1},Y_{-2},\ldots$
the distribution of each $Y_{k},k\geq1$ is the same as before, and
so $x=-1+2\sum3^{-i}Y_{i}$ is distributed according to the normalized
Hausdorff measure $\mu$ on an affine image of the standard Cantor
set:
\[
C=\{-1+2x\,:\, x\in[0,1]\mbox{ can be written in base }3\mbox{ using only the digits }0,2\}
\]
The CP-distribution we get is $\delta_{\mu}\times\mu$. Note that
the first component doesn't change when $M_{3}$ is applied.

\subsection{\label{sub:invariant-measures-on-unterval}Example: $\times m$-invariant
measures on the $[0,1]$}

Let $b\geq2$ be an integer and let $\mu$ be an ergodic probability
measure on $[0,1]$ which is invariant under the $b$-to-1 map $f_{b}:x\mapsto bx\bmod1$.
Such a measure can be identified with a stationary process which is
the sequence digits in the base-$b$ expansion of $x\sim\mu$. As
such, $\mu$ gives rise to a shift-invariant measure on $\{0,\ldots,b-1\}^{\mathbb{N}}$
and its natural extension may be realized as a shift-invariant measure
$\widetilde{\mu}$ on $\{0,\ldots,b-1\}^{\mathbb{Z}}$ which projects
to $\mu$ on the positive coordinates. It is then not hard to see
that $\mu=\int\nu_{\omega}d\widetilde{\mu}(\omega)$, where $\nu_{\omega}$
is the distribution of $\sum_{n=1}^{\infty}b^{-i}\omega_{i}$ given
$(\omega_{j})_{j\leq0}.$ As we saw above the distribution on $\nu_{\omega}$
for $\omega\sim\widetilde{\mu}$ is an ergodic restricted CP-distribution
$P$ (after we identify $[0,1]$ with $B_{1}=[-1,1]$). 
\begin{prop}
\label{pro:invariant-measures-are-USMs}Every $\mu$ as above is a
USM for an ergodic EFD $P$, supported on measures whose dimension
is $\frac{1}{\log b}h(\mu,f_{b})$, where $h(\cdot,\cdot)$ is the
Kolmogorov-Sinai entropy.
\end{prop}
For a detailed proof see \cite{Hochman2010}.

\subsection{\label{sub:Self-similar-measures}Example: Self-similar measures}

Let $\Lambda$ be a finite set and $\{f_{i}\}_{i\in\Lambda}$ a system
of contracting similarities of $\mathbb{R}^{d}$, i.e.
\[
f_{i}(x)=r_{i}U_{i}(x)+v_{i}
\]
for some $0<r_{i}<1$, $v_{i}\in\mathbb{R}^{d}$ and $U_{i}$ an orthogonal
transformation of $\mathbb{R}^{d}$. For $a=a_{1}\ldots a_{n}\in\Lambda^{n}$
write
\[
f_{a}=f_{a_{1}}\circ f_{a_{2}}\circ\ldots\circ f_{a_{n}}
\]
It is well known that for every $x\in\mathbb{R}^{d}$ and $a\in\Lambda^{\mathbb{N}}$
the sequence $f_{a_{1}\ldots a_{n}}(x)$ converges, as $n\rightarrow\infty$,
to a point $\varphi(a)\in\mathbb{R}^{d}$ which is independent of
$x$. Furthermore, $\varphi:\Lambda^{\mathbb{N}}\rightarrow\mathbb{R}^{d}$
is continuous, and its image $X$ is the unique non-empty compact
set with the property $X=\bigcup_{i\in\Lambda}f_{i}(X)$ (the attractor
of the IFS). A set $X$ which arises in this way is called a self
similar set.

The IFS $\{f_{i}\}_{i\in\Lambda}$ satisfies the strong separation
condition if the sets $f_{i}(X),i\in\Lambda$ are pairwise disjoint;
this implies that $\varphi$ is an injection. If strong separation
holds then, by replacing $\Lambda$ with $\Lambda^{k}$ for some large
$k$ and the system $\{f_{i}\}_{i\in\Lambda}$ with $\{f_{a}\}_{a\in\Lambda^{k}}$,
and then applying a compactness argument, one may assume that there
is an open set $A\subseteq\mathbb{R}^{d}$ such that $X\subseteq A$
and $f_{i}(A),i\in\Lambda$ are pairwise disjoint and contained in
$A$.

Let $\widetilde{\mu}$ be a product measure on $\Lambda^{\mathbb{N}}$
with marginal $(p_{i})_{i\in\Lambda}$. The image $\mu=\varphi\widetilde{\mu}$
on $X$ is called a self-similar measure. We claim that a self similar
measure is uniformly scaling, and furthermore that the associated
FD is supported on measures which on every bounded set are equivalent,
up to a similarity, to a restriction of $\mu$.

Let $G<GL_{n}$ denote the group generated by the orthogonal maps
$U_{i}$, $i\in\Lambda$ and let $\gamma$ denote Haar measure on
$G$. Consider the distribution $P$ on $\mathcal{M}\times\mathbb{R}^{d}\times G$
obtained by selecting $y\sim\mu$ and $U\sim\gamma$ independently,
and forming the triple $(U\mu,Uy,U)$. 

For $(\nu,y,V)\in\mathcal{M}\times\mathbb{R}^{d}\times G$ which satisfies
$V^{-1}\nu=\mu$ and $V^{-1}y\in\supp\mu$ define $M(\nu,y,U)=(\nu',y',V')\in\mathcal{M}\times\mathbb{R}^{d}\times G$
as follows. Let $i\in\Lambda$ be the unique index such that $V^{-1}y\in f_{i}A$.
Let $A_{i}=f_{i}A$. Now define 
\begin{eqnarray*}
\nu' & = & \frac{1}{p_{i}}Vf_{i}^{-1}((V^{-1}\nu)|_{A_{i}})\\
y' & = & Vf_{i}^{-1}V^{-1}y\\
V' & = & U_{i}V
\end{eqnarray*}
It is easy to verify that (i) $M$ is defined $P$-a.e., (ii) $P$
is invariant under $M$, (iii) $P$ is adapted in the sense that conditioned
on the first component $\nu,$the second component $y$ is distributed
according to $\nu$. Indeed (i) is trivial, (iii) is by construction,
and (ii) is a direct consequence of self-similarity of $\mu$ and
invariant of $\gamma$ under left multiplication by $U_{i}$.

One can now repeat almost verbatim the construction of a FD from a
CP-distribution to obtain a FD from $P$. One starts with a two-sided
stationary process $(\nu_{n},y_{n},V_{n})_{n\in\mathbb{Z}}$ whose
marginals are $P$ and such that $M(\nu_{n},y_{n},V_{n})=(\nu_{n+1},y_{n+1},V_{n+1})$.
Let $i_{n}\in\Lambda$ be the index such that $V_{n}^{-1}y_{n}\in f_{i_{n}}A$,
and define $g_{n}:\mathbb{R}^{d}\rightarrow\mathbb{R}^{d}$ by $g_{n}=V_{n}f_{i_{n}}^{-1}V_{n}^{-1}$.
Then if for $n\geq0$ we set
\[
\theta_{n}=\frac{1}{p_{i_{-1}}\cdot p_{i_{-2}}\cdot\ldots\cdot p_{i_{-n}}}g_{-1}\ldots g_{-n+1}g_{-n}\nu_{-n},
\]
then $\theta_{n}$, $n\leq0$ have the property that $\theta_{n}$
is the restriction of $\theta_{n+1}$ to an appropriate open set.
Hence the random measure $\theta=\lim_{n\rightarrow\infty}\theta_{n}$
is well defined. 

Form the random measure $\widetilde{\theta}=S_{r}^{*}T_{y_{0}}^{*}\theta$,
where $r\in[0,r_{i_{0}})$ is chosen uniformly. Then the distribution
$\widetilde{P}$ of $\widetilde{\theta}$ is a FD supported on measures
which, on any bounded set, are absolutely continuous with respect
to an image of $\mu$ by a similarity whose orthogonal part comes
from $G$. The distribution of $\widetilde{\theta}$ is $S^{*}$-invariant,
and is a factor of the suspension of measure preserving system $(\mathcal{M}^{*}\times\mathbb{R}^{d}\times G,P,M)$
by the function whose height at $(\nu,y,V)$ is $r_{i}$, where $i\in\Lambda$
is the index such that $V^{-1}y\in f_{i}A$. 

Finally, let us show that $\widetilde{P}$ is ergodic and $\mu$ generates
$\widetilde{P}$. This can be deduced from the ergodicity of $(\mathcal{M}^{*}\times\mathbb{R}^{d}\times G,P,M)$,
but we give an alternative argument. Let $\widetilde{P}'$ be an ergodic
component of $\widetilde{P}$, and note that $\widetilde{P}'$ arises
as the factor of the suspension as above of an ergodic component $P'$
of $P$. By Theorem \ref{thm:EFD-ergodic-components}, $\widetilde{P}'$
is an EFD, and by Theorem \ref{thm:FD-to-USM} a $\widetilde{P}'$-typical
measure $\widetilde{\theta}'$ is a USM generating $\widetilde{P}'$.
On the other hand $\widetilde{\theta}'$ contains as an absolutely
continuous component a scaled copy of $\mu$ rotated by an element
of $G$. Hence $\mu$ generates a rotated version of $\widetilde{P}'$
and in particular is a USM. We also find that all ergodic components
of $\widetilde{P}$ are rotations of each other by elements of $G$.
The fact that $\widetilde{P}$ is ergodic is completed by showing
that the EFD generated by $\mu$ is in fact invariant under such rotation;
this follows from the fact that $\mu$ contains as absolutely continuous
components measures which are (scaled copies of) rotations of $\mu$
by every composition $U_{i_{1}}U_{i_{2}}\cdot\ldots\cdot U_{i_{n}}$,
$i_{1},\ldots,i_{n}\in\Lambda$, and these compositions generate $G$.
We note that this argument is related to Theorem \ref{thm:properties-of-homogeneous-measures},
although if $G$ is infinite then $\mu$ is not homogeneous.

\subsection{\label{sub:A-random-fractal-example}Example: Random fractals}

As another demonstration of the versatility of EFDs we give an example
of a fairly general random fractal construction that gives rise to
an EFD. One can easily generalize this further. 

Let $(I_{n},J_{n},w_{n})_{n=1}^{\infty}$ be an ergodic stationary
process in which each pair $(I_{n},J_{n})$ consists of closed, disjoint%
\footnote{It suffices to assume that the interiors are disjoint, but for simplicity
we assume disjointness. %
} sub-intervals of $[-1,1]$ and $0<w_{n}<1$. We assume that 
\begin{equation}
\mathbb{E}(w_{1}\log|I_{1}|+(1-w_{1})\log|J_{1}|)<\infty\label{eq:random-fractal-integrability-condition}
\end{equation}
where $|I|$ is the length of $I$. 

Construct a random measure $\mu$ as follows: set $\mu(I_{1})=w_{1}$
and $\mu(J_{1})=1-w_{1}$, and continue recursively for each of the
sub-intervals $I_{1},J_{1}$ using $(I_{2},J_{2},w_{2})$, i.e. in
the next step we define $\mu(T_{I_{1}}^{-1}(I_{2}))=w_{1}w_{2}$,
$\mu(T_{I_{1}}^{-1}(J_{2}))=w_{1}(1-w_{2})$ and similarly $\mu(T_{J_{1}}^{-1}(I_{2}))=(1-w_{1})w_{2}$,
$\mu(T_{J_{1}}^{-1}(J_{2}))=(1-w_{1})(1-w_{2})$. As before, we use
the notation $T_{I}$ for the homothety mapping $I$ onto $[-1,1]$. 

Unlike in our previous examples, the measure $\mu$ and $x\in\supp\mu$
do not necessarily determine the sequence of intervals $W_{n}\in\{I_{n},J_{n}\}$
such that $x\in\bigcap W_{n}$. Thus, in defining our CP-like process
we will encode this information directly in the state space. Let $\mathcal{U}$
denote the set of closed intervals of $[-1,1]$, considered as a topological
space e.g. by identifying $[a,b]$ with the pair $(a,b)\in\mathbb{R}^{2}$.
Give $\mathcal{U}^{\mathbb{N}}$ the product topology and associated
Borel structure. Consider the distribution $P$ on triples $(\nu,y,(U_{n})_{n=1}^{\infty})\in\mathcal{M}^{\sqr}\times[-1,1]\times\mathcal{U}^{\mathbb{N}}$
obtained by choosing a realization of the process $(I_{n},J_{n},w_{n})$,
constructing $\mu$ as above, choosing $x\sim\mu$, and selecting
the sequence $U_{n}\in\{I_{n},J_{n}\}$ in such a way that $T_{I_{n}}T_{I_{n-1}}\ldots T_{I_{1}}x\in[-1,1]$,
which defines the choice uniquely. Let $M$ be the map given by 
\[
M(\nu,y,(U_{n})_{n\in\mathbb{N}})=(T_{U_{1}}^{\sqr}\nu\,,\, T_{U_{i}}y\,,\,(U_{n+1})_{n\in\mathbb{N}})
\]
which is defined whenever $y\in\supp\nu$ and $T_{I_{n}}T_{I_{n-1}}\ldots T_{I_{1}}x\in[-1,1]$
for all $n$, and hence is defined $P$-a.e. One then verifies that
$P$ is $M$-invariant. By definition, given $\nu$, the distribution
of $y$ is $\nu$.

One now proceeds as before to define an extended version of this distribution,
and a corresponding EFD by centering. The only caveat here is that
one does not introduce a discrete-time system first, because the amount
of magnification at each step if different. Instead, one constructs
a continuous-time distribution directly as the suspension of the process
by a function whose height at $(\nu,y,(U_{n})_{n\in\mathbb{N}})$
is $\log|U_{n}|$. One must ensure that the mean of this height function
is finite for the resulting distribution to be finite. This is the
reason for the integrability condition \eqref{eq:random-fractal-integrability-condition}
above, which ensures that the roof function is in $L^{1}$. One may
verify that when the integrability condition fails, $\mu$ has dimension
0 a.s.

\subsection{Invariant measures for nonlinear iterated function systems}

Suppose that $\mathcal{I}=\{f_{i}\}_{i\in\Lambda}$ is a finite family
of $C^{1+\varepsilon}$ contractions of $[0,1]$, and that $f_{i}[0,1]$
are intervals that are disjoint except possibly at the endpoints.
For simplicity we assume that the contractions are orientation-preserving.
The attractor $X$ of such a system and the symbolic coding $\Phi:\Lambda^{\mathbb{N}}\rightarrow X$
are defined as in the linear case. A measure on $\mu\in\mathcal{P}(X)$
is called invariant if it is the image under $\Phi$ of a shift-invariant
measure on $\Lambda^{\mathbb{N}}$. For simplicity we consider measures
$\mu$ that are images of Gibbs measures for Holder potentials on
$\Lambda^{\mathbb{N}}$. The main property we use is that for such
a measure $\mu$, $f_{i}\mu\sim\mu_{f_{i}[0,1]}$ for all $i\in\Lambda^{*}=\bigcup_{k=1}^{\infty}\Lambda^{k}$,
and the Radon-Nikodym derivative is bounded away from $0$ and $\infty$
uniformly in $i$. The measure, being invariant, extends to an invariant
measure on $\Lambda^{\mathbb{Z}}$ and one can see easily that the
conditional measures given $\omega\in\Lambda^{\mathbb{Z}_{\leq0}}$,
the conditional measure $\widetilde{\mu}_{\omega}$ on $\Lambda^{\mathbb{N}}$
satisfy that $\Phi\widetilde{\mu}_{\omega}\sim\mu$, and the derivative
is again bounded away from $0,\infty$.

To keep things short we shall borrow freely from the notation and
definitions in \cite[Section 11.2]{HochmanShmerkin09}. For $\omega=(\ldots,\omega_{1},\omega_{0})\in\Lambda^{\mathbb{Z}_{\leq0}}$
let $F_{\omega}$ denote the limit diffeomorphism (we have changed
the indexing relative to \cite{Hochman2010} so that the coordinates
of $\omega$ begin at $0$, not $-1$). Let $\mathcal{I}_{\omega}=\{f_{i,\omega}\}_{i\in\Lambda}$,
$f_{i,\omega}=F_{\omega}f_{i}F_{\omega}^{-1}$, be the conjugate IFS
and $\Phi_{\omega}=F_{\omega}\Phi$ the coding map of $\mathcal{I}_{\omega}$. 

Now let $\nu_{\omega}=\Phi_{\omega}\widetilde{\mu}_{\omega}$ and
consider the set of triples $(\nu_{\omega},\omega,y)$, where $\omega\in\Lambda^{-\mathbb{N}}$,
$y=\Omega_{\omega}\eta$ for $\eta\in\Lambda^{\mathbb{N}}$. For $y\in F_{\omega}f_{i}([0,1))$,
define the magnification 
\begin{eqnarray*}
M(\nu,\omega,y) & = & (T_{f_{i,\omega}}^{\sqr}\nu_{\omega},\omega i,T_{f_{i}\omega}y)\\
 & = & (\nu_{\omega i},\omega i,\Phi_{\omega i}(\sigma\eta))
\end{eqnarray*}
where $\sigma$ is the shift map on $\Lambda^{\mathbb{N}}$. Let $Q$
be the distribution in such triples where $(\omega,\eta)$ are chosen
according to $\widetilde{\mu}$. Then, arguing as in section \ref{sub:invariant-measures-on-unterval},
it is routine to verify that $M$ preserves $Q$, and that an appropriate
adaptedness condition holds, i.e. conditioned on $\nu_{\omega}$ (in
fact, conditioned on $\omega$), $y$ is distributed according to
$\nu_{\omega}$.

\subsection{Brownian motion}

As observed by Gavish in \cite{Gavish09}, in dimension $d\geq3$
the occupation measure of a Brownian path, normalized to have mass
$1$ on $B_{1}$, is a USM generating the distribution on the random
occupation measure; hence it is an EFD.

\section{\label{sec:Equivalence-of-the}Equivalence of the models}

So far we have shown that (a) the centering of a CP-distribution is
a FD, and (b) a.e. measure of a FD is a USM. In this section we prove
two more such connections, which together imply Theorem \ref{thm:FD-to-CP}
and Proposition \ref{pro:CP-distribution-in-arbitrary-coords}. Recall
that $\lambda$ denotes Lebesgue measure.
\begin{thm}
\label{thm:USM-to-FD-and-CP}Let $b\geq2$ and $\mu\in\mathcal{M}$.
For $\mu$-a.e. $x$, every accumulation point $P$ of $\left\langle \mu\right\rangle _{x,T}^{\sqr}$
as $T\rightarrow\infty$ is the centering of a base-$b$ CP-distribution
$Q$, and in particular is a FD. Furthermore $Q$ may be chosen so
that $\int\theta^{\sqr}\, dQ(\theta)=\lambda^{\sqr}$.
\end{thm}

\begin{thm}
\label{thm:FDs-come-from-CPs}Given $b\geq2$, every EFD $P$ is the
centering of a base-$b$ CP-distribution $Q$ which can be chosen
so that $\int\theta^{\sqr}\, dQ(\theta)=\lambda^{\sqr}$.
\end{thm}
The second theorem above follows from the first by taking $\mu$ to
be a $P$-typical measure, since, as we saw, this is a USM generating
$P$. As for the part of Theorem \ref{thm:FD-to-CP} asserting that
every ergodic component of $P$ with respect to $S_{\log b}^{*}$
is the discrete centering of a base-$b$ CP-distribution, this is
evident from the fact that, if $Q$ is a CP-distribution whose continuous
centering is $P$, then the ergodic components of $Q$ have discrete
centerings that are $S_{\log b}^{*}$-invariant, ergodic (being factors
of ergodic CP-distributions) and integrate to $P$, hence they are
the ergodic components of $P$. 

It remains to prove the first theorem above. Throughout this section
we fix an integer $b\geq2$ and for brevity write 
\begin{eqnarray*}
\widetilde{M} & = & M_{b}^{\sqr}\\
\widetilde{S}_{t} & = & S_{t\log b}^{\sqr}
\end{eqnarray*}
Note that we have changed the time scale for $S_{t}$ so that $S_{1}$
and $M$ scale by the same factor $b$. We similarly modify the definition
of $\left\langle \mu\right\rangle _{x,T}^{\sqr}$ using $\widetilde{S}_{t}$:
\[
\left\langle \mu\right\rangle _{x,T}^{\sqr}=\frac{1}{T}\int_{0}^{T}\delta_{\widetilde{S}_{t}(T_{x}\mu)}\, dt=\frac{1}{T\log b}\int_{0}^{T\log b}\delta_{S_{t}^{\sqr}(T_{x}\mu)}\, dt
\]
We introduce analogous notation for CP-sceneries:
\begin{equation}
\left\langle \mu,x\right\rangle _{N}=\frac{1}{N}\sum_{n=1}^{N}\delta_{\widetilde{M}^{n}(\mu,x)}\in\mathcal{P}(\mathcal{M}^{\square}\times B_{1})\label{eq:CP-orbit-distribution}
\end{equation}
These operations are distinguished notationally by the number of arguments
inside the brackets, which correspond to the space on which the resulting
distribution is defined: $\left\langle \mu\right\rangle _{x,T}^{\sqr}$
and $\left\langle \mu,x\right\rangle _{N}$ are distributions on $\mathcal{M}^{\sqr}$
and on $\mathcal{M}^{\sqr}\times B_{1}$, respectively.

\subsection{\label{sub:Outline-of-equivalence}Outline of the argument}

Let $\mu\in\mathcal{M}^{\sqr}$ and fix a $\mu$-typical $x$. Let
$P$ be an accumulation point of the scenery distributions at $x$,
i.e. for some sequence $N_{k}\rightarrow\infty$,
\[
P=\lim_{k\rightarrow\infty}\left\langle \mu\right\rangle _{x,N_{k}}^{\sqr}
\]
We may assume without loss of generality that $N_{k}$ are integers
since $\left\langle \mu\right\rangle _{x,N_{k}}^{\sqr}$ and $\left\langle \mu\right\rangle _{x,[N_{k}]}^{\sqr}$
are asymptotic in $\mathcal{P}(\mathcal{M}^{\sqr})$. Our strategy
is to construct a CP-distribution whose centering is $P$. A first
attempt, which fails only narrowly, is the most direct approach: pass
to a subsequence $N_{k(i)}$ such that $\left\langle \mu,x\right\rangle _{N_{k(i)}}\rightarrow Q$
for a distribution $Q$ on $\mathcal{M}^{\sqr}\times B_{1}$, and
show that (i) $Q$ is a CP-distribution, (ii) the extended version
of $P$ is the centering of the extended version of $Q$.

The argument is complicated by the fact that the transformation $\widetilde{M}$
is not continuous, and by some related issues. For (ii), note that
if it were not for the restriction to $B_{1}$ implicit in the definition
of $\widetilde{M}$, we would have the identity 
\[
\cnt\left\langle \mu,x\right\rangle _{N_{k(i)}}=\left\langle \mu\right\rangle _{x,N_{k(i)}}^{\sqr}
\]
Indeed, if $(\mu',x')=\widetilde{M}^{n}(\mu,x)$ and $\mu''=T_{x'}^{*}\mu'$
then $\mu''=(S_{n\log b}^{*}T_{x}\mu)|_{B}$, where $B$ is a translate
of $B_{1}$ and the restriction occurs because of the restriction
in the definition of $\widetilde{M}$. Thus the distribution on the
left hand side of the above equation is supported on measures which
are restrictions of the measures from the distribution on the right
hand side. However, after taking the limiting distributions and passing
to extended versions this difference should disappear. As for (i),
there are two parts: $\widetilde{M}$-invariance would follow easily
if $\widetilde{M}$ were continuous (this is similar to the argument
in \ref{sub:Generic-points}). If we could do this, it would then
remain to show that $Q$ is adapted, i.e. that conditioned on the
first component $\nu$, the distribution of the second component is
$\nu$.

It turns out to be more efficient to prove adaptedness first. Recall
that a distribution $Q$ on $\mathcal{M}^{\sqr}\times B_{1}$ is adapted
if, given that the first component of a $Q$-realization is $\nu$,
the second component is distributed according to $\nu$. We shall
use the following characterization. We refer the reader to Section
\ref{sub:Measure-valued-integration} for a discussion of measure-valued
integration. 
\begin{lem}
$Q\in\mathcal{P}(\mathcal{M}^{\sqr}\times B_{1})$ is adapted if and
only if for every $f\in C(\mathcal{M}^{\sqr})$,\textup{
\begin{equation}
\int f(\nu)\cdot\nu\, dQ(\nu,x)=\int f(\nu)\cdot\delta_{x}\, dQ(\nu,x)\label{eq:1}
\end{equation}
}\end{lem}
\begin{proof}
If $Q$ is adapted then for $Q$-a.e. $\nu$ we have $\nu=\int\delta_{x}d\nu(x)$,
hence $f(\nu)\cdot\nu=\int f(\nu)\cdot\delta_{x}d\nu(x)$ and integrating
over $\nu$ and using adaptedness gives \eqref{eq:1}. In the other
direction suppose that \eqref{eq:1} holds for all $f\in C(\mathcal{M}^{\sqr})$.
A standard approximation argument shows that then for every $E\subseteq\mathcal{M}^{\sqr}$
we have
\begin{equation}
\int1_{E}(\nu)\, dQ(\nu,x)=\int1_{E}(\nu)\cdot\delta_{x}\, dQ(\nu,x)\label{eq:2}
\end{equation}
Choose a sequence of finite partitions $\mathcal{E}_{n}$ of $\mathcal{M}^{\sqr}$
refining to the partition into points (we can do so because $\mathcal{M}^{\sqr}$
is compact metric), and let $\widetilde{\mathcal{E}}_{n}=\{E\times B_{1}\,:\, E\in\mathcal{E}_{n}\}$
be the induced partition on $\mathcal{M}^{\sqr}\times B_{1}$. When
convenient we use the same notation for a partition and the smallest
$\sigma$-algebra containing it. Let $\theta_{1},\theta_{2}:\mathcal{M}^{\sqr}\times B_{1}\rightarrow\mathcal{M}^{\sqr}$
denote the measure-valued random variables $\theta_{1}(\nu,x)=\nu$
and $\theta_{2}(\nu,x)=\delta_{x}$, where $(\nu,x)\sim Q$. Then
\eqref{eq:2} means that $\mathbb{E}(\theta_{1}|\widetilde{\mathcal{E}}_{n})=\mathbb{E}(\theta_{2}|\widetilde{\mathcal{E}}_{n})$.
By the martingale convergence theorem,  $\mathbb{E}(\theta_{1}|\widetilde{\mathcal{E}}_{n})\rightarrow\theta_{1}$
a.s., and therefore also $\mathbb{E}(\theta_{2}|\widetilde{\mathcal{E}}_{n})\rightarrow\theta_{1}$
a.s. But by the martingale theorem we also have $\mathbb{E}(\theta_{2}|\widetilde{\mathcal{E}}_{n})\rightarrow\mathbb{E}(\theta_{2}|\widetilde{\mathcal{E}})$,
where $\widetilde{\mathcal{E}}=\{E\times B_{1}\,:\, E\subseteq\mathcal{M}^{\sqr}\mbox{ is Borel}\}$.
Thus $\mathbb{E}(\theta_{2}|\widetilde{\mathcal{E}})=\theta_{1}$,
which is the same as adaptedness.\end{proof}
\begin{prop}
\label{pro:adaptedness-of-typical-sceneries}Let $\nu\in\mathcal{M}^{\sqr}$.
Then for $\nu$-a.e. $y$ and any $f\in C(\mathcal{M}^{\sqr})$, writing
$Q_{N}=\left\langle \nu,y\right\rangle _{N}$,
\[
\lim_{N\rightarrow\infty}\left(\int f(\theta)\cdot(\theta-\delta_{z})\, dQ_{N}(\theta,z)\right)=0
\]
in the weak-{*} topology on $\mathcal{M}^{\sqr}$. In particular,
any accumulation point of $\left\langle \nu,y\right\rangle _{N}$
as $N\rightarrow\infty$ is adapted.\end{prop}
\begin{proof}
For $\theta\in\mathcal{M}^{\sqr}$ let the level-$m$ discretization
of $\theta$ be 
\[
\theta^{(m)}=\sum_{k\in\mathbb{Z}^{d}}\theta(\mathcal{D}_{b^{m}}(\frac{k}{b^{m}}))\delta_{k/b^{m}}
\]
Thus $\theta^{(m)}\rightarrow\theta$ uniformly in $\theta$, as $m\rightarrow\infty$.
It therefore suffices to prove that for all $m$ and $f\in C(\mathcal{M}^{\sqr})$,
for $\nu$-a.e. $y$ and $Q_{N}=\left\langle \nu,y\right\rangle _{N}$,
\[
\lim_{N\rightarrow\infty}\left(\int f(\theta^{(m)})\cdot(\theta-\delta_{z})^{(m)}\, dQ_{N}(\theta,z)\right)=0
\]
Since $C(\mathcal{M}^{\sqr})$ is separable, it suffices to prove
this for a countable dense set of $f\in C(\mathcal{M}^{\sqr})$. Thus,
expanding the definition of $Q_{N}$, we must show that for a fixed
$f\in C(X)$, for every $m$ and $\nu$-a.e. $y$,

\begin{equation}
\lim_{N\rightarrow\infty}\left(\frac{1}{N}\sum_{n=1}^{N}f(\nu_{n}^{(m)})\nu_{n}^{(m)}-f(\nu_{n}^{(m)})\delta_{y_{n}}^{(m)}\right)=0\label{eq:asymptotic-adaptedness}
\end{equation}
where $(\nu_{n},y_{n})=\widetilde{M}^{n}(\nu,y)$ (note that $\nu_{n}$
depends implicitly on $y$).

Fix $m$ and define measure-valued functions $F_{n},G_{n}:B_{1}\rightarrow\mathcal{M}^{\sqr}$
by
\[
F_{n}(y)=f(\nu_{n}^{(m)})\cdot\nu_{n}^{(m)}
\]
and
\[
G_{n}(y)=f(\nu_{n}^{(m)})\cdot\delta_{y_{n}}^{(m)}
\]
Now, proving \eqref{eq:asymptotic-adaptedness} is equivalent to showing
that for $\nu$-a.e. $y$,
\[
\lim_{N\rightarrow\infty}\frac{1}{N}\sum_{n=1}^{N}\left(F_{n}(y)-G_{n}(y)\right)=0
\]
Evidently (conditioning with respect to $\nu$), 
\[
\mathbb{E}(G_{n}|\mathcal{D}_{b^{n}})=F_{n}.
\]
Also, $y\mapsto\delta_{y_{n}}^{(m)}$ is $\mathcal{D}_{b^{n+m}}$-measurable.
Thus, for $0\leq i<m$, the averages
\[
\lim_{N\rightarrow\infty}\frac{1}{N}\sum_{n=1}^{N}\left(F_{nm+i}(y)-G_{nm+i}(y)\right)
\]
are averages of (measure-values) martingale differences, and converge
a.e. to $0$ by the law of large numbers for Martingale differences
\cite[Chapter 9, Theorem 3]{Feller71}. Hence their sums from $i=0$
to $m-1$ also converge to $0$, which is what we wanted to prove. 
\end{proof}
We return to Theorem \ref{thm:USM-to-FD-and-CP}. Suppose we begin
with a $\mu$-typical point $x$, an accumulation point $P=\lim_{k\rightarrow\infty}\left\langle \mu\right\rangle _{x,N_{k}}^{\sqr}$,
and, passing to a further subsequence, choose an accumulation points
$Q=\lim_{i\rightarrow\infty}\left\langle \mu,x\right\rangle _{N_{k(i)}}$.
We have just seen that $Q$ is adapted and we would like to show that
$Q$ is $M^{\sqr}$-invariant, and that the centering of the extended
version of $Q$ is the extended version of $P$. We shall not quite
prove these statements for $\mu$, but instead perturb $\mu$ by a
random translation and prove them for the perturbed measure, replacing
$Q$ above with an accumulation point of the $b$-sceneries of the
perturbed measure. We may first assume, without loss of generality,
that $\mu$ is supported on $B_{1/2}(0)$, since it suffices to prove
the result for an affine image of $\mu$, and we may apply the map
$x\mapsto\frac{1}{2}x$. 
\begin{prop}
\label{pro:USM-to-FD-and-CP}Let $\mu$ be a probability measure on
$B_{1/2}$. Fix a $\lambda$-typical $y\in B_{1/2}$ and $\mu$-typical
$x$ and let $\mu'=T_{y}\mu$ and $x'=x+y$. Suppose $N_{k}\rightarrow\infty$
and that $\left\langle \mu\right\rangle _{x,N_{k}}^{\sqr}\rightarrow P$
and $\left\langle \mu',x'\right\rangle _{N_{k}}\rightarrow Q$. Then
a.s. (with respect to the choice of $y$ and $x$),
\begin{enumerate}
\item $\int\theta\, dQ(\theta)=\lambda^{\sqr}$
\item \label{pro:invariance-of-CP-accumulation-pts}$Q$ is $\widetilde{M}$-invariant.
\item \label{pro:CP-acc-pts-center-to-S-acc-pts} $P=(\cnt\widetilde{Q})^{\sqr}$,
where $\widetilde{Q}$ is the extended versions of $Q$.
\end{enumerate}
In particular, with probability $1$, $Q$ is a CP-distribution (and
$P$ is a FD).
\end{prop}
The final conclusion of the proposition follows by combining (2) and
(3) above with the previous proposition, and implies Theorem \ref{thm:USM-to-FD-and-CP}.
Indeed, let $y\in B_{1/2}$ be a $\lambda$-typical point. Fix a set
$A$ of full $\mu$-measure so that for $x\in A$ the last proposition
holds for $\mu'=T_{y}\mu$ and $x'=x+y$. Fix $x\in A$ and suppose
that $P=\lim_{M_{k}\rightarrow\infty}\left\langle \mu\right\rangle _{x,M_{k}}$
for some sequence $M_{k}\rightarrow\infty$. Pass to a subsequence
$N_{k}$ of $M_{k}$ such that $\left\langle \mu',x'\right\rangle _{N_{k}}\rightarrow Q$,
and note that $\left\langle \mu'\right\rangle _{x',N_{k}}^{\sqr}\rightarrow P$.
So by the proposition $Q$ is a CP-distribution with the desired properties
and $P$ is a centering of the extended version of $Q$.

We prove (1)--(3) in the sections below. Throughout the proof let
$y,x,\mu',x'$, $N_{i}$, $Q$ and $P$ be as in the proposition.
Write 
\begin{equation}
Q_{i}=\left\langle \mu',x'\right\rangle _{N_{i}}\label{eq:3}
\end{equation}
so $Q_{i}\rightarrow Q$. The statements that follow hold a.s. for
the choice of $x,y$.

\subsection{Proof of Proposition \ref{pro:USM-to-FD-and-CP} (1)}

Note that the distributions of the second component of $Q_{i}$ and
$Q$ are determined completely by $x'$ and the sequence $(N_{i})_{i=1}^{\infty}$,
and do not depend on $\mu'$. Since $x'=x+y$ where $y$ is a $\lambda$-typical
point, this is the same as the distribution of the second coordinate
of $\left\langle \lambda_{B_{1}},z\right\rangle _{N_{i}}$ for a $\lambda_{B_{1}+x}$-typical
point $z$, and by Proposition \ref{pro:adaptedness-of-typical-sceneries}
or the law of large numbers, these distributions converge to $\lambda_{B_{1}}$.
Therefore we have
\[
\int\delta_{z}\, dQ(\theta,z)=\lambda_{B_{1}}
\]
But by adaptedness of $Q$, we also have
\[
\int\delta_{z}\, dQ(\theta,z)=\int\theta\, dQ(\theta,z)
\]
and (1) follows.

\subsection{Proof of Proposition \ref{pro:USM-to-FD-and-CP} (2)}

Let 
\begin{equation}
E=\bigcup_{D\in\mathcal{D}_{b}}\partial D\label{eq:b-adic-boundaries}
\end{equation}
where $\partial$ denotes topological boundary, and define 
\begin{eqnarray*}
\mathcal{N}_{1} & = & \{(\theta,z)\in\mathcal{M}^{\sqr}\times B_{1}\,:\,\theta(E)=0\}\\
\mathcal{N}_{2} & = & \{(\theta,z)\in\mathcal{M}^{\sqr}\times B_{1}\,:\,\theta(\mathcal{D}_{b}(z))>0\}
\end{eqnarray*}
and
\[
\mathcal{N}=\mathcal{N}_{1}\cap\mathcal{N}_{2}
\]
These are Borel sets. $\mathcal{N}_{2}$ is the domain of $M$, and
$\mathcal{N}_{1}$ is a superset of the discontinuity points of $M$.
\begin{lem}
\label{lem:CP-acc-pts-supported-on-N}$Q$ is supported on $\mathcal{N}$.\end{lem}
\begin{proof}
From Proposition \ref{pro:USM-to-FD-and-CP} (1) we know that $\int\theta(E)\, dQ(\theta)=\lambda(E)=0$,
hence $\theta(E)=0$ for $Q$-a.e. $\theta$, giving that $Q$ is
supported on $\mathcal{N}_{1}$.

To prove that $Q$ is supported on $\mathcal{N}_{2}$ we argue as
follows. First, we claim that for a.e. choice of $x,y$ in the definition
of $\mu',x'$,
\begin{equation}
\liminf_{n\rightarrow\infty}\left(-\frac{1}{n\log b}\log\mu'(\mathcal{D}_{b^{n}}(x'))\right)\leq d\label{eq:dimension-bound-on-scenery-decay}
\end{equation}
Indeed, given $y$ and $\mu'=T_{y}\mu$, for $\mu'$-a.e. $x'$ (equivalently,
$\mu$-a.e. $x$ and $x'=x+y$), this bound follows from the fact
that the Hausdorff dimension of $\mathbb{R}^{d}$ is $d$, and hence
any measure supported on it has upper pointwise dimension at most
$d$. See e.g. Section \ref{sub:Entropy-and-entropy-dim}.

Next note that, writing again $(\mu'_{n},x'_{n})=M^{n}(\mu',x')$
\[
\mu'(\mathcal{D}_{b^{n}}(x'))=\prod_{j=0}^{n-1}\frac{\mu'(\mathcal{D}_{b^{j+1}}(x'))}{\mu'(\mathcal{D}_{b^{j}}(x'))}=\prod_{j=1}^{n}\mu'_{j}(\mathcal{D}_{b}(x'_{j}))
\]
Taking logarithms and using \eqref{eq:dimension-bound-on-scenery-decay}
we see that
\[
\frac{1}{n\log b}\sum_{j=1}^{n}\left(-\log\mu'_{j}(\mathcal{D}_{b}(x'_{j}))\right)\leq d+1
\]
for all large enough $n$. For $n=N_{i}$ and $Q_{i}$ defined as
in \eqref{eq:3}, this can be written as
\[
\int\left(-\log\theta(\mathcal{D}_{b}(z))\right)\, dQ_{i}(\theta,z)\leq(d+1)\log b
\]
Since the integrand is non-negative, we conclude e.g. by Chebychev
that for every $r>0$,
\[
Q_{i}\left\{ (\theta,z)\,:\,\theta(\mathcal{D}_{b}(z))<e^{-r}\right\} \leq\frac{(d+1)\log b}{r}
\]
From this we conclude (by approximating $\theta(\mathcal{D}_{b}(z))$
by a continuous function of $\theta$) that a similar bound holds
for $Q=\lim Q_{i}$, and hence 
\[
Q\left\{ (\theta,z)\,:\,\theta(\mathcal{D}_{b}(z))>0\right\} =1
\]
so $Q$ is supported on $\mathcal{N}_{2}$. 
\end{proof}
Recall that, writing \eqref{eq:3} explicitly, $Q_{i}=\frac{1}{N_{i}}\sum_{j=1}^{N_{i}}\delta_{\widetilde{M}^{j}(\mu',x')}$.
Let $\widehat{Q}_{i}$ denote the distribution obtained by replacing
each summand $\delta_{(\theta,z)}$ in this average with $\delta_{(\widehat{\theta},z)}$
where $\widehat{\theta}$ is the measure
\[
\widehat{\theta}=\theta_{B_{1}\setminus E}
\]
Then  $\widehat{Q}_{i}$ is supported on $\mathcal{N}$. Clearly the
fact that $Q$ is supported on $\mathcal{N}_{1}$ implies that 
\[
\widehat{Q}_{i}-Q_{i}\rightarrow0
\]
Since $Q_{i}=\frac{1}{N_{i}}\sum_{j=1}^{N_{i}}\widetilde{M}^{j}\left\langle \mu',x'\right\rangle $
it is also clear that 
\[
\widetilde{M}Q_{i}-Q_{i}\rightarrow0
\]
Write
\[
\widetilde{M}\widehat{Q}_{i}-\widehat{Q}_{i}=(\widetilde{M}\widehat{Q}_{i}-\widetilde{M}Q_{i})+(\widetilde{M}Q_{i}-Q_{i})+(Q_{i}-\widehat{Q}_{i})
\]
We have seen that the last two terms on the right tend to $0$ as
$i\rightarrow\infty$. The first term, $\widetilde{M}\widehat{Q}_{i}-\widetilde{M}Q_{i}$,
also tends to $0$, again using Proposition \ref{pro:USM-to-FD-and-CP}
(1). We conclude that
\[
\widetilde{M}\widehat{Q}_{i}-\widehat{Q}_{i}\rightarrow0
\]
We have already shown that $\widehat{Q}_{i}\rightarrow Q$. Since
$\widehat{Q}_{i},Q$ are supported on $\mathcal{N}$, and $\widetilde{M}$
is continuous on $\mathcal{N}$, we conclude that $\widetilde{M}\widehat{Q}_{i}\rightarrow\widetilde{M}Q$,
so by the above 
\[
\widetilde{M}Q-Q=0
\]
as desired.

\subsection{Proof of Proposition \ref{pro:USM-to-FD-and-CP} (3)}

We continue with the previous notation. We wish to show that the extended
version of $P$ is the centering of the extended version $\widetilde{Q}$
of $Q$, or equivalently, that $(\cnt\widetilde{Q})^{\sqr}=P$. We
continue to write $\widetilde{M}=M_{b}^{\square}$ and $\widetilde{S}_{t}=S_{t\log b}^{\sqr}$.

For $m\in\mathbb{N}$ let 
\[
P_{m}=\widetilde{S}_{m}\cnt Q
\]
or more explicitly, the push-forward of $Q$ through the map
\begin{equation}
(\theta,z)\mapsto\int_{0}^{1}\delta_{\widetilde{S}_{m+t}(T_{z}\theta)}\, dt\label{eq:continuous-centering-integral}
\end{equation}
This map is continuous where it is defined. Indeed, the maps $\widetilde{S}_{t}=S_{t\log b}^{\sqr}$
are discontinuous only when there is mass on the boundary of $B_{e^{-t\log b}}$.
For a given $\theta$, positive mass on the boundary can occur only
for countable many $t$. Therefore if $(\theta_{n},z_{n})\rightarrow(\theta,z)$,
then for all but countable many $t$, $\widetilde{S}_{m+t}T_{z}\theta_{n}\rightarrow\widetilde{S}_{m+t}T_{z}\theta$.
By bounded convergence (applied after integrating against a continuous
test function), we have $\int_{0}^{1}\delta_{\widetilde{S}_{m+t}(T_{z}\theta_{n})}\, dt\rightarrow\int_{0}^{1}\delta_{\widetilde{S}_{m+t}(T_{z}\theta)}\, dt$. 
\begin{lem}
$P_{m}\rightarrow P$ as $m\rightarrow\infty$.
\end{lem}
From the construction of the extended version $\widetilde{Q}$ of
$Q$ and the centering operation, it is easy to see that $P_{m}\rightarrow(\cnt\widetilde{Q})^{\sqr}$
, so the lemma implies $P^{\sqr}=(\cnt Q)^{\sqr}$. and hence $P=\cnt Q$.
\begin{proof}
[Proof of the Lemma] Let $P_{m,i}=\widetilde{S}_{m}\cnt Q_{i}$.
Our work will be completed by showing that
\[
\lim_{i\rightarrow\infty}P_{m,i}=P_{m}
\]
and
\[
\lim_{m\rightarrow\infty}\lim_{i\rightarrow\infty}P_{m,i}=P
\]
Indeed, the first of these follows from the continuity if the map
\eqref{eq:continuous-centering-integral} discussed prior to the lemma.
As for the second statement, we make the following observation. With
$m$ fixed and $n\in\mathbb{N}$, write $(\mu'_{n},x'_{n})=\widetilde{M}^{n}(\mu',x')$,
and notice that as long as the distance of $x'_{n}$ from $E$ is
at least $b^{-m}$, we have 
\[
\widetilde{S}_{n+m}(T_{x'}\mu')=\widetilde{S}_{m}(T_{x'_{n}}\mu'_{n})
\]
(note that if the operators $\widetilde{S}_{t}$ were defined using
the $*$-variant instead of $\sqr$ this would be an identity without
any assumption on $x'_{n}$). Thus under these assumptions,
\[
\int_{0}^{1}\delta_{\widetilde{S}_{n+m+t}(T_{x'}\mu')}\, dt=\widetilde{S}_{m}\cnt\delta_{(\mu'_{n},x'_{n})}
\]
Now, we can write 
\begin{eqnarray*}
P_{m,i}-\left\langle \mu'\right\rangle _{x',N_{i}} & = & \frac{1}{N_{i}}\sum_{n=1}^{N_{i}}\widetilde{S}_{m}\cnt\delta_{(\mu'_{n},x'_{n})}-\frac{1}{N_{i}}\sum_{n=0}^{N_{i}-1}\int_{0}^{1}\delta_{\widetilde{S}_{t}(T_{x'}\mu')}\, dt\\
 & = & \frac{1}{N_{i}}\sum_{n=1}^{N_{i}-m}\left(\widetilde{S}_{m}\cnt\delta_{(\mu'_{n},x'_{n})}-\int_{0}^{1}\delta_{\widetilde{S}_{n+m+t}(T_{x'}\mu')}dt\right)\\
 &  & +\frac{m}{N_{i}}\theta_{m,i}
\end{eqnarray*}
where $\theta_{m,i}$ is a probability measure, so that with $m$
fixed, we have $\frac{m}{N_{i}}\theta_{m,i}\rightarrow0$ as $i\rightarrow\infty$.
On the other hand, as we have seen, in the average above the summands
vanish whenever the distance of $x'_{n}$ from $E$ is at least $b^{-m}$.
Therefore, as $i\rightarrow\infty$ the right hand side is a probability
measure whose total mass is asymptotic to
\begin{equation}
\frac{1}{N_{i}}\#\{1\leq n\leq N_{i}\,:\, d(x_{n},E)>b^{-m}\}\label{eq:frequency-of-CP-scenery-approaching-bdry}
\end{equation}
Using again the fact that $\frac{1}{N_{i}}\sum_{n=1}^{N_{i}}\delta_{x_{n}}\rightarrow\lambda_{B_{1}}$
as $i\rightarrow\infty$ (Proposition \ref{pro:USM-to-FD-and-CP}
(1)), we see that \eqref{eq:frequency-of-CP-scenery-approaching-bdry}
converges as $i\rightarrow\infty$ to $1-c_{m}$ for a sequence $c_{m}$
satisfying $\lim_{m\rightarrow\infty}c_{m}=0$. This completes the
proof. 
\end{proof}

\subsection{\label{sub:FD-to-CP-and-erg-decomposition}CP-distribution from FDs}
\begin{lem}
If $(\Omega,\mathcal{B})$ and $(\Omega',\mathcal{B}')$ are standard
Borel spaces and $\pi:\Omega\rightarrow\Omega'$ an onto, measurable
map, and if $\mu'$ is a probability measure on $(\Omega',\mathcal{B}')$,
then there is a measure $\mu$ on $\Omega$ such that $\pi\mu=\mu'$.\end{lem}
\begin{proof}
See \cite[Lemma 2.2]{Varadarajan1963}.
\end{proof}

\begin{proof}
[Proof of Theorem \ref{thm:FDs-come-from-CPs}] We show that if $P$
is a FD then there is a CP-process $Q$ with $\cnt Q=P$. If $P$
is ergodic choose a $P$-typical $\mu$. By Theorem \ref{thm:FD-to-USM}
we know that $\mu$ generates $P$, so  by Proposition \ref{pro:USM-to-FD-and-CP}
\eqref{pro:CP-acc-pts-center-to-S-acc-pts}, $P$ is the centering
of a CP-distribution with the desired properties.

In the non-ergodic case let $\mathcal{E}\mathcal{CP}$ denote the
set of ergodic CP-distributions and $\mathcal{E}\mathcal{FD}$ the
set of EFDs. One may verify that both sets are measurable subsets
of the Borel spaces in which they reside, and therefore the restriction
of the $\sigma$-algebra to them is again a standard Borel space.
The map $\cnt:\mathcal{ECP\rightarrow\mathcal{EFD}}$ is measurable
and the paragraph above shows that it is onto. Thus by the preceding
lemma, any probability measure $\tau$ on $\mathcal{EFD}$ lifts to
a probability measure $\tau'$ on $\mathcal{ECP}$ with $\cnt\tau'=\tau$.
If $P$ is a FD we may identify it with a probability measure $\tau\in\mathcal{P}(\mathcal{EFD})$
corresponding to its ergodic decomposition. Let $\tau'\in\mathcal{P}(\mathcal{ECP})$
be the lift of $\tau$. Then $Q=\int R\, d\tau'(R)$ is a CP-distribution
and $\cnt Q=P$.

In order to obtain the additional property of $Q$ replace $\mathcal{EFD}$
in the argument above with the set of ergodic CP-distributions $Q$
with $\int\theta\, dQ(\theta)=\lambda^{\sqr}$. We have seen that
$Q\mapsto\cnt Q$ is still onto $\mathcal{EFD}$. Continue as before.
\end{proof}

\subsection{\label{sub:Changing-coordinates}Changing coordinates}

Let us quickly prove Corollary \ref{cor:CP-change-of-coords}. Let
$Q$ be a FD. Let $L\in GL_{n}(\mathbb{R})$ be a linear map and consider
the new coordinates system whose basis is $\{Le_{i}\}$, where $\{e_{i}\}$
is the standard basis of $\mathbb{R}^{d}$. Viewed in the new coordinate
system, $Q$ is the same as $(L^{-1})^{*}Q$ in the standard coordinate
system, and by Proposition \ref{prop:linear-images-of-FDs} this is
an FD. Therefore, given $b$, by Theorem \ref{thm:FDs-come-from-CPs},
there is a CP-distribution $Q'$ (in standard coordinates) whose centering
is $(L^{-1})^{*}Q$. Then $LQ'$, defined in the obvious way, is a
CP-distribution in the new coordinate system and its centering is
$L^{*}(L^{-1})^{*}Q=Q$. Now, since centering involves only pushing
a distribution forward through a map that translates, scales and normalizes,
for any Borel set $\mathcal{A}\subseteq\mathcal{M}$ which is invariant
under these operations, we have $Q(\mathcal{A})=Q'(\mathcal{A})$
(more precisely, $Q'(\mathcal{A})$ denotes $(\pi Q')(\mathcal{A})$,
where $\pi$ is projection to the first component of $\mathcal{M}^{\sqr}([-1,1]^{d})$).

Finally, if instead of an FD $Q$ we began with a CP-distribution
$P$, we could carry out the above for $Q=\cnt P$. Since $Q(\mathcal{A})=P(\mathcal{A})$
we obtain again the desired result.

\section{\label{sec:Projections-of-USMs}Geometry of FDs and USMs}

\subsection{\label{sub:Entropy-and-entropy-dim}\label{sub:-padic-cells-and-filtrations}Dimension
and entropy}

Recall that $\mathcal{D}_{b}$ is the partition of $\mathbb{R}^{d}$
into cubes $\times_{i=1}^{d}I_{i}$ for intervals of the form $I_{i}=[\frac{k}{b},\frac{k+2}{b})$,
for $k\in\mathbb{Z}$ with $k=b\bmod2$. This partitions $B_{1}$
into $b^{d}$ homothetic cubes. For a partition $\mathcal{D}$ of
$E\subseteq\mathbb{R}^{d}$ and $x\in E$, we write $\mathcal{D}(x)$
for the unique partition element containing $x$. The proof of the
following can be found in \cite[Theorem 15.3]{Pesin97}: 
\begin{lem}
\label{lem:decay-along-filtration}Let $\mu$ be a measure on $\mathbb{R}^{d}$
and $b$ an integer. Then for $\mu$-a.e. $x$,
\begin{eqnarray}
\overline{D}_{\mu}(x) & = & \limsup_{b\rightarrow\infty}\frac{\log\mu(\mathcal{D}_{b}(x))}{\log(1/b)},\nonumber \\
\underline{D}_{\mu}(x) & = & \liminf_{b\rightarrow\infty}\frac{\log\mu(\mathcal{D}_{b}(x))}{\log(1/b)}.\label{eq:4}
\end{eqnarray}
In particular, $\dim\mu=\alpha$ if and only if 
\[
\lim_{n\rightarrow\infty}\frac{\log\mu(\mathcal{D}_{b}(x))}{\log(1/b)}=\alpha\qquad\mbox{for }\mu\mbox{-a.e. }x
\]

\end{lem}
If $\mu$ is a probability measure and $\mathcal{Q}$ is a finite
or countable partition, then the Shannon entropy of $\mathcal{Q}$
with respect to $\mu$ is 
\[
H(\mu,\mathcal{Q})=-\sum_{Q\in\mathcal{Q}}\mu(Q)\log\mu(Q),
\]
with the convention that $0\log0=0$. This quantity measures how spread
out $\mu$ is among the atoms of $Q$. For the basic properties of
Shannon entropy see \cite{CoverThomas06}.
\begin{lem}
\label{lem:lower-dim-lower-boundes-entropy}If $\ldim\mu\geq\beta$
then
\[
\liminf_{b\rightarrow\infty}\frac{1}{\log b}H(\mu,\mathcal{D}_{b})\geq\beta
\]
\end{lem}
\begin{proof}
Since $H(\mu,\mathcal{D}_{b})=\int\log\mu(\mathcal{D}_{b^{n}}(x))d\mu(x)$,
the claim follows by integrating \eqref{eq:4} and applying Fatou's
lemma. See also \cite{FanLauRao02}.
\end{proof}
With $b$ fixed, the function $H(\cdot,\mathcal{D}_{b})$ is discontinuous
on the space of probability measures on $\mathbb{R}$, but only mildly
so:
\begin{lem}
\label{lem:continuous-approximation-of-entropy}There is a constant
$C$ (depending only on the dimension $d$) such that, for every $b$,
there is a continuous function $f_{b}$ on the space of probability
measures on $\mathbb{R}^{d}$, such that for any probability measure
$\mu$ on $\mathbb{R}^{d}$,
\[
\left|f_{b}(\mu)-H(\mu,\mathcal{D}_{b})\right|<C
\]
\end{lem}
\begin{proof}
Choose a countable partition of unity $\varphi_{u}$, $u\in\mathbb{Z}^{d}$
such that each $\varphi_{u}$ is continuous and supported on a cube
$B_{1/b}(\frac{1}{b}u)$; it is not hard to see that such a partition
exists. We claim that
\[
f_{b}(\mu)=-\sum_{u\in\mathbb{Z}^{d}}(\int\varphi_{u}\, d\mu)\log(\int\varphi_{u}\, d\mu)
\]
has the desired properties.

Given $\mu$ define a probability measure $\nu$ on $\mathbb{R}^{d}\times\mathbb{Z}^{d}$
by
\[
\nu(A\times\{i\})=\int_{A}\varphi_{i}(x)\, d\mu(x)
\]
Let $\mathcal{F}_{1}$ be the partition of $\mathbb{R}^{d}\times\mathbb{Z}^{d}$
induced from the first coordinate by $\mathcal{D}_{b}$, that is $\mathcal{F}_{1}=\{I\times\mathbb{Z}^{d}\,:\, I\in\mathcal{D}_{b}\}$.
Then 
\[
H(\mu,\mathcal{D}_{b})=H(\nu,\mathcal{F}_{1})
\]
Also, let $\mathcal{F}_{2}$ be the partition of $\mathbb{R}^{d}\times\mathbb{Z}^{d}$
according to the second coordinate. Then 
\[
H(\nu,\mathcal{F}_{2})=f_{b}(\mu)
\]
Let $\mathcal{E}$ denote the partition $\mathcal{E}=\{I\times\{j\}\,:\, I\in\mathcal{D}_{b}\,,\, j\in\mathbb{Z}^{d}\}$.
Notice that $\mathcal{E}$ refines $\mathcal{F}_{1}$ and each $A\in\mathcal{F}_{1}$
contains at most $b^{d}$ atoms $A'\in\mathcal{E}$ of positive $\nu$-mass
(because each $I\in\mathcal{D}_{b}$ intersects at most $b^{d}$ balls
$B_{2/b}(\frac{1}{b}u)$, $u\in\mathbb{Z}^{d}$, and the $\varphi_{u}$
are supported on such balls). Therefore,
\[
H(\nu,\mathcal{F}_{1})\leq H(\nu,\mathcal{E})=H(\nu,\mathcal{F}_{1})+H(\nu,\mathcal{E}|\mathcal{F}_{1})\leq H(\nu,\mathcal{F}_{1})+d\log b
\]
Similarly $\mathcal{E}$ refines $\mathcal{F}_{2}$ and each atom
$A\in\mathcal{F}_{2}$ contains at most $b^{d}$-atoms of $\mathcal{E}$
of positive $\nu$-mass (since each $\varphi_{u}$ is supported on
at most $b^{d}$ cells $I\in\mathcal{D}_{b}$ ) , so
\[
H(\nu,\mathcal{F}_{2})\leq H(\nu,\mathcal{E})=H(\nu,\mathcal{F}_{2})+H(\nu,\mathcal{E}|\mathcal{F}_{2})\leq H(\nu,\mathcal{F}_{2})+d\log b
\]
Combining the last four equations we have
\[
\left|H(\mu,\mathcal{D}_{b})-f_{b}(\mu)\right|=\left|H(\nu,\mathcal{F}_{1})-H(\nu,\mathcal{F}_{2})\right|\leq d\log b^{2}
\]

\end{proof}

\subsection{\label{sub:Dimension-of-FDs}Dimension of FDs and USMs}
\begin{proof}
[Proof of Lemma \ref{lem:dimension-of-FDs}] Let $P$ be an EFD and
$0<r<1$. Set
\[
F(\mu)=\frac{\log\mu(B_{r}(0))}{\log r}
\]
and notice that
\begin{eqnarray}
\frac{\log\mu(B_{r^{N}})}{\log(r^{N})} & = & \frac{1}{N\log r}\sum_{n=1}^{N}\log\left(\mu(B_{r^{n}}(0))/\mu(B_{r^{n-1}}(0))\right)\nonumber \\
 & = & \frac{1}{N}\sum_{n=1}^{N}\frac{\log(S_{-n\log r}^{*}\mu(B_{r}(0)))}{\log r}\nonumber \\
 & = & \frac{1}{N}\sum_{N=1}^{n}F(S_{n\log(1/r)}^{*}\mu)\label{eq:5}
\end{eqnarray}
This is an ergodic average for the transformation $S_{\log(1/r)}^{*}$,
so the limit exists $P$-a.e. Although $S_{\log(1/r)}^{*}$ may not
be ergodic, it is easy to see directly that the limit is invariant
under $S_{t}^{*}$ for every $t$ (e.g. because it is the local dimension
at $0$, which is invariant under staling), so it is $P$-a.e. constant
and equal to $\alpha=\int\frac{\log\mu(B_{r}(0))}{\log r}\, dP(\mu)$.
Thus the local dimension at $0$ satisfies $D_{\mu}(0)=\alpha$ for
$P$-a.e. $\mu$ and by the quasi-Palm property, for $P$-a.e. $\mu$
we have $D_{\mu}(x)=\alpha$ for $\mu$-a.e. $x$, so $\dim\mu=\alpha$.
\end{proof}
Setting 
\[
\widetilde{F}(\mu)=\int_{0}^{-\log r}F(S_{t}^{*}\mu)\, dt
\]
it follows from the lemma above that $\dim P=\int\widetilde{F}\, dP$,
and the same formula holds for non-ergodic FDs (recall that the dimension
of a FD is the mean dimension of its measures, so the non-ergodic
case follows from the ergodic one). The advantage of $\widetilde{F}$
is that, as a function $\widetilde{F}:\mathcal{M}^{\sqr}\rightarrow\mathbb{R}$,
it is continuous. 
\begin{proof}
[Proof of Proposition \ref{pro:dimension-via-accumulation-pts}] Let
$\mu\in\mathcal{M}$. Recall that the set of accumulation points of
$\left\langle \mu\right\rangle _{x,t}^{\sqr}$ as $t\rightarrow\infty$
is denoted $\mathcal{V}_{x}$. 

Let us show for example that for $\mu$-a.e. $x$, 
\[
\underline{D}_{\mu}(x)\geq\essinf_{y\sim\mu}\inf_{Q\in\mathcal{V}_{y}}\dim Q
\]
Denote the right hand side by $\beta$. By Theorem \ref{thm:limit-distributions-are-FDs}
there is a set of full $\mu$-measure of points $x$ such that every
accumulation point of $\left\langle \mu\right\rangle _{x,T}^{\sqr}$
is a FD. For such a point $x$, for every $Q\in\mathcal{V}_{x}$ we
have $\int\widetilde{F}\, dQ=\dim Q\geq\beta$, and so, since $\widetilde{F}$
is continuous we must have 
\[
\liminf_{T\rightarrow\infty}\frac{1}{T}\int_{0}^{T}\widetilde{F}(\mu_{x,t}^{\sqr})\, dt\geq\beta
\]
since otherwise there would be accumulation points $Q$ of $\left\langle \mu\right\rangle _{x,T}^{\sqr}$
with mean $\widetilde{F}$-value smaller than $\beta$. Since 
\[
\frac{1}{N}\sum_{n=1}^{N}\widetilde{F}(\mu_{x,n\log b}^{\sqr})=\int_{0}^{1}\frac{1}{N}\sum_{n=1}^{N}F(S_{t}^{\sqr}\mu_{x,n\log b})dt
\]
By the identity \eqref{eq:5}, the integrand of the right hand side
is just $\mu(B_{b^{-n-t}}(x))/\mu(B_{b^{-t}}(x))$, and so the integrand
is constant up to a multiplicative factor. Also clearly $\frac{1}{T}\int_{0}^{T}\widetilde{F}(\mu_{x,t}^{\sqr})\, dt$
has the same asymptotics whether we allow $T\rightarrow\infty$ along
the reals or multiples of $\log b$. It follows that 
\[
\liminf_{N\rightarrow\infty}\frac{1}{N\log b}\sum_{n=1}^{N}F(\mu_{x,n\log b+t}^{\sqr})\geq\beta
\]
and hence $\underline{D}_{\mu}(x)\geq\beta$.
\end{proof}
We next show that there is a unique FD of maximal dimension. A similar
statement for Z\"{a}hle distributions was proved in \cite{Zahle88}.
\begin{prop}
\label{pro:uniquness-of-d-dimensional-FD}The only $d$-dimensional
FD on $\mathbb{R}^{d}$ is the point mass at $\lambda^{*}$.\end{prop}
\begin{proof}
By the ergodic decomposition (Theorem \ref{thm:EFD-ergodic-components})
it suffices to prove this for EFDs. Let $P$ be a $d$-dimensional
EFD on $\mathbb{R}^{d}$, and let $Q$ be an ergodic CP-distribution
with $\cnt Q=P$. By \cite{Furstenberg08}, the dimension of $Q$
(defined as the $Q$-a.s. dimension of a measure) is given by
\[
d=\int\frac{1}{\log b}H(\mu,\mathcal{D}_{b})\, dP(\mu)
\]
Since the maximal entropy of any measure on $[-1,1]^{d}$ with respect
to $\mathcal{D}_{b}$ is $d\log b$, we conclude that $P$-a.e. $\mu$
gives equal mass to each of the $\mathcal{D}_{b}$-cells of $[-1,1]^{d}$.
Iterating $M_{b}$ for a $Q$-typical $\mu$ and $\mu$-typical point,
we find that every cell $D\in\mathcal{D}_{b}^{n}$ of positive $\mu$-mass
has mass $b^{-n}$. It follows that $\mu_{B_{1}}=\lambda_{B_{1}}$.
Thus, a.s. over the choice of $\mu$, we find that $\mu_{B_{1}}$
generates $\delta_{\lambda^{*}}$; but it also generates $P$. Therefore
$P=\delta_{\lambda^{*}}$.
\end{proof}
USMs do not enjoy this property. In fact any measure $\mu$ on $\mathbb{R}^{d}$
which has exact dimension $d$ is a USM and generates the unique $d$-dimensional
FD. Indeed at $\mu$-a.e. $x$ every $Q\in\mathcal{V}_{x}$ has dimension
$d$ (Proposition \ref{pro:dimension-via-accumulation-pts}), and
hence $\mathcal{V}_{x}=\{\delta_{\lambda^{*}}\}$. Since there is
only one accumulation point it is in fact the limit of the scenery
distribution, so $\mu$ generates $\delta_{\lambda^{*}}$ at $x$.
However, there exist many different measures $\mu$ on $[-1,1]^{d}$
with dimension $d$, including ones that are singular with respect
to Lebesgue measure.
\begin{prop}
The only $0$-dimensional FD on $\mathbb{R}^{d}$ is the point mass
at $\delta_{0}$.\end{prop}
\begin{proof}
Let $P$ be an FD with $\dim(P)=0$. Then for every $n\in\mathbb{N}$,
by Lemma \ref{lem:dimension-of-FDs}, $\int\frac{\log\mu(B_{1/n}(0))}{\log(1/n)}dP(\mu)=\dim P=0$.
It follows that $\mu(B_{1/n}(0))=1$ for $P$-a.e. $\mu$ and all
$n$, which implies that $\mu=\delta_{0}$ a.s., as claimed.\end{proof}
\begin{cor}
If $\mu\in\mathcal{M}$ is exact dimensional then $\dim\mu=0$ if
and only if $\mu$ generates $\delta_{\delta_{0}}$, and $\ldim\mu=0$
if and only if $\delta_{\delta_{0}}\in\mathcal{V}_{x}$ for $\mu$-a.e.
$x$.\end{cor}
\begin{proof}
Immediate from the last proposition and Proposition \ref{pro:dimension-via-accumulation-pts}.
\end{proof}

\subsection{\label{sub:Dimension-of-projected-measures}Dimension of projected
measures }

In this section we prove Theorem \ref{thm:projection-of-SMs}. The
theorem is a generalization of \cite[Theorems 8.1 and 8.2]{HochmanShmerkin09}
from the case of a measure which ``generates'' an ergodic CP distribution
along $b$-adic sceneries, to the case of a measure which does not
display such regular behavior. However, the methods are much the same.
One of the main tools is the ''local entropy averages'' lemma from
\cite{HochmanShmerkin09}, which we give here in a more convenient
form.
\begin{lem}
[Local entropy averages lemma] \label{lem:local-entropy-averages}Let
$\mu\in\mathcal{P}(\mathbb{R}^{d})$ and let $b\geq2$ and $m\geq1$
be integers. Then for $\mu$-a.e. $x$,
\[
\underline{D}(\mu,x)=\liminf_{N\rightarrow\infty}\frac{1}{N}\sum_{n=0}^{N-1}\frac{1}{m\log b}H(\mu_{\mathcal{D}_{b^{n}}(x)},\mathcal{D}_{b^{n+m}})
\]
\end{lem}
\begin{proof}
Let 
\[
I_{n}(x)=-\log\frac{\mu(\mathcal{D}_{b^{n+m}}(x))}{\mu(\mathcal{D}_{b^{n}}(x))}
\]
Then for $i\geq0$ we have
\[
\frac{1}{N}\sum_{n=0}^{N-1}I_{nm+i}(x)=-\frac{1}{N}\log\prod_{n=0}^{N-1}\frac{\mu(\mathcal{D}_{b^{nm+i+m}}(x))}{\mu(\mathcal{D}_{b^{nm+i}}(x))}=-\frac{1}{N}\log\mu(\mathcal{D}_{b^{(N-1)m+i}}(x))+o(1)
\]
as $N\rightarrow\infty$. Averaging this above over $i=0,\ldots,m$
we have
\[
\frac{1}{N}\sum_{n=0}^{N-1}I_{nm+i}(x)=-\frac{1}{N}\log\mu(\mathcal{D}_{b^{(N-1)m+i}})+o(1)
\]
and taking $n\rightarrow\infty$, 
\[
\underline{D}(\mu,x)=\liminf_{n\rightarrow\infty}\frac{1}{N}\sum_{n=0}^{N-1}\frac{1}{m\log b}I_{nm+i}(x)
\]
As usual, we use $\mathcal{D}_{n}$ also to denote the $\sigma$-algebra
generated by $\mathcal{D}_{n}$. Then for $i=0,\ldots,m-1$, 
\[
\frac{1}{N}\sum_{n=0}^{N-1}\left(\mathbb{E}(I_{nm+i}\,|\,\mathcal{D}_{b^{nm+i}})(x)-I_{nm+i}(x)\right)
\]
is an average of uniformly $L^{2}$-bounded martingale differences.
Indeed, to verify $L^{2}$ boundedness, note that the function $x\log^{2}x$,
which arises when integrating the second power of the information
function, is bounded on $[0,1]$. Also, the $n$-th term in the sum
is clearly $\mathcal{D}_{(n+m)m+i}$-measurable, giving the \,martingale
property. So by the Law of Large Numbers for martingale differences
\cite[Theorem 3 in Section 9]{Feller71}, it converges to $0$ a.e.
Therefore, by the discussion above and Lemma \ref{lem:decay-along-filtration},
for $\mu$-a.e. $x$,
\[
\underline{D}(\mu,x)=\liminf_{N\rightarrow\infty}\frac{1}{N}\sum_{n=0}^{N-1}\mathbb{E}(I_{n}\,|\,\mathcal{D}_{b^{n}})(x)
\]
Finally, observe that 
\[
\mathbb{E}(I_{n}\;|\;\mathcal{D}_{b^{n}})(x)=H(\mu_{\mathcal{D}_{b^{n}}(x)},\mathcal{D}_{b^{n+m}})
\]
Substituting this into the last equation gives the lemma.
\end{proof}
Returning to Theorem \ref{thm:projection-of-SMs}, we first establish
the case in which the projecting map is linear, and remove this assumption
later. Thus, let $\mu\in\mathcal{M}$ and $\pi\in\Pi_{d,k}$, and
suppose that $\ldim\pi\mu<\alpha$. Our objective is to show that
there is a set of $x$ of positive $\mu$-measure for which there
exits $P\in\mathcal{V}_{x}$ with $E_{\pi}(P)\geq\alpha$, where $\mathcal{V}_{x}$
is the set of accumulation points of $\left\langle \mu\right\rangle _{x,t}$.

We begin with some reductions. Since $\ldim\pi\mu=\essinf_{y\sim\pi\mu}\underline{D}(\pi\mu,y)$,
there is an $\varepsilon>0$ and a set $E\subseteq\mathbb{R}^{k}$
with $\pi\mu(E)>0$ and $\underline{D}(\pi\mu,y)<\alpha-\varepsilon$
for all $y\in E$. Let $\mu'=\mu_{\pi^{-1}E}$. Then for $\mu'$-a.e.
$x$ the accumulation points of $\left\langle \mu'\right\rangle _{x,t}$
are precisely $\mathcal{V}_{x}$ (Proposition \ref{pro:equivalent-measures-have-same-sceneries}),
and $\ldim\pi\mu'\leq\alpha-\varepsilon$. Thus, we may assume from
the start that this is true for $\mu$. Restricting $\mu$ further,
we may also assume that it has bounded support, and we can normalize
it to be a probability measure, since normalization affects neither
$\ldim\pi\mu$ nor $\mathcal{V}_{x}$. Finally, we may choose our
coordinate system so that $\pi(x_{1},\ldots,x_{d})=(x_{1},\ldots,x_{k})$
is the coordinate projection and $\mu\in\mathcal{M}^{\sqr}$. 

Next, in order to relate $\ldim\pi\mu$ to the behavior of $H(\pi(\mu_{A}),\mathcal{D}_{b^{n+m}})$
for $A\in\mathcal{D}_{b^{n}}^{d}$, we use the following device, taken
from \cite{HochmanShmerkin09}.
\begin{lem}
\label{lem:projection-decomposition}Let $\mu\in\mathcal{M}^{\sqr}$.
For every integer $b\geq2$ there exists a standard Borel probability
space $(\Omega,\mathcal{F},P)$ and measurable family of measures
$\{\mu_{\omega}\}_{\omega\in\Omega}\subseteq\mathcal{M}^{\sqr}$ such
that
\begin{enumerate}
\item $\mu=\int\mu_{\omega}\, dP(\omega)$.
\item For $P$-a.e. $\omega$, for every $n$ and $A\in\mathcal{D}_{b^{n}}^{d}$
with $\mu_{\omega}(A)>0$, we have $H(\pi((\mu_{\omega})_{A}),\mathcal{D}_{b^{n+1}}^{k})=H(\pi(\mu_{A}),\mathcal{D}_{b^{n+1}}^{k})$.
\end{enumerate}
\end{lem}
\begin{proof}
We define a random sequence $\mu_{n}$, $n=1,2,\ldots$, where $\mu_{n}$
is a probability measure on $[-1,1]^{d}$ with and the algebra of
sets generated by $\mathcal{D}_{b^{n}}^{d}$. We construct the sequence
so that 
\begin{enumerate}
\item [(0')] $\mu_{n+1}$ extends $\mu_{n}$ ( $\mu_{n+1}|_{\mathcal{D}_{b^{n}}^{d}}=\mu_{n}$).
\item [(1')] $\mathbb{E}\mu_{n}=\mu|_{\mathcal{D}_{b^{n}}^{d}}$.
\item [(2')] With probability one, for every $A\in\mathcal{D}_{b^{n}}^{d}$
with $\mu_{n}(A)>0$, the measures $\mu_{A}$ and $(\mu_{n})_{A}$
agree on the sets $\pi^{-1}E$ for $E\in\mathcal{D}_{b^{n+1}}^{k}$. 
\end{enumerate}
Let $(\Omega,\mathcal{F},P)$ be the underlying probability space,
which clearly can be realized as a standard Borel space. By (0') we
have $\mathbb{E}(\mu_{n+1}|\mathcal{\mu}_{n})=\mu_{n}$, and by (1')
we have $\mathbb{E}(\mu_{n})=\mu$, so by the martingale theorem for
measures, $\widetilde{\mu}=\lim_{n\rightarrow\infty}\mu_{n}$ exists
a.s. and $\mathbb{E}(\widetilde{\mu})=\mu$. Finally, for $A\in\mathcal{D}_{b^{n}}^{d}$,
\[
H(\pi\widetilde{\mu}_{A},\mathcal{D}_{b^{n}}^{k})=H(\pi(\mu_{n})_{A},\mathcal{D}_{b^{n}}^{k})=H((\mu_{n})_{A},\pi^{-1}\mathcal{D}_{b^{n}}^{k})=H(\mu_{A},\pi^{-1}\mathcal{D}_{b^{n}}^{k})
\]
where the last inequality is because by (2') $\mu_{A}$ and $(\mu_{n})_{A}$
agree on every element of $\pi^{-1}\mathcal{D}_{b^{n}}^{k}$. This
gives (2).

We define $\mu_{n}$ by a random recursive construction. Let $\mu_{0}$
be the unique probability measure on $\mathcal{D}_{1}^{d}$. Assume
we have defined $\mu_{n}$. For each cell $A\in\mathcal{D}_{b^{n}}^{d}$
with $\mu_{n}(A)>0$ we must specify how to split the mass $\mu_{n}(A)$
between the cells $A'\in\mathcal{D}_{b^{n+1}}^{d}$, $A'\subseteq A$.
Choose $E\in\mathcal{D}_{b^{n+1}}^{k}$ with $\mu_{n}(\pi^{-1}E)>0$.
Choose a cell $A'\in\mathcal{D}_{b^{n+1}}^{d}$ with $A'\subseteq E$
according to its relative $\mu$-mass (so we choose $A'\in\mathcal{D}_{b^{n}}^{d}$
with probability $\mu(A'\cap A)/\mu(E\cap A)$). Make these choices
independent for such each $E$. Finally, for each $E$ and $A'$ define
\[
\mu_{n}(A')=\mu_{n}(A)\cdot\frac{\mu(E\cap A)}{\mu(A)}
\]
Observe that in each ``column'' $E$ we have chosen a single cell
and given it all the mass of that column; thus (2') is satisfied.
On the other hand, (1') is also satisfied, since the distributions
of the sets $E$ was chosen the same as $\mu$, and given $E$, so
was the mass of $A'\subseteq E$. 
\end{proof}
Fix an arbitrary integer base $b$ and let $(\Omega,\mathcal{F},P)$
and $\mu=\int\mu_{\omega}dP(\omega)$ be as in the lemma above. Then
$\pi\mu=\int\pi\mu_{\omega}dP(\omega)$, hence 
\[
\essinf_{\omega\sim P}\ldim\pi\mu_{\omega}\leq\ldim\pi\mu\leq\alpha-\varepsilon
\]
In other words, for $P$-a.e. $\omega$, for $\mu_{\omega}$-a.e.
$x$, $\underline{D}(\pi\mu_{\omega},\pi x)\leq\alpha-\varepsilon$.
Using Lemma \ref{lem:local-entropy-averages}, this means that for
$P$-a.e. $\omega$, for $\mu_{\omega}$-a.e. $x$, for every $m$,
\begin{equation}
\liminf_{N\rightarrow\infty}\frac{1}{N}\sum_{n=1}^{N}\frac{1}{m\log b}H((\pi\mu_{\omega})_{\mathcal{D}_{b^{n}}^{k}(\pi x)},\mathcal{D}_{b^{n+m}})\leq\alpha-\varepsilon\label{eq:12}
\end{equation}

Consider the random pair $(\omega,x)\in\omega\times\mathbb{R}^{d}$
generated by first choosing $\omega\sim P$ and then $x\sim\mu_{\omega}$.
Let $\widetilde{P}=\int\delta_{\omega}\times\mu_{\mu}$ denote the
measure on $\Omega\times\mathbb{R}^{d}$ corresponding to this random
pair. The marginal distribution of $x$ is then $\int\mu_{\omega}dP(\omega)=\mu$,
and we can disintegrate \texttt{$\widetilde{P}$} according to the
second coordinate, obtaining $\widetilde{P}=\int\widetilde{P}_{x}d\mu(x)$
where $\widetilde{P}_{x}$ is a probability measure on $(\Omega,\mathcal{F})$,
defined $\mu$-a.e. (we use here the fact that $\Omega$ can be taken
to be a standard Borel space, as noted in Lemma \ref{lem:projection-decomposition},
and that the coordinate projections on $\Omega\times\mathbb{R}^{d}$
are measurable). Now the statement at the end of the previous paragraph
becomes the following: For $\mu$-a.e. $x$, for $\widetilde{P}_{x}$-a.e.
$\omega$, \eqref{eq:14} holds. Using property (2) of Lemma \ref{lem:projection-decomposition},
we have shown that 
\begin{equation}
\forall m\qquad\liminf_{N\rightarrow\infty}\frac{1}{N}\sum_{n=1}^{N}\frac{1}{m\log b}H(\mu_{\mathcal{D}_{b^{n}}^{d}(\pi x)},\pi^{-1}\mathcal{D}_{b^{n+m}})\leq\alpha-\varepsilon\qquad\mu\mbox{-a.e. }x\label{eq:13}
\end{equation}

As in Section \ref{sec:Equivalence-of-the}, denote 
\[
\left\langle \mu,x\right\rangle _{N}=\frac{1}{N}\sum_{n=1}^{N}\delta_{(M_{b}^{\sqr})^{n}(\mu,x)}
\]
Let $\mathcal{W}_{x}$ denote the set of accumulation points of $\left\langle \mu,x\right\rangle _{N}$.
We can assume that the coordinate system has been chosen (equivalently,
$\mu$ translated) so that for $\mu$-a.e $x$, every $Q\in\mathcal{W}_{x}$
is a base-$b$ CP-distribution, and every $Q\in\mathcal{V}_{x}$ is
the centering of some $Q\in\mathcal{W}_{x}$. Thus, our assumption
is that for $\mu$-a.e. $x$, for every $Q\in\mathcal{W}_{x}$, $E_{\pi}(Q)\geq\alpha$.
Then we can re-write \eqref{eq:13} as
\[
\liminf_{N\rightarrow\infty}\int\frac{1}{m\log b}H(\theta,\pi^{-1}\mathcal{D}_{b^{m}})d\left\langle \mu,x\right\rangle _{N}(\theta)\leq\alpha-\varepsilon\qquad\mu\mbox{-a.e. }x
\]
where we use the usual convention that a measure on $\mathcal{P}(B_{1})\times B_{1}$
such as $\left\langle \mu,x\right\rangle _{N}$ is implicitly identified
with its marginal on $\mathcal{P}(B_{1})$. Replacing $H(\cdot,\mathcal{D}_{b^{m}}^{k})$
in \eqref{eq:13} by a continuous approximation $f_{m}$ as in Lemma
\ref{lem:continuous-approximation-of-entropy}, we conclude that for
$\mu$-a.e. $x$ there is some $Q\in\mathcal{W}_{x}^{b}$ such that
\[
\forall m\qquad\int\frac{1}{m\log b}H(\theta,\pi^{-1}\mathcal{D}_{b^{m}})dQ(\theta)\leq\alpha-\varepsilon+\frac{C}{m\log b}
\]
which implies that
\[
\limsup_{m\rightarrow\infty}\int\frac{1}{m\log b}H(\pi\theta,\mathcal{D}_{b^{m}})dQ(\theta)\leq\alpha-\varepsilon
\]
To conclude the proof, we rely on the fact from \cite[Theorems 8.1 and 8.2]{HochmanShmerkin09}
that if $Q$ is ergodic then $E_{\pi}(Q)$ is precisely the left-hand
side of the last equation (and the limsup is in fact a limit). Thus
if $Q$ is ergodic we have obtained our contradiction, finding, for
almost every $x$, a $Q\in\mathcal{W}_{x}$ (and hence $\cnt Q\in\mathcal{V}_{x}$)
with $E_{\pi}(Q)<\alpha$. If $Q$ is not ergodic, it has an ergodic
decomposition $Q=\int Q_{\theta}dQ(\theta)$ with respect to $M$.
Since $0\leq\frac{1}{m\log b}H(\theta,\mathcal{D}_{b^{m}}^{k})\leq d$
for $\theta\in\mathcal{M}^{\sqr}(\mathbb{R}^{k})$, by the above and
bounded convergence,
\begin{eqnarray*}
E_{\pi}(Q) & = & \int E_{\pi}(Q_{\theta})dQ(\theta)\\
 & = & \int\lim_{m\rightarrow\infty}\int\frac{1}{m\log b}H(\pi\theta,\mathcal{D}_{b^{m}}^{k})\, dQ_{\theta}(\nu)dQ(\theta)\\
 & = & \lim_{m\rightarrow\infty}\int\int\frac{1}{m\log b}H(\pi\nu,\mathcal{D}_{b^{m}}^{k})\, dQ_{\theta}(\nu)dQ(\theta)\\
 & = & \lim_{m\rightarrow\infty}\int\frac{1}{m\log b}H(\pi\theta,\mathcal{D}_{b^{m}}^{k})\, dQ(\theta)
\end{eqnarray*}
which by the above is $\leq\alpha-\varepsilon$, and we have again
reached our goal.

Finally let us treat the case that instead of a linear map we have
a regular map $\varphi\in C^{1}(\mathbb{R}^{d},\mathbb{R}^{k})$.
This case can be reduced to the linear one by ``straightening out''
the map at the expense of ``twisting'' the measure. First, we note
that it suffices to prove the claim locally: if $U$ is open and $\mu(U)>0$,
then we may consider $\mu_{U}$, since $\mu_{U}$-a.e. the sceneries
have the same accumulation points as for $\mu$, and $\ldim\varphi\mu=\inf\ldim\varphi\mu_{U}$
where the infimum is taken over a fixed countable open cover of $\mathbb{R}^{d}$
by open sets of finite measure. Thus we may assume that $\mu=\mu_{U}$
and $U$ is small enough that we can apply the implicit function theorem
in it and find a $C^{1}$ local diffeomorphism $\psi:U\rightarrow\mathbb{R}^{d}$,
so that $\varphi=\pi\psi$, where $\pi$ is a fixed linear map (e.g.
projection onto the first $k$ coordinates). Thus it suffices to bound
$\ldim\pi\nu$, where $\nu=\psi\mu$. By Proposition \ref{pro:smooth-images-of-usms},
the set $\mathcal{V}'_{\psi x}=\{D\psi(x)^{*}P\,:\, P\in\mathcal{V}_{x}\}$
is the set of distribution generated by $\nu$ at $\psi x$. Thus
by the theorem above (in the linear case), 
\[
\ldim\varphi\mu=\ldim\pi\nu\geq\essinf_{y\sim\nu}\inf_{P\in\mathcal{V}'_{\psi x}}E_{Q}(\pi)=\essinf_{x\sim\mu}\inf_{P\in\mathcal{V}_{x}}E_{D\psi(x)^{*}P}(\pi)
\]
Now, clearly $E_{D\psi(x)^{*}P}(\pi)=E_{P}(\pi\circ D\psi(x))=E_{P}(D\varphi(x))$,
since $\varphi=\pi\psi$ implies by the chain rule that $D\varphi(x)=\pi\circ D\psi(x)$.
Thus the formula above is the desired one.

The theorem remains valid under the weaker assumption of differentiability
of $\varphi$ rather than $C^{1}$, but then one must work directly
with suitable partitions of $\mathbb{R}^{d}$, chosen so that locally,
neighborhoods of the level sets of $Df$ are unions of atoms of the
partition (as is the case for $\mathcal{D}_{b}$ when $f$ is a coordinate
projection). We omit the details.

\subsection{\label{sub:Exact-dimension-of-of-conditionals}Conditional measures }

Let $\pi\in\Pi_{d,k}$. The conditional measures of $\mu\in\mathcal{P}(\mathbb{R}^{d})$
on the fibers of $\pi$ were discussed in Section \ref{sub:Conditional-measures}.
We take a moment to discuss in more detail the case of general (non-probability)
measure $\mu\in\mathcal{M}$. We need the following technical lemma
whose proof we provide for completeness.
\begin{lem}
\label{lem:restriction-of-conditionals}Suppose $\mu\in\mathcal{P}(\mathbb{R}^{d})$
and let $A$ be a measurable set with $\mu(A)>0$. Then $(\mu_{A})_{[x]_{\pi}}=(\mu_{[x]_{\pi}})_{A}$
for $\mu_{A}$-a.e. $x$. \end{lem}
\begin{proof}
For any measurable set $E$, by definition
\begin{eqnarray*}
\mu_{A}(E) & = & \int(\mu_{A})_{[x]_{\pi}}(E)\, d\mu_{A}(x)
\end{eqnarray*}
On the other hand,
\begin{alignat*}{2}
\mu_{A}(E)\quad & =\quad\frac{1}{\mu(A)}\mu(A\cap E) & \quad & \mbox{definition of }\mu_{A}\\
 & =\quad\frac{1}{\mu(A)}\int\mu_{[x]_{\pi}}(A\cap E)\, d\mu(x) &  & \mu=\int\mu_{[x]_{\pi}}d\mu(x)\\
 & =\quad\frac{1}{\mu(A)}\int\mu_{[x]_{\pi}}(A)\cdot(\mu_{[x]_{\pi}})_{A}(E)\, d\mu(x) &  & \mbox{definition of }(\mu_{[x]_{\pi}})_{A}\\
 & =\quad\frac{1}{\mu(A)}\int(\int1_{A}(y)\, d\mu_{[x]_{\pi}}(y))\cdot(\mu_{[x]_{\pi}})_{A}(E)\, d\mu(x) &  & \mu_{[x]_{\pi}}(A)=\int1_{A}d\mu_{[x]_{\pi}}\\
 & =\quad\frac{1}{\mu(A)}\int\int1_{A}(y)\cdot(\mu_{[y]_{\pi}})_{A}(E)\, d\mu_{[x]_{\pi}}(y)\, d\mu(x) &  & \mu_{[x]_{\pi}}=\mu_{[y]_{\pi}}\;,\;\mu_{[x]_{\pi}}\mbox{-a.e. }y\\
 & =\quad\frac{1}{\mu(A)}\int1_{A}(y)\cdot(\mu_{[y]_{\pi}})_{A}(E)\, d\mu(y) &  & \mu=\int\mu_{[x]_{\pi}}d\mu(x)\\
 & =\quad\int(\mu_{[y]_{\pi}})_{A}(E)\, d\mu_{A}(y) &  & \mbox{definition of }\mu_{A}
\end{alignat*}

Since this holds for every $E$, and since $(\mu_{[y]_{\pi}})_{A}$
is supported on $[y]_{\pi}$ $\mu_{A}$-a.s., by uniqueness of the
decomposition of $\mu_{A}$ we have $(\mu_{[y]_{\pi}})_{A}=(\mu_{A})_{[y]_{\pi}}$
$\mu_{A}$-a.e. \end{proof}
\begin{cor}
For $\mu\in\mathcal{P}(\mathbb{R}^{d})$ and Borel sets $A\subseteq B$
with $\mu(A)>0$ we have $(\mu_{A})_{[x]_{\pi}}=c\cdot(\mu_{B})_{[x]_{\pi}}|_{A}$
for $\mu$-a.e. $x\in A$.\end{cor}
\begin{proof}
Apply the previous lemma to $\mu_{B}$. 
\end{proof}
Returning to $\mu\in\mathcal{M}$, for large enough $R$ that $\mu(B_{R})>0$,
the conditional measures of $\mu_{B_{R}}$ with respect to $\pi$
are defined $\mu_{B_{R}}$-a.e. and by the corollary, if $R'>R$ then
$(\mu_{B_{R'}})_{[x]_{\pi}}|_{B_{R}}=c\cdot(\mu_{B_{R}})_{[x]_{\pi}}$
for a constant $c$ (see Lemma \ref{lem:restriction-of-conditionals}).
Therefore it is possible to choose constants $c_{R}=c_{R}(x)$ such
that $c_{R}\cdot(\mu_{B_{R}})_{[x]_{\pi}}$ are a consistent, in the
sense that, for $R'>R$ we have $c_{R'}\cdot(\mu_{B_{R'}})_{[x]_{\pi}}|_{B_{R}}=c_{R}\cdot(\mu_{B_{R}})_{[x]_{\pi}}$,
and $c_{R}\cdot(\mu_{B_{R}})_{[x]_{\pi}}$ form an increasing family
of measures, allowing us to define 
\[
\mu_{[x]_{\pi}}=\lim_{R\rightarrow\infty}c_{R}\cdot(\mu_{B_{R}})_{[x]_{\pi}}
\]
(the limit is in the sense that $\mu_{[x]_{\pi}}(A)=\lim c_{R}(\mu_{B_{R}})_{[x]_{\pi}}(A)$
for all measurable $A$; the limit exists by monotonicity and a standard
verification show that it is a measure). The measure we obtain depends,
up to a constant, on our choice of $c_{R}$. Similarly, for $\mu$-a.e.
$x$ and $\mu_{[x]_{\pi}}$-a.e. $y$ we will have $\mu_{[y]_{\pi}}=c\cdot\mu_{[x]_{\pi}}$,
since by Lemma \ref{lem:restriction-of-conditionals} again, the conditionals
measures of $\mu_{B_{R}(x)\cap B_{R}(y)}$ agree with those of $\mu_{B_{R}(x)}$
and $\mu_{B_{R}(y)}$ on $B_{r}(x)\cap B_{R}(y)$. In summary, up
to a scalar, $\mu_{[x]_{\pi}}$ is well defined and depends only on
$[x]_{\pi}$ (rather than $x$). From now on we abuse notation and
treat it as a true measure. Furthermore, we note that if $\mu_{[x]_{\pi}}(B_{1})>0$
then $(\mu_{[x]_{\pi}})^{*}$ is a bona-fide measure, even though
$\mu_{[x]_{\pi}}$ is not.

Now, the conditional measures of $\mu$ exist a.e., but they may not
defined for a particular $x$. There is a more constructive way to
define conditional measures, namely, as limits of conditionals on
``thickened'' fibers:
\[
\mu_{[x]_{\pi}}=\lim_{R\rightarrow\infty}(\lim_{\varepsilon\rightarrow0}c_{R,\varepsilon}\cdot\mu_{B_{R}\cap[x]_{\pi}^{(\varepsilon)}})
\]
where $A^{(\varepsilon)}$ is the $\varepsilon$-neighborhood of $A$
and $c_{R,\varepsilon}$ are appropriate constants. But this still
does not have to be well-defined at a given $x$. 

However, if $P$ is an EFD then by the quasi-Palm property, since
for $P$-typical $\mu$ the conditional measures are defined a.e.,
for $P$-typical $\mu$ it is also defined for $x=0$. Furthermore,
we claim that $(\mu_{[0]_{\pi}})^{*}$ is well-defined. This follows
from the following lemma and the quasi-Palm property:
\begin{lem}
If $\nu$ is a probability measure the $y\in\supp\nu_{[y]_{\pi}}$
for $\nu$-a.e. $y$.\end{lem}
\begin{proof}
Let $Y=\{y\,:\, y\notin\supp\nu_{[y]_{\pi}}\}$. Then, since $\nu_{[y]_{\pi}}$
is a.s. supported on $[y]_{\pi}$,
\[
\nu(Y)=\int\nu_{[y]_{\pi}}(Y)\, d\nu(y)=\int\nu_{[y]_{\pi}}(Y\cap[y]_{\pi})\, d\nu(y)\leq\int\nu_{[y]_{\pi}}(\mathbb{R}^{d}\supp\nu_{[y]_{\pi}})\, d\nu(y)=0
\]
as claimed.
\end{proof}
Thus, by uniqueness (up to normalization) of the conditional measures
it is clear that 
\[
(S_{t}^{*}\mu)_{[0]_{\pi}}=S_{t}^{*}(\mu_{[0]_{\pi}})
\]
Hence the push-forward $P_{[0]_{\pi}}$ of $\mu$ via $\mu\mapsto\mu_{[0]_{\pi}}$
is an $S^{*}$-invariant distribution (which is defined on $\mathbb{R}^{d}$
but its measures are supported on $[0]_{\pi}=\ker\pi$, which may
be identified with $\mathbb{R}^{d-k}$).
\begin{claim}
$P_{[0]_{\pi}}$ is an an EFD. \end{claim}
\begin{proof}
We already saw that it is $S^{*}$-invariant so we must only show
the quasi-Palm property. Let $\mathcal{A}\subseteq\mathcal{M}$ be
an event with $P_{[0]_{\pi}}(\mathcal{A})=0$. Let $\mathcal{A}'$
denote the set of measures $\nu$ such that $\nu_{[0]_{\pi}}$ is
defined and $\nu_{[0]_{\pi}}\in\mathcal{A}$. Since $P_{[0]_{\pi}}(\mathcal{A})=0$,
by definition $P(\mathcal{A}')=0$. By the quasi-Palm property of
$P$, we know that for $P$-a.e. $\mu$, for $\mu$-a.e. $x$, $T_{x}^{*}\mu\notin\mathcal{A}'$,
equivalently $(T_{x}^{*}\mu)_{[0]_{\pi}}\notin\mathcal{A}$. Using
the basic properties of the conditional measures, this means that
for $\mu$-a.e. $x$, for $\mu_{[x]_{\pi}}$-a.e. $y$, $(T_{y}^{*}\mu)_{[0]_{\pi}}\notin\mathcal{A}$.
Equivalently, that for $P$-a.e. $\mu$ and $\mu$-a.e. $x$, for
$(T_{x}^{*}\mu)_{[0]_{\pi}}$-a.e. $y$ we have $T_{y-x}^{*}((T_{x}^{*}\mu)_{[0]_{\pi}})\notin\mathcal{A}$.
By the quasi-Palm property of $P$ again, this means the for $P$-a.e.
$\mu$, for $\mu_{[0]_{\pi}}$-a.e. $y$, $T_{y}^{*}(\mu_{[0]_{\pi}})\notin\mathcal{A}$.
This is precisely the quasi-Palm property of $P_{[0]_{\pi}}$.
\end{proof}
Since EFDs are exact dimensional this establishes Proposition \ref{pro:exact-dimension-of-conditionals}.
As noted in the introduction, Theorem \ref{thm:fibers-for-EFDs-1}
on dimension conservation follows from Furstenberg's corresponding
result for coordinate projections and CP-processes, and the ability
to change coordinate systems (Corollary \ref{cor:CP-change-of-coords}).

We turn to the proof of Theorem \ref{thm:conditional-measures-of-USMs}
on the conditional measures of USMs. Arguing in the same way as for
projections, for regular $C^{1}$ maps $f:\mathbb{R}^{d}\rightarrow\mathbb{R}^{k}$,
by a (non-linear) change of coordinates it is enough to consider linear
projections; we only consider this case. With another (linear) change
of coordinates, we may further assume that the coordinate projection
is $f(x_{1},\ldots,x_{d})=(x_{1},\ldots,x_{k})$. We use superscripts
to denote the ambient space $\mathbb{R}^{d}$ and $\mathbb{R}^{k}$
as needed. Next, we need the following well-known fact.
\begin{lem}
\label{lem:Fubini-for-dimension}Suppose that $\mu\in\mathcal{P}([-1,1]^{d})$
and, for $\pi:\mathbb{R}^{d}\rightarrow\mathbb{R}^{k}$ the coordinate
projection, $\mu$ is such that $\ldim\pi\mu\geq\alpha$ and $\ldim\mu_{[x]_{\pi}}\geq\beta$
for $\mu$-a.e. $x$. Then $\ldim\mu\geq\alpha+\beta$. \end{lem}
\begin{proof}
We give the easy proof for completeness. Write $\nu=\pi\mu$. We must
show, for $\mu$-typical $x$, that $\mu(\mathcal{D}_{2^{n}}^{d}(x))\leq2^{-(\alpha+\beta-o(1))n}$.
Using Egorov's theorem, we can find a compact set $A\subseteq\mathbb{R}^{k}$
of $\nu$-measure arbitrarily close to $1$ on which $\nu(\mathcal{D}_{2^{n}}^{k}(y))<2^{-(\alpha-o(1))n}$
uniformly for $y\in A$; and a compact set $B\subseteq\mathbb{R}^{d}$
of measure arbitrarily close to $1$ such that $\mu_{[x]_{\pi}}(\mathcal{D}_{2^{n}}^{d}(x))<2^{-(\beta-o(1))n}$
uniformly for $x\in B$. The set $C=B\cap\pi^{-1}A$ still can be
made to have has measure arbitrarily close to $1$, and it suffices
to prove that $\mu(\mathcal{D}_{2^{n}}^{d}(x))\leq2^{-(\alpha+\beta-o(1))n}$
for $\mu$-a.e. $x\in C$. By the Besicovitch density theorem (Theorem
\ref{lem:Besicovitch-lemma}) this is equivalent to showing that $\mu(C\cap\mathcal{D}_{2^{n}}^{d}(x))\leq2^{-(\alpha+\beta-o(1))n}$
for $\mu$-a.e. $x\in C$. This is what we will show.

Fix $\mu$-typical $x\in C$ and let $I_{x,n}=\mathcal{D}_{2^{n}}^{d}(x)$
and $J_{x,n}=\pi I_{x,n}=\mathcal{D}_{2^{n}}^{k}(\pi x)$. By definition
of of the conditional measures, 
\begin{eqnarray*}
\mu(C\cap\mathcal{D}_{2^{n}}^{d}(x)) & = & \int\mu_{[z]_{\pi}}(C\cap I_{x,n})\, d\mu(z)\\
 & = & \int\mu_{\pi^{-1}(y)}(C\cap I_{x,n})\, d\nu(y)\\
 & = & \nu(J_{x,n})\cdot\int\mu_{\pi^{-1}(y)}(C\cap I_{x,n})\, d\nu_{J_{x,n}}(y)\\
 & < & 2^{-(\alpha-o(1))n}\cdot\int\mu_{\pi^{-1}(y)}(C\cap I_{x,n})\, d\nu_{J_{x,n}}(y)\\
 & < & 2^{-(\alpha-o(1))n}\cdot\int2^{-(\beta-o(1))n}\, d\nu_{J_{x,n}}(y)\\
 & = & 2^{-(\alpha+\beta-o(1))n}
\end{eqnarray*}
where the first inequality is because $J_{x,n}=\mathcal{D}_{2^{n}}^{k}(\pi x)$
and $\pi x\in A$, and the second is because $\mu_{\pi^{-1}(y)}(C\cap I_{x,n})$
is equal either to $\mu_{\pi^{-}(y)}(\emptyset)$ or else to $\mu_{\pi^{-1}(y)}(C\cap\mathcal{D}_{2^{n}}^{d}(z))$
for some $z\in B$, and in both cases the mass is bounded by $2^{-(\beta-o(1))n}$.
\end{proof}
Proceeding with the proof of Theorem \ref{thm:conditional-measures-of-USMs},
we argue as follows. Fix $\mu$ and the coordinate projection $\pi$.
For any measurable set $A$, if $\mu(A)>0$ then by the Besicovitch
density theorem, $\underline{D}(\mu,x)=\underline{D}(\mu_{A},x)$
$\mu_{A}$-a.e. On the other hand, by the lemma above and the density
theorem again, $\underline{D}((\mu_{A})_{[x]_{\pi}},y)=\underline{D}((\mu_{[x]_{\pi}})_{A},y)=\underline{D}(\mu_{[x]_{\pi}},y)$
for $\mu_{A}$-a.e. $x$ and $(\mu_{A})_{[x]_{\pi}}$-a.e. $y$ (equivalently
$(\mu_{[x]_{\pi}})_{A}$-a.e. $y$). Since $\int(\mu_{A})_{[x]_{\pi}}d\mu_{A}(x)=\mu_{A}$,
this implies that $\underline{D}((\mu_{A})_{[x]_{\pi}},x)=\underline{D}(\mu_{[x]_{\pi}},x)$
for $\mu_{A}$-a.e. $x$. Finally, by Proposition \ref{pro:equivalent-measures-have-same-sceneries}
we know that for $\mu_{A}$-a.e. $x$ the accumulation points of $\left\langle \mu_{A}\right\rangle _{x,t}^{\sqr}$
are the same as the accumulation points of $\left\langle \mu\right\rangle _{x,t}^{\sqr}$,
that is, the elements of $\mathcal{V}_{x}$. Thus, it will be enough
to prove the theorem for $\mu_{A}$ for sets $A$ of arbitrarily large
measure.

Let $\varepsilon,\alpha,\beta>0$ and consider the set 
\[
A=A_{\varepsilon,\alpha,\beta}=\{x\,:\,|\underline{D}(\mu,x)-\alpha|<\varepsilon\mbox{ and }|\underline{D}(\mu_{[x]_{\pi}},x)-\beta|<\varepsilon\}
\]
Then by the previous remarks, $\ldim\mu_{A}>\alpha-\varepsilon$ and
$\ldim(\mu_{A})_{[x]_{\pi}}<\beta+\varepsilon$ for $\mu_{A}$-a.e.
$x$. Thus by Lemma \ref{lem:Fubini-for-dimension}, $\ldim\pi\mu_{A}<\alpha-\beta+2\varepsilon$,
so by Theorem \ref{thm:projection-of-SMs},
\[
\inf_{P\in\mathcal{V}_{X}}E_{\pi}(P)<\alpha-\beta+2\varepsilon\qquad\mu_{A}\mbox{-a.e.}
\]
Thus we have shown that for given $\varepsilon,\alpha,\beta$ and
$A=A_{\varepsilon,\alpha,\beta}$ as above,
\[
\underline{D}(\mu_{[x]_{\pi}},x)\leq\alpha-\inf_{P\in\mathcal{V}_{X}}E_{\pi}(P)+2\varepsilon\qquad\mu_{A}\mbox{-a.e.}
\]
Since $\mu(\bigcup_{\alpha,\beta\in\mathbb{Q}}A_{\varepsilon,\alpha,\beta})=1$,
this implies that the same thing holds $\mu$-a.e., and since $\varepsilon$
was arbitrary, this completes the proof of Theorem \ref{thm:conditional-measures-of-USMs}.

\section{\label{sec:FDs-with-additional-invariance}Fractal distributions
with additional invariance}

\subsection{\label{sub:Linear-images-of-FDs}Images of FDs and USMs}

We begin by proving Proposition \ref{prop:linear-images-of-FDs}.
Let $L\in GL_{n}(\mathbb{R})$ and $L^{*}:\mathcal{M}^{*}\rightarrow\mathcal{M}^{*}$
the induced map on measures, which itself induces a map on distributions
which we denote in the same way. Let $P$ be an (extended) FD; we
wish to show that $P'=L^{*}P$ is an FD as well; and if $P$ was ergodic
so is $P'$. Since $L$ is linear it commutes with scaling: $L(ax)=a(Lx)$
for $a\in\mathbb{R}$. This means that $L^{*}$ commutes with $S_{t}^{*}$
for all $t$, hence $P'=L^{*}P$ is $S^{*}$ invariant, and if $P$
was $S^{*}$-ergodic so is $L^{*}P$, being a factor of $P$ by the
factor map $L^{*}$. To see that $P'$ is an EFD we must show that
it is quasi-Palm. Note that if $U\subseteq\mathbb{R}^{d}$ is an open
neighborhood of $0$ then so is $U'=L^{-1}U$, and that $T_{Lx}\circ L=L\circ T_{x}$.
Hence

\begin{alignat*}{2}
\int\left\langle \nu\right\rangle _{U}^{*}\, dL^{*}P(\nu)\quad & =\quad\int\left\langle L^{*}\nu\right\rangle _{U}^{*}dP(\nu) & \quad & \mbox{by definition of }L^{*}P\\
 & =\quad\int\left(\int_{U}\delta_{T_{x}^{*}L\nu}\, dL^{*}\nu(x)\right)dP(\nu) &  & \mbox{by definition of }\left\langle \cdot\right\rangle _{U}\\
 & =\quad\int\left(\int_{L^{-1}U}\delta_{T_{Lx}^{*}L\nu}\, d\nu(x)\right)dP(\nu) &  & \mbox{by definition of }L^{*}\nu\\
 & =\quad\int\left(\int_{U'}\delta_{L^{*}T_{x}^{*}\nu}\, d\nu(x)\right)dP(\nu) &  & U'=L^{-1}U\,,\, LT_{x}=T_{Lx}L\\
 & =\quad\int\left(\int_{U'}L^{*}\delta_{T_{x}^{*}\nu}\, d\nu(x)\right)dP(\nu) &  & L^{*}\delta_{y}=\delta_{L^{*}y}\\
 & =\quad L^{*}\int\left(\int_{U'}\delta_{T_{x}^{*}\nu}\, d\nu(x)\right)dP(\nu) &  & \mbox{by linearity of }L^{*}\\
 & =\quad L^{*}\int\left\langle \nu\right\rangle _{U'}dP(\nu) &  & \mbox{by definition of }\left\langle \cdot\right\rangle _{U}\\
 & \sim\quad L^{*}P
\end{alignat*}
where in the last equality we used that fact that $\int\left\langle \nu\right\rangle _{U'}dP(\nu)\sim P$
because $P$ is quasi-Palm, and that $P\sim Q$ implies $L^{*}P\sim L^{*}Q$,
which is trivial.

With regard to restricted FDs, we have the following. 
\begin{lem}
\label{lem:image-of-restricted-FD}Suppose that $P$ is a restricted
FD with extended version $\widetilde{P}$. Then for any $L\in GL_{n}(\mathbb{R})$
there is a $t_{0}>0$ such that
\begin{equation}
S_{t_{0}}^{\sqr}(LP)=(L^{*}\widetilde{P})^{\sqr}\label{eq:8}
\end{equation}
\end{lem}
\begin{proof}
Choose $t_{0}$ such that $B_{e^{-t}}(0)\subseteq L(B_{1}(0))$. Then
for any $\theta\in\mathcal{M}$ with $\theta(B_{1})>0$, the restriction
of the measures $L\theta$ and $L(\theta^{\sqr})$ to $B_{e^{-t}}(0)$
are proportional, and so $S_{t}^{\sqr}(L^{*}\theta)=S_{t}^{\sqr}(L(\theta^{\sqr}))$
for all $t\geq t_{0}$. Consider the map $\pi:\mu\mapsto S_{t_{0}}^{\sqr}(L(\mu^{\sqr}))$.
By definition, the push-forward of $\widetilde{P}$ by $\pi$ is the
left hand side of \eqref{eq:8}. On the other hand, by the previous
discussion, $\pi(\mu)=S_{t_{0}}^{\sqr}(L^{*}\mu)$. We can write $\pi=\pi_{3}\circ\pi_{2}\circ\pi_{1}$,
where $\pi_{1}(\theta)=L^{*}\theta$, $\pi_{2}(\theta)=S_{t_{0}}^{*}\theta$,
and $\pi_{3}(\theta)=\theta^{\sqr}$. Then $\pi_{1}\widetilde{P}=L^{*}\widetilde{P}$,
and since $L^{*}\widetilde{P}$ is $S_{t_{0}}^{*}$-invariant, $\pi_{2}\pi_{1}\widetilde{P}=\pi_{2}L^{*}\widetilde{P}=L^{*}\widetilde{P}$.
Finally, $\pi(\widetilde{P})=\pi_{3}(L^{*}\widetilde{P})=(L^{*}\widetilde{P})^{\sqr}$
which is the right hand side of \eqref{eq:8}, as claimed.\end{proof}
\begin{lem}
If $V<\mathbb{R}^{d}$ is an affine subspace with $0\notin V$ and
$P$ a FD (restricted or extended). Then $P(\mu\,:\,\mu(V)>0)=0$.\end{lem}
\begin{proof}
Suppose $P$ is an FD. Choose $\mu\sim P$ and suppose that $P(\mu(V)>0)>0$.
Write $V_{r}=V\cap B_{r}(0)$. Then there is an an $\varepsilon>0$
such that $P(\mu(V_{\varepsilon})>\varepsilon)>\varepsilon$. By $S^{*}$-invariance
of $P$ also $P(\mu(tV_{\varepsilon})>\varepsilon)>\varepsilon$ for
every $t>0$, whence
\[
P(\mu\,:\,\mu(\frac{1}{n}V_{\varepsilon})>\varepsilon\mbox{ for infinitely many }n\in\mathbb{N})>0
\]
But $V$ is an affine subspace with $0\notin V$, so the sets $\frac{1}{n}V$,
$n\in\mathbb{N}$, are disjoint, and so are $\frac{1}{n}V_{\varepsilon}$.
Since $\frac{1}{n}V_{\varepsilon}\subseteq B_{1/\varepsilon}(0)$,
for a $\mu$ as in the event above, $\mu(B_{1/\varepsilon}(0))\geq\sum\mu(\frac{1}{n}V_{\varepsilon})=\infty$,
contradicting $\mu\in\mathcal{M}$.\end{proof}
\begin{cor}
For $V$ as above, for $P$ a FD and every $\varepsilon>0$, there
is a $\delta>0$ such that 
\[
P(\mu\,:\,\mu^{\sqr}(V^{(\delta)})<\varepsilon)>1-\varepsilon
\]
Here $V^{(\delta)}$ is the $\delta$-neighborhood of $V$.\end{cor}
\begin{proof}
By the previous lemma, $\lim_{n\rightarrow\infty}\mu^{\sqr}(V^{(1/n)})=\mu^{\sqr}(V)=0$
for $P$-a.e $\mu$, and hence the convergence is also in probability. 
\end{proof}
We write $C(X,Y)$ for the space of continuous maps $X$$\rightarrow Y$
with the uniform metric, $d(f,g)=\sup_{x\in X}d(f(x),g(x))$. On spaces
of measures we use the metrics introduced in Section \ref{sub:Metrics-on-spaces-of-measures}
and the associated notation, e.g. writing $\nu=\mu+O(1)$ for $d(\mu,\nu)=O(1)$.
\begin{lem}
\label{lem:image-approximation}Let $L,t_{0}$ be as in the previous
lemma and $P$ a restricted FD. Then for every $\varepsilon>0$ there
is a $\delta>0$ and an open set $\mathcal{V}\subseteq\mathcal{M}^{\sqr}$,
such that 
\begin{enumerate}
\item $P(\mathcal{V})>1-\varepsilon$.
\item If $\mu,\nu\in\mathcal{V}$ with $d(\mu,\nu)<\delta$ and $L'\in C(B_{1}(0),\mathbb{R}^{d})$
with $d(L.L')<\delta$, then $d(S_{t_{0}}^{\sqr}L\mu,S_{t_{0}}^{\sqr}L'\nu)<\varepsilon$.
\end{enumerate}
\end{lem}
\begin{proof}
Write $E=B_{e^{-t_{0}}}(0)$. Consider first the operation $\nu\mapsto S_{t_{0}}^{\sqr}\nu=S_{t_{0}}(\frac{1}{\nu(E)}\nu|_{E})$
for $\nu\in\mathcal{P}(B_{2\left\Vert L\right\Vert }(0))$, where
the space $\mathcal{P}(B_{2\left\Vert L\right\Vert }(0)))$ was chosen
because $L'\mu$ belongs to it for all $\mu\in\mathcal{M}^{\sqr}$
and $L'$ sufficiently close to $L$. If this operation were weakly
continuous the lemma would be trivial. Unfortunately, there are discontinuities
when $\nu(\partial E)>0$. Thus we proceed more cautiously.

Observe that $\lip(f\circ S_{t_{0}})=e^{t_{0}}\lip(f)$, so from our
definition of the metric on measures, $d(S_{t_{0}}\nu_{1},S_{t_{0}}\nu_{2})=e^{t_{0}}\cdot d(\nu_{1},\nu_{2})$
for $\nu_{1},\nu_{2}\in\mathcal{P}(B_{2\left\Vert L\right\Vert }(0))$.
Our aim is to find an open set $\mathcal{V}'\subseteq\mathcal{P}(B_{2\left\Vert L\right\Vert }(0))$
and $\delta'>0$ such that
\begin{enumerate}
\item [(1')] $P(\mu\in\mathcal{M}\,:\, L\mu\in\mathcal{V}')>1-\varepsilon$
\item [(2')] If $\mu\in\mathcal{V}'$ and $\nu\in\mathcal{P}(B_{2\left\Vert L\right\Vert }(0))$
with $d(\mu,\nu)<\delta'$, then $d(\mu_{E},\nu_{E})<e^{-t_{0}}\varepsilon$.
\end{enumerate}
Assuming this, let $\mathcal{V}=L^{-1}\mathcal{V}'\cap\mathcal{M}^{\sqr}$.
Then (1') implies (1). Also, since $\nu\mapsto L'\nu$ is weakly continuous,
there is a $\delta$ such that whenever $d(L,L')<\delta$ and $\mu,\nu\in\mathcal{V}$
satisfy $d(\mu,\nu)<\delta$, we have $d(L\mu,L'\nu)<\delta'$. By
(2') and the previous discussion this is (2).

We turn to the construction of $\mathcal{V}'$ and $\delta'$. For
measurable $f:\mathbb{R}^{d}\rightarrow[0,1]$ and $\mu\in\mathcal{P}(B_{2\left\Vert L\right\Vert }(0))$
with $\int fd\mu>0$, write $\mu_{f}$ for the measure $\mu_{f}(A)=\int_{A}fd\mu/\int fd\mu\in\mathcal{P}(B_{2\left\Vert L\right\Vert }(0))$.
When $f$ is continuous, $\nu\mapsto\nu_{f}$ is weakly continuous. 

Fix a parameter $\rho>0$ and choose a continuous function $f$ with
$1_{(1-\rho)E}\leq f\leq1_{E}$. For $\mu\in\mathcal{P}(B_{2\left\Vert L\right\Vert }(0))$
with $\mu(E)>0$, 
\[
\mu_{f}=\alpha\cdot\mu_{1_{(1-\rho)E}}+(1-\alpha)\cdot\mu_{f-1_{(1-\rho)E}}
\]
for $\alpha=\frac{\mu((1-\delta)E)}{\mu(E)}$. Similarly
\[
\mu_{E}=\alpha\cdot\mu_{1_{(1-\rho)E}}+(1-\alpha)\cdot\mu_{1_{E\setminus(1-\rho)E}}
\]
Thus $\mu_{f}-\mu_{E}=(1-\alpha)\theta$, where $\theta=\mu_{f-1_{(1-\rho)E}}-\mu_{1_{E\setminus(1-\rho)E}}$,
so 
\[
d(\mu_{E},\mu_{f})\leq1-\alpha
\]

For $\nu\in\mathcal{P}(B_{2\left\Vert L\right\Vert }(0))$ close enough
to $\mu$ we will certainly have $\int fd\nu>0$, so $\nu_{f}$ is
well defined, and, writing $\beta=\frac{\nu((1-\rho)E)}{\nu(E)}$,
\[
d(\mu_{E},\nu_{E})\leq d(\mu_{E},\mu_{f})+d(\mu_{f},\nu_{f})+d(\nu_{f},\nu_{E})\leq(1-\alpha)+(1-\beta)+d(\mu_{f},\nu_{f})
\]
If $\frac{\mu((1-2\rho)E)}{\mu((1+\rho)E)}$ is sufficiently close
to $1$ then, trivially, $\alpha$ is close to $1$, and if $d(\mu,\nu)$
is small enough in a manner depending on $\rho$, then also $\beta$
is close to $1$. Also, since $\mu\mapsto\mu_{f}$ is continuous,
the last term on the right in the equation above can be made arbitrarily
small by making $d(\mu,\nu)$ small. We conclude that, given $\rho$,
there is a $\delta'>0$ such that $d(\mu_{E},\nu_{E})<e^{-t_{0}}\varepsilon$
as long as $d(\mu,\nu)<\delta'$ and $\mu$ belongs to the set 
\[
\mathcal{V}'=\{\mu\in\mathcal{P}(B_{2\left\Vert L\right\Vert }(0))\,:\,\frac{\mu((1-2\rho)E^{\circ})}{\mu((1+\rho)E)}>1-\delta'\}
\]
(here $E^{\circ}$ is the interior of $E$, and appears above to ensure
that the set is open). To complete the proof we must show that there
exists a $\rho>0$ such that $P(L\mu\in\mathcal{V}')>1-\varepsilon$.
This is immediate from Corollary, which implies that for $P$-a.e.
$\mu$, 
\[
\frac{\mu((1-2\rho)L^{-1}E^{\circ})}{\mu((1+\rho)L^{-1}E)}\geq1-\frac{\mu((\partial L^{-1}E)^{(2\rho)})}{\mu((1+\rho)L^{-1}E)}\rightarrow0\qquad\mbox{as }\rho\rightarrow0
\]
because $\partial L^{-1}E$ is a finite union of affine subspaces
that do not contain $0$.
\end{proof}
We now turn to Proposition \ref{pro:smooth-images-of-usms}. Let $\mu$
be a probability measure on an open set $U\subseteq\mathbb{R}^{d}$
and $f:U\rightarrow V\subseteq\mathbb{R}^{d}$ a $C^{1}$-diffeomorphism
to another open set $V$. Suppose that $\left\langle \mu\right\rangle _{x,T_{i}}^{\sqr}\rightarrow P$
for a FD $P$ with extended version $\widetilde{P}$, let $\nu=f\mu$
and $y=fx$. We claim that $\left\langle \nu\right\rangle _{y,T_{i}}^{\sqr}\rightarrow L^{\sqr}\widetilde{P}$,
where $L=Df(x)$. Without loss of generality we can assume that $x=y=0$
(if not, replace $f(z)$ by $f(z+x)-y$). 

Choose $t_{0}>0$ as in Lemma \ref{lem:image-of-restricted-FD}. The
linearization of $f$ at $x=0$ is 
\[
f(w)=L(w)+o(\left\Vert w\right\Vert )
\]
For large enough $t$, and $w\in B_{e^{-t}}(x)$, this implies that
\[
B_{e^{-(t+t_{0})}}(0)\subseteq f(B_{e^{-t}}(0))
\]
and therefore
\begin{eqnarray*}
S_{t+t_{0}}^{\sqr}\nu & = & S_{t+t_{0}}^{\sqr}(f(\mu_{B_{e^{-t}(x)}}))\\
 & = & S_{t+t_{0}}^{\sqr}((L+o(e^{-t}))(\mu_{B_{e^{-t}(x)}})\\
 & = & S_{t_{0}}^{\sqr}((L+o(1))\mu_{x,t}^{\sqr})\qquad\mbox{as }t\rightarrow\infty
\end{eqnarray*}
Let $\varepsilon>0$ and let $\mathcal{U},\mathcal{V},\delta$ be
as in Lemma \ref{lem:image-approximation}. The map $L+o(1)$ in the
last equation belongs to $\mathcal{U}$ for all large enough $t$,
and if for such a $t$ we have $\mu_{x,t}^{\sqr}\in\mathcal{V}$ then
\[
d(S_{t+t_{0}}^{\sqr}\nu,S_{t_{0}}^{\sqr}L(\mu_{x,t}^{\sqr}))<\varepsilon
\]
Since $\mathcal{V}$ is open, $P(\mathcal{V})>1-\varepsilon$ and
$\left\langle \mu\right\rangle _{x,T_{i}}^{\sqr}\rightarrow P$, for
large enough $i$ we will have
\begin{equation}
\frac{1}{T_{i}}\lambda\left(\left\{ t\in[0,T_{i}]\,:\,\mu_{x,t}\in\mathcal{V}\right\} \right)>1-\varepsilon\label{eq:14}
\end{equation}
Therefore, as $i\rightarrow\infty$,
\begin{alignat*}{2}
\left\langle \nu\right\rangle _{y,T_{i}}^{\sqr}\quad & =\quad\frac{1}{T_{i}}\int_{0}^{T_{i}}\delta_{\nu_{y,t}^{\sqr}}\, dt\\
 & =\quad\frac{1}{T_{i}}\int_{t_{0}}^{T_{i}+t_{0}}\delta_{\nu_{y,t}^{\sqr}}\, dt+o(1)\\
 & =\quad\frac{1}{T_{i}}\int_{0}^{T_{i}}\delta_{\nu_{y,t+t_{0}}^{\sqr}}\, dt+o(1)\\
 & =\quad\frac{1}{T_{i}}\int_{0}^{T_{i}}\delta_{S_{t_{0}}^{\sqr}L\mu_{x,t}^{\sqr}}\, dt+O(\varepsilon) &  & \mbox{by \eqref{eq:14}}\\
 & =\quad S_{t_{0}}^{\sqr}L\left(\frac{1}{T_{i}}\int_{0}^{T}\delta_{\mu_{x,t}^{\sqr}}\, dt\right)+O(\varepsilon) & \qquad & \mbox{by linearity}\\
 & =\quad S_{t_{0}}^{\sqr}L\left(\left\langle \mu\right\rangle _{x,T_{i}}^{\sqr}\right)+O(\varepsilon)
\end{alignat*}
Write $P_{i}=\left\langle \mu\right\rangle _{x,T_{i}}^{\sqr}$. Since
$P_{i}\rightarrow P$ as $i\rightarrow\infty$ and $\mathcal{V}$
is open, we find that $\liminf P_{i}(\mathcal{V})\geq P(\mathcal{V})>1-\varepsilon$.
Thus, writing 
\begin{eqnarray*}
P_{i} & = & (P_{i})|_{\mathcal{V}}+(P_{i})|_{\mathcal{M}^{\sqr}\setminus\mathcal{V}}\\
 & = & (1-O(\varepsilon))(P_{i})_{\mathcal{V}}+O(\varepsilon)(P_{i})_{\mathcal{M}^{\sqr}\setminus\mathcal{V}}
\end{eqnarray*}
we find that 
\[
\left\langle \nu\right\rangle _{y,T_{i}}^{\sqr}=S_{t_{0}}^{\sqr}L((P_{i})_{\mathcal{V}})+O(\varepsilon)\qquad\mbox{as }i\rightarrow\infty
\]
Splitting $P$ in the same way as $P=(1-O(\varepsilon))P_{\mathcal{V}}+O(\varepsilon)P_{\mathcal{M}^{\sqr}\setminus\mathcal{V}}$,
we also have
\[
S_{t_{0}}^{\sqr}L(P_{\mathcal{V}})=S_{t_{0}}^{\sqr}L(P)+O(\varepsilon)
\]
From the definition of $\mathcal{V}$ ((2) in Lemma \ref{lem:image-approximation}),
upon pushing $(P_{i})_{\mathcal{V}}$ forward through $\theta\mapsto S_{t_{0}}^{\sqr}L\theta$
and applying Lemma \ref{lem:modulous-of-continuity-for-measures},
we find that $S_{t_{0}}^{\sqr}L((P_{i})_{\mathcal{V}})=S_{t_{0}}^{\sqr}L(P_{\mathcal{V}})+O(\varepsilon)$,
hence 
\begin{eqnarray*}
\left\langle \nu\right\rangle _{y,T_{i}}^{\sqr} & = & S_{t_{0}}^{\sqr}LP+O(\varepsilon)\qquad\mbox{as }i\rightarrow\infty
\end{eqnarray*}
But $\varepsilon$ was arbitrary, so this means that $\left\langle \nu\right\rangle _{y,T_{i}}^{\sqr}\rightarrow S_{t_{0}}^{\sqr}LP$.
By Lemma \ref{lem:image-of-restricted-FD}, $S_{t_{0}}^{\sqr}LP=L^{*}\widetilde{P}$,
and we are done.

\subsection{\label{sub:Homogeneous-measures}Homogeneous measures}
\begin{proof}
[Proof of Proposition \ref{pro:homogeneous-implies-usm}] Let $\mu$
be a homogeneous measure. Then we may choose a point $x$ which is
homogeneous (Definition \ref{def:homogeneous-measures}) and such
that every accumulation point of the scenery distributions at $x$
are FDs (Theorem \ref{thm:limit-distributions-are-FDs}). Let $P$
be such an accumulation point (we think of $P$ as a restricted FD)
and choose a $P$-typical $\nu$. Then $\nu$ is an accumulation point
of $\mu_{x,t}$ as $t\rightarrow\infty$, and so by homogeneity $\mu\ll T_{B}^{\sqr}\nu$
for some ball $B$. Since $\nu$ is a USM generating the ergodic component
$P_{\nu}$ of $P$ we find that $\mu$ is a USM. In particular $P=P_{\nu}$
and $\mu$-a.s is independent of $x$. 

It remains to show that a $P$-typical $\nu$ is homogeneous. Let
$\widetilde{P}$ denote the extended version of $P$. By our previous
remark a $P$-typical $\nu$ contains a positive-measure ball $E\subseteq B_{1}$
on which $\nu$ is equivalent, after re-scaling, to $\mu$. It follows
from $S^{*}$-invariance of $\widetilde{P}$ that for $\widetilde{P}$-typical
$\widetilde{\nu}$ there is a ball $E\supseteq B_{1}$ on which $\widetilde{\nu}$
is equivalent, after re-scaling, to $\mu$. Since $\widetilde{\nu}^{\sqr}$
is equivalent to $\widetilde{\nu}$ on $B_{1}$, we find that $\widetilde{\nu}^{\sqr}$
is $\widetilde{P}$-a.s. equivalent to $\mu$ on some ball $A$, so
the same is true for $P$-typical $\nu$ instead of $\widetilde{P}$-typical
$\widetilde{\nu}$. Now, if $\nu$ is fixed and $\nu'$ is an accumulation
point of the scenery of $\nu$ at a $\nu$-typical point $y$, we
know by homogeneity of $\mu$ that $\mu\ll T_{C}\nu'$ for a ball
$C$ and hence $\nu\ll T_{A_{\nu}}T_{C}\nu'$. This shows that $\nu$
is homogeneous.
\end{proof}

\begin{proof}
[Proof of Theorem \ref{thm:properties-of-homogeneous-measures}]Let
$\mu$ be homogeneous and $P$ the generated EFD. For each $\pi\in\Pi_{d,k}$,
for $P$-a.e. $\mu'$ we have that $\pi\mu'$ is exact dimensional
and $\dim\pi\mu'=E_{\pi}(P)$. Since for $P$-a.e. $\mu'$ we have
$\mu\ll A\mu'$ for some homothety $A$, we conclude that $\pi\mu$
is exact dimensional and $\dim\pi\mu=E_{\pi}(P)$. Since $\pi\mapsto E_{\pi}(P)$
is lower semi-continuous by Theorem \ref{thm:semi-continuity-for-CP-processes}
(and the fact that $P$ is the centering of an ergodic CP-process
by Theorem \ref{thm:FDs-come-from-CPs}), the map $\pi\mapsto\dim\pi\mu$
is as well. The same argument applies to conditional measures: since
$P$-a.e. $\mu'$ has a.s. exact dimensional fibers and whose dimensions
satisfy equation \eqref{eq:dimension-conservation}, and $\mu\ll A\mu'$,
the same holds for $\mu$.
\end{proof}

\subsection{\label{sub:EFDs-invariant-under-linear transformations}EFDs invariant
under groups of linear transformations}

Next we turn to Proposition \ref{pro:measures-with-additional-invariance}.
In fact this is an immediate consequence of lower semi-continuity
of the function $E_{P}(\cdot)$, and the last part uses also Marstrand's
theorem (Theorem \ref{thm:Hunt-Kaloshin}).

This proposition sheds light on the main results of \cite{HochmanShmerkin09}.
Let us recall one of these:
\begin{thm}
[\cite{HochmanShmerkin09}]\label{thm:HS-2}Let $f_{1},\ldots,f_{r}$
be contracting similarities of $\mathbb{R}^{d}$ satisfying the strong
separation condition. Suppose the orthogonal parts of $f_{i}$ generate
a dense subgroup of the orthogonal group. Let $X$ be the attractor
of the IFS $\{f_{i}\}$ and $\mu$ the $\dim X$-dimensional Hausdorff
measure on $X$. Then for every $\pi\in\Pi_{d,k}$ we have 
\[
\dim\pi\mu=\min\{k,\dim\mu\}
\]

\end{thm}
Let us see how this follows from our methods. The measure $\mu$ generates
an EFD $P$, and since $f_{i}\mu\ll\mu$, we see that this EFD will
be invariant under the orthogonal part of $f_{i}$ for each $i=1,\ldots,r$.
The group $H$ generated by the orthogonal parts of the $f_{i}$ is
by assumption dense in the orthogonal group. Thus the orbit of $GL_{k}(\mathbb{R})\times H$
on $\Pi_{d,k}$ is the entire space $\Pi_{d,k}$ and Theorem \ref{thm:HS-2}
follows from Proposition \ref{pro:measures-with-additional-invariance}.
\begin{proof}
[Proof of Theorem \ref{thm:dissonance-with-Ahlfors}] Let $\mu$ be
a self-similar measure on $\mathbb{R}$ defined by an IFS satisfying
strong separation and such that the contraction ratios generate a
dense multiplicative subgroup of $\mathbb{R}^{+}$, and with dimension
$\alpha$. Let $\nu$ be an Ahlfors-regular measure of dimension $\beta$.
We wish to show that $\dim\mu*\nu=\min\{1,\alpha+\beta\}$. We assume
$\alpha,\beta>0$, since otherwise the statement is trivial.

Write $\gamma=\min\{1,\alpha+\beta\}$ and let $\theta=\mu\times\nu$.
By Theorem \ref{thm:projection-of-SMs}, we must show that for $\theta$-a.e.
$z=(x,y)$, and for any distribution $P$ generated by $\theta$ at
$z$, we have $E_{\pi}(P)=\gamma$ for $\pi(x,y)=x+y$. By definition
of $E_{\pi}(P)$ it suffices to show that $E_{\pi}(P')=\gamma$ for
a.e. ergodic component $Q$ of $P$. We fix such a component $Q$.

Let $\theta'$ be a $Q$-typical measure. Then $\theta'$ arises as
a weak limit of measures $\theta_{z,t}=\mu_{x,t}\times\nu_{x,t}$
for a sequence of $t$'s tending to $\infty$, and hence $\theta'=\mu'\times\nu'$
for $\mu'$, $\nu'$ that are accumulation points of $\mu_{x,t}$
and $\nu_{x,t}$ respectively. Notice that $\mu',\nu'$ are equal
to the conditional measures of $\theta'$ on the $x$ and $y$ axes,
respectively, and so we may assume that they are typical measures
for some EFDs and hence exact dimensional.

Since $z=(x,y)$ was $\theta$-typical, $x$ is $\mu$-typical, so
by homogeneity we can assume that $\mu\ll T_{I}^{\sqr}\mu'$ for some
interval $I$ and in particular $\dim\mu'=\dim\mu'_{I}=\dim\mu=\alpha$
(the first equality is due to exact dimensionality of $\mu'$). From
the definition of Ahlfors-regularity, it is clear that $\nu_{y,t}$
is Ahlfors-regular with the same constants, and hence $\nu'$ is too.
In particular $\dim\nu'=\beta$. Thus $\dim\theta'=\alpha+\beta$,
and this is also $\dim Q$ (this is the only place where Ahlfors regularity
is used. Without it we could conclude that $\dim P=\beta$ but not
that this passes to its ergodic components).

Since $\theta'$ is $Q$-typical, a.e. point generates $Q$. Fix an
interval $J$ such that $\nu'(J)>0$ and $|J|=|I|$. Then, since $\mu\ll T_{I}^{\sqr}\mu'$,
we have $\mu\times T_{J}^{\sqr}\nu'\ll(T_{I}\times T_{J})^{\sqr}\theta'$.
Since $T_{I}\times T_{J}=T_{I\times J}$ is a homothety, $\theta''=\mu\times T_{J}^{\sqr}\nu'\ll(T_{I}\times T_{J})^{\sqr}\theta'$
generates $Q$ at a.e. point. Write $\nu''=T_{J}^{\sqr}\nu$.

Now, suppose that $U$ is an interval with $T_{U}^{\sqr}\mu=\mu$.
Then $(T_{U}\times\id)^{\sqr}\theta''=\theta''$, and so by Proposition
\ref{pro:smooth-images-of-usms}, the (extended) distribution $Q$
generated at a $\theta''$-typical $z\in U\times[-1,1]$ is invariant
under the map $(T_{U}\times\id)^{*}$. By our assumption on $\mu$,
the maps $T_{U}^{*}$ for $U$ as above generate a dense subgroup
of $(S_{t}^{*})_{t\in\mathbb{R}}$, so that $Q$ is invariant under
$((S_{t}\times\id)^{*})_{t\in\mathbb{R}}$. Since it is also invariant
under $(S_{t}^{*}\times S_{t}^{*})_{t\in\mathbb{R}}$ (here the maps
$S_{t}^{*}$ act partially on $\mathcal{M}(\mathbb{R})$), for all
$t$, we find that $Q$ is invariant under $((S_{t_{1}}\times S_{t_{2}})^{*})_{t_{1},t_{2}\in\mathbb{R}}$.
It then follows from Proposition \ref{pro:measures-with-additional-invariance}
that 
\[
E_{\pi}(Q)=\min\{1,\dim Q\}=\min\{1,\alpha+\beta\}=\gamma
\]
This completes the proof.
\end{proof}

\section{Constructions}

\subsection{\label{sub:non-conservation-example}Counterexamples to dimension
conservation}

In this section we construct a self-similar set $A\subseteq\mathbb{R}^{2}$,
defined by an IFS with uniform contraction, no rotations and strong
separation; and a $C^{\infty}$ diffeomorphism $f_{0}:\mathbb{R}^{2}\rightarrow\mathbb{R}^{2}$
such that (i) $Df|_{A}=\id$, (ii) the coordinate projection $\pi(x,y)=x$
is injective on the set $E=f_{0}(A)$, and (iii) $\dim\pi E<\min\{1,\dim A\}$.
Taking $\mu$ to be the natural (normalized $\dim A$-dimensional
Hausdorff) measure on $A$, we know that $\mu$ is a USM for an EFD
$P$ of the same dimension as $\mu$, namely $\dim A$. Because of
(i) and Proposition \ref{pro:smooth-images-of-usms}, the measure
$\nu=f_{0}\mu$ is also a USM for $P$. Also, $\nu$ is supported
on $E$ and $\pi\nu$ is supported on $\pi E$ so $\dim\pi\nu\leq\min\{1,\dim E\}$.
Since the conditional measure $\nu_{[x]_{\pi}}$ are $\nu$-a.s. supported
on singletons, we find that for $\nu$-a.e. $x$,
\[
\dim\pi\nu+\dim\nu_{[x]_{\pi}}<\min\{1,\dim A\}+0<\min\{1,\dim\nu\}
\]
This shows that dimension conservation fails for the USM $\nu$. Similarly,
taking $f=\pi\circ f_{0}$ we find that for every $y\in f(A)$,
\[
\dim f(A)+\dim(A\cap f^{-1}(y))<\dim A
\]
The set $A$ is homogeneous in the sense of Furstenberg \cite{Furstenberg08},
and this shows that dimension conservation for homogeneous sets, which
Furstenberg proves for linear maps, does not hold in general for $C^{\infty}$
maps.

We turn to the construction. The set $A$ will be the self-product
$A=C\times C$ of the Cantor set 
\[
C=\{\sum a_{n}10^{-n}\,:\, a_{n}\in\{1,8\}\}
\]
(Our choices of base and digits is essentially arbitrary as will be
seen in the construction). Clearly $A=\bigcup_{u\in U}(\frac{1}{10}A+u)$
where $U=\{(\frac{a}{10},\frac{b}{10})\,:\, a,b\in\{1,8\}\}\subseteq\mathbb{R}^{2}$,
and so it satisfies the properties claimed above.

\begin{figure}
\input{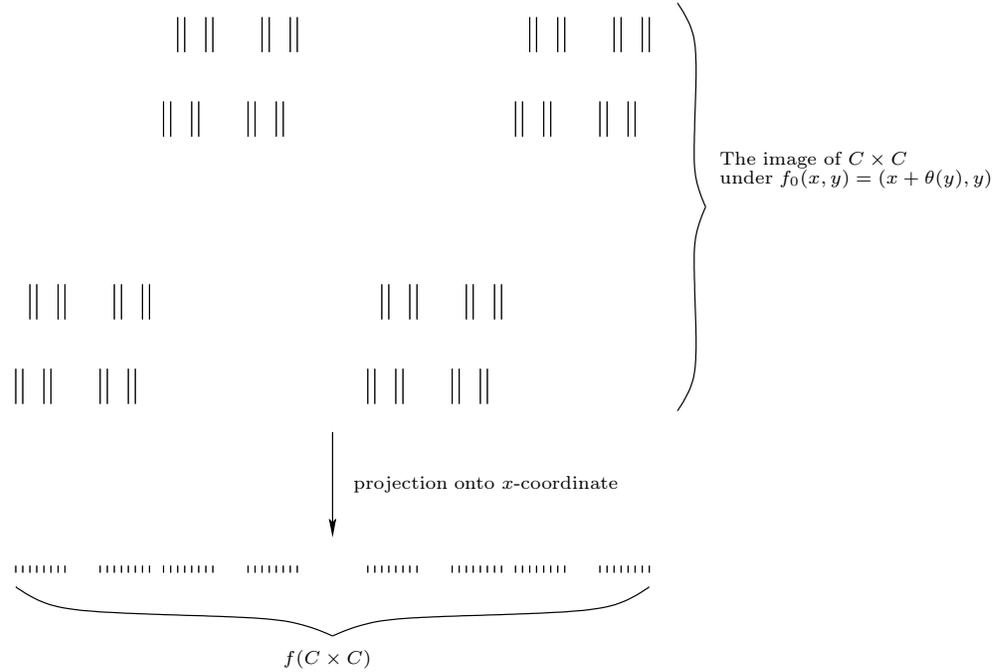}

\caption{Schematic representation of $f$ as composition of a plane map and
a projection.\label{fig:smooth-map-1}}
\end{figure}

It will be convenient to expand points $x\in[0,1]$ in base $10$,
written $x=0.x_{1}x_{2}\ldots$; note that points in $C$ have a unique
expansion. We associate to $x=0.x_{1}x_{2}\ldots\in C$ a sequence
$(c_{n}(x))_{n=1}^{\infty}\in\{0,1\}^{\mathbb{N}}$ defined by 
\[
c_{n}(x)=1_{\{x_{n}=8\}}
\]
We first define $f_{0}$ on $C\times C$ and later extend it to $[0,1]^{2}$.
Fix a rapidly increasing sequence of integers $n_{k}\rightarrow\infty$
(we will specify its properties later) and define
\[
\theta(y)=\sum_{k=1}^{\infty}c_{k}(y)\cdot10^{-n_{k}}
\]
Finally, set
\[
f_{0}(x,y)=(x+\theta(y),y)
\]
see Figure \ref{fig:smooth-map-1}.

\begin{lem}
For $(x,y)\in C\times C$ the point $x+\theta(y)$ determines both
$x$ and $y$, thus $f=\pi f_{0}$ is injective on \textup{$C\times C$
(equivalently, $\pi$ is injective on $f_{0}(C\times C)$).}\end{lem}
\begin{proof}
Write $z=x+\theta(y)=0.z_{1}z_{2}\ldots$. The decimal digits of $x$
are all $1$ or $8$, and the digits of $\theta(y)$ are $0$ or $1$,
so $z_{i}\in\{1,2,8,9\}$ for all $i$ (note that we have arranged
that in the addition $x+\theta(y)$ no ``carries'' occur). Furthermore
$z_{i}\in\{2,9\}$ can only occur if $i=n_{k}$ for some $k$, in
which case $\theta_{k}(y)=1$, equivalently $y_{k}=8$. Thus, setting
\[
x_{i}=\left\{ \begin{array}{cc}
1 & z_{i}\in\{1,2\}\\
8 & z_{i}\in\{8,9\}
\end{array}\right.
\]
we have $x=0.x_{1}x_{2}\ldots$, and setting 
\[
y_{k}=\left\{ \begin{array}{cc}
0 & z_{n_{k}}\in\{1,8\}\\
1 & z_{n_{k}}\in\{2,9\}
\end{array}\right.
\]
we have $y=0.y_{1}y_{2}\ldots$, so $z$ determines $x,y$.\end{proof}
\begin{lem}
If $n_{k}\rightarrow\infty$ rapidly enough then $\dim\pi f_{0}(C\times C)=\dim C$.\end{lem}
\begin{proof}
Note that ~$\pi f_{0}(C\times C)=\{x+\theta(y)\,:\,(x,y)\in C\times C\}$.
As we saw in the previous lemma, the digits of $x+\theta(y)$ can
be $1,8$ at indices $i\neq n_{k}$, and are among $1,2,8,9$ at the
indices $n_{1},n_{2},\ldots$. Let $k(i)$ denote the unique $k$
such that $n_{k}\leq i<n_{k+1}$. Then the number $N_{i}$ of sequences
$a_{1}\ldots a_{i}$ that occur as the initial  segments of expansions
of $x+\theta(y)$ is (at most, but in fact equal to) $2^{i-k(i)}\cdot4^{k(i)}$,
so the box dimension of $\pi f_{0}(C\times C)$ is bounded by
\[
\limsup_{i\rightarrow\infty}\frac{\log(2^{i-k(i)}4^{k(i)})}{\log10^{-i}}=\limsup_{i\rightarrow\infty}\frac{(i-k(i))+2k(i)}{i\log10}=\frac{1}{\log10}+\limsup_{i\rightarrow\infty}\frac{k(i)}{i\log10}
\]
This shows that if $n_{k}$ grows super-linearly (i.e. $k(i)/i\rightarrow0$)
then $\bdim\pi f_{0}(C\times C)\leq\frac{1}{\log10}=\dim C$. On the
other hand, $\pi f_{0}(C\times C)\supseteq\pi f_{0}(C\times\{0.1111\ldots\})=C$,
so $\dim\pi f_{0}(C\times C)\geq\dim C$, and we have equality.
\end{proof}
Finally, the functions $c_{k}(\cdot)$ may be extended as follows.
First, set $c_{k}(x)=0$ on the interval consisting of points $x=0.x_{1}x_{2},\ldots$
with $x_{k}=1$, and $c_{k}(x)=1$ on the interval where $x_{k}=8$.
Now $c_{k}$ is defined on two disjoint intervals, on each of which
it is constant, so there is a $C^{\infty}$ extension to $[0,1]$.
Having fixed the extension of $c_{k}$, choose $n_{k}$ large enough
that the first $k$ derivatives of $10^{-n_{k}}c_{k}$ are less than
$2^{-k}$ everywhere on $[0,1]$; clearly this is satisfied for all
large enough $n_{k}$. It follows that $\theta(\cdot)=\sum c_{k}(\cdot)10^{-n_{k}}$
is $C^{\infty}$ on $[0,1]$, and therefore so will $f_{0}$. The
fact that $Df_{0}=\id$on $C\times C$ follows directly from the fact
that $c_{k}'(x)=0$ for $x\in C$.

\subsection{\label{sub:recipe-for-constructions}A framework for constructing
USMs}

In this section we present a procedure for constructing a USM from
a sequence of USMs by ``splicing'' them together at different scales.
It will be convenient to adopt the metrics on measures and the notation
of Section \ref{sub:Metrics-on-spaces-of-measures}.

The basic operation will be the following. Let $\mu,\nu\in\mathcal{M}^{\sqr}$.
For an integer $n$, the $(\nu,n)$-\emph{discretization of }$\mu$
is the measure $\dsc(\mu,\nu,n)\in\mathcal{M}^{\sqr}$ given by 
\[
\dsc(\mu,\nu,n)=\sum_{E\in\mathcal{E}}\mu(E)\cdot T_{E}^{-1}\nu
\]
Note that $T_{E}^{-1}\nu$ is just a copy of $\nu$ scaled down from
$[-1,1]^{d}$ to $E$. Thus, $\dsc(\mu,\nu,n)$ is defined on each
cell $E\in\mathcal{D}_{n}$ to be a scaled copy of $\nu$ with total
mass equal to $\mu(E)$. 

Fix an integer base $b\geq2$ and, as before, write $\widetilde{M}=M_{b}^{\sqr}$.
Let $\mu'=\dsc(\mu,\nu,b^{N})$. We make a number of observations.
First, by definition $\mu(E)=\mu'(E)$ for all $E\in\mathcal{D}_{b^{N}}$,
so by \eqref{eq:23}, 
\begin{equation}
d(\mu',\mu)\leq c\cdot b^{-N}\label{eq:dsc-distance}
\end{equation}

Second, for $\mu'$-a.e. $x$ we have $\mu(\mathcal{D}_{b^{n}}(x))=\mu'(\mathcal{D}_{b^{n}}(x))>0$
for $n\leq N$ and $\mu'(\mathcal{D}_{b^{n}}(x))>0$ for $n>N$. Thus,
the $b$-adic scenery $(\mu'_{n},x_{n})=\widetilde{M}^{n}(\mu',x)$
of $\mu'$ is well defined for all $n$, and $(\mu_{n},x_{n})=\widetilde{M}^{n}(\mu,x)$
is well defined for $n\leq N$. 

Next, fixing $\mu'$s-typical $x$, for $n\leq N$ it is clear that
$\mu'_{n}=\dsc(\mu_{n},\nu,b^{N-n})$, so by \eqref{eq:dsc-distance},
\[
d(\mu_{n},\mu'_{n})<cb^{_{-(N-n)}}
\]
and by \eqref{eq:distance-between averages}, for $0\leq k\leq N$,
\begin{eqnarray}
d(\left\langle \mu,x\right\rangle _{k},\left\langle \mu',x\right\rangle _{k}) & \leq & \frac{1}{k}\sum_{n=0}^{k-1}d(\mu_{n},\mu'_{n})\nonumber \\
 & = & O_{b}(\frac{1}{k})\label{eq:17}
\end{eqnarray}
On the other hand, by the definition of the discretization operation,
for $n>N$ we have $(\mu'_{n},x_{n})=M^{n-N}(\nu,T_{\mathcal{D}_{b^{N}}(x)}x)$.
Combining these two observations, for any $L\geq N$,
\begin{equation}
\left\langle \mu',x\right\rangle _{L}=\frac{N}{L}\left\langle \mu,x\right\rangle _{N}+\frac{L-N}{L}\left\langle \nu,x\right\rangle _{L-N}+O_{b}(\frac{1}{L})\label{eq:dsc-time-average-distance}
\end{equation}
Note that the error term on the right hand side is uniform in $\mu,\nu$. 

We make two more observations. First, for measures $\mu,\nu,\theta$
and $M>N$ we have the identity
\begin{equation}
\dsc(\mu,\dsc(\nu,\theta,b^{M}),b^{N})=\dsc(\dsc(\mu,\nu,b^{N}),\theta,b^{M+N})\label{eq:dsc-recursion}
\end{equation}
Second, if $\mu'$ is a weak limit of $\mu_{n}=$$\dsc(\mu,\nu_{n},b^{N_{n}})$,
then, if $N_{n}\rightarrow\infty$ we have $\mu'=\mu$; otherwise
$\mu'=\dsc(\mu,\nu,b^{N})$ for some $N$ and $\nu$ an accumulation
point of the $\nu_{n}$. We omit the verification of these two easy
facts.

Our basic construction is the following. Begin with a sequence $\nu_{1},\nu_{2},\ldots$
of measures such that $\nu_{i}$ generate a CP-distribution $P_{i}$
long $b$-adic cubes. Choose an increasing sequence $N_{n}$ and define
\begin{eqnarray*}
\mu_{1} & = & \nu_{1}\\
\mu_{n+1} & = & \dsc(\mu_{n},\nu_{n+1},b^{N_{n}})
\end{eqnarray*}
Now, it is easy to see that $\mu_{n}\rightarrow\mu$ weakly for some
measure $\mu$, since by \eqref{eq:dsc-distance},
\begin{eqnarray*}
d(\mu_{n},\mu_{n+1}) & = & d(\mu_{n},\dsc(\mu_{n},\nu_{n+1},b^{N_{n}}))\\
 & \leq & c\cdot b^{-N_{n}}
\end{eqnarray*}
which is summable, and $\mathcal{M}^{\sqr}$ is complete. Also, for
all $n$, 
\begin{equation}
\mu=\dsc(\mu_{n},\mu^{(n)},b^{N_{n}})\label{eq:19}
\end{equation}
where $\mu^{(n)}$ is the measure constructed in the same way as $\mu$
using from the data $\{\nu_{i}\}_{i=n+1}^{\infty}$ and $\{N_{i}-N_{n}\}_{i=n+1}^{\infty}$. 

Fix $\mu$-typical $x$ and let $x_{n}$ denote the second coordinate
of $M^{n}(\mu,x)$. Then, iterating \eqref{eq:dsc-time-average-distance}
using \eqref{eq:19},
\begin{eqnarray}
\left\langle \mu,x\right\rangle _{k} & = & \sum_{i=1}^{n-1}\frac{N_{i}-N_{i-1}}{k}\left(\left\langle \nu_{i},x_{N_{i-1}}\right\rangle _{N_{i}-N_{i-1}}+O_{b}(\frac{1}{N_{i}-N_{i-1}})\right)\nonumber \\
 &  & \quad+\;\frac{k-N_{n-1}}{k}\left(\left\langle \nu_{n},x_{N_{n-1}}\right\rangle _{k-N_{n-1}}+O_{b}(\frac{1}{k})\right)\nonumber \\
 & = & \sum_{i=1}^{n-1}\frac{N_{i}-N_{i-1}}{k}\left\langle \nu_{i},x_{N_{i-1}}\right\rangle _{N_{i}-N_{i-1}}\nonumber \\
 &  & \quad+\;\frac{k-N_{n-1}}{k}\left\langle \nu_{n},x_{N_{n-1}}\right\rangle _{k-N_{n-1}}+O_{b}(\frac{n}{k})\label{eq:dsc-iteratived-distance}
\end{eqnarray}
With a little care, this estimate is enough to control what distributions
$\mu$ generates along $b$-adic sceneries:
\begin{prop}
\label{prop:splicing-measures}Suppose that $\nu_{n}\in\mathcal{M}^{\sqr}$
is such that $\left\langle \nu_{n},y\right\rangle _{N}\rightarrow P_{n}$
as $N\rightarrow\infty$ for $\nu_{n}$-a.e. $y$. Suppose that $0\leq\delta_{n}\rightarrow0$
and $N_{n}$ are integers such that 
\begin{enumerate}
\item $N_{n}-N_{n-1}\rightarrow\infty$ as $n\rightarrow\infty$
\item $N_{n}\geq(1+\delta_{n})N_{n-1}$ 
\item For all $\varepsilon>0$, setting
\[
q_{n}(\varepsilon)=\nu_{n}\left(\left\{ y\,:\,\exists k\geq\delta_{n}N_{n-1}\quad d(\left\langle \nu_{n},y\right\rangle _{k},P_{n})>\varepsilon\right\} \right)
\]
we have $\sum q_{n}(\varepsilon)<\infty$
\end{enumerate}
Let $\mu$ be constructed from $\nu_{n},N_{n}$ as above. Then for
$\mu$-a.e. $x$, writing $n(k)$ for the unique $n$ such that $k\in(N_{n-1},N_{n}]$,
we have
\begin{equation}
\left\langle \mu,x\right\rangle _{k}=\sum_{i=1}^{n(k)-1}\frac{N_{i}-N_{i-1}}{k}P_{i}+\frac{k-N_{n(k)-1}}{k}P_{n(k)}+o_{x}(1)\qquad\mbox{ as }k\rightarrow\infty\label{eq:21}
\end{equation}
\end{prop}
\begin{proof}
Consider the expression \eqref{eq:dsc-time-average-distance} for
$\left\langle \mu,x\right\rangle _{k}$. We can ignore the error term,
since this is $o(1)$ as $k\rightarrow\infty$. the rest is an average
of distributions; our main task goal will be to compare the $i$-th
distribution to $P_{i}$. Fix $\varepsilon>0$ and for $x\in\mathbb{R}^{d}$
write $x_{i}=T_{\mathcal{D}_{b^{i}}}x$ as before. For each $i$ let
$I_{i}=((1+\delta_{i})N_{i-1},N_{i}]\cap\mathbb{N}$ and consider
the event 
\[
A_{i}=\left\{ x\,:\,\exists k\in I_{i}\quad d(\left\langle \nu_{i},x_{N_{i-1}}\right\rangle _{k},P_{i})>\varepsilon\right\} 
\]
We shall show that $\sum\mu(A_{i})<\infty$, and claim that this implies
the proposition because. The implication is proved as follows: by
Borel-Cantelli, it implies that for $\mu$-a.e $x\notin A_{i}$ there
is an $i_{0}$ with $x\notin A_{i}$ for $i\geq i_{0}$. This means
that in the expression \eqref{eq:dsc-time-average-distance} for $\left\langle \mu,x\right\rangle _{k}$,
the terms corresponding to $i_{0}\leq i\leq n-1$ are within $\varepsilon$
of $P_{i}$. Now, the total weight of the first $i_{0}-1$ terms is
$o(1)$ as $k\rightarrow\infty$. As for the last term, if $k<(1+\delta_{n})N_{n-1}$
then its weight is $\delta_{n}$, which by assumption is $o(1)$ as
$n\rightarrow\infty$; and if $k\geq(1+\delta_{n})N_{n}$ then that
term is also within $\varepsilon$ of $P_{n}$. Thus $\left\langle \mu,x\right\rangle _{k}$
is within $\varepsilon+o(1)$ of the expression in \eqref{eq:21},
and $\varepsilon$ is arbitrary.

We turn to the proof that $\sum\mu(A_{i})<\infty$. By \eqref{eq:18},
$\mu=\dsc(\mu_{i-1},\mu^{(i-1)},b^{N_{i-1}})$, so
\begin{equation}
\mu(A_{i})=\mu^{(i-1)}\left(y\,:\,\exists k\in I_{i}-N_{i-1}\quad d(\left\langle \nu_{i},y\right\rangle _{k},P_{i})>\varepsilon\right)\label{eq:22}
\end{equation}
We wish to compare the right hand side to the quantity $q_{i}(\varepsilon/2)$
in the statement. Indeed, $\mu^{(i-1)}$ and $\nu_{i}$ are nearly
the same, since $\mu^{(i-1)}=\dsc(\nu_{i},\mu^{(i)},b^{N_{i}-N_{i-1}})$.
More precisely, the cells in $\mathcal{D}_{b^{N_{i}-N_{i-1}}}$ have
the same mass under $\mu^{(i-1)}$ and $\nu_{i}$. We claim that the
property $d(\left\langle \nu_{i},y\right\rangle _{k},P_{i})<\varepsilon$
is almost entirely determined by the cell $\mathcal{D}_{b^{N_{i}-N_{i-1}}}$
of $y$ and not $y$ itself. Indeed, let $E=\mathcal{D}_{b^{N_{i}-N_{i-1}}}$,
let $y,y'\in E$, let 
\[
J_{i}=[\delta_{i}N_{i-1},N_{i}-N_{i-1}]\cap\mathbb{N}
\]
and let $k\in J_{i}$. The first (measure) component of $\widetilde{M}^{k}(\nu_{i},y)$
and $\widetilde{M}^{k}(\nu_{i},y')$ is the same, and the second components
are $T_{\mathcal{D}_{b^{k}}(y)}(y),T_{\mathcal{D}_{b^{k}}(y')}(y')$,
respectively. Since $y,y'$ are in the same cell $E\in\mathcal{D}_{b^{N_{i}-N_{i-1}}}$
we have $\mathcal{D}_{b^{N_{i}-N_{i-1}}}(y)=\mathcal{D}_{b^{N_{i}-N_{i-1}}}(y')$
and $\left\Vert y-y'\right\Vert <c\cdot b^{-(N_{i}-N_{i-1})}$, and
so $\left\Vert T_{\mathcal{D}_{b^{k}}(y)}(y)-T_{\mathcal{D}_{b^{k}}(y)}(y')\right\Vert <c\cdot b^{-(N_{i}-N_{i-1}-k)}$.
Taken together, this gives 
\begin{eqnarray}
\left\langle \nu_{i},y\right\rangle _{k} & = & \left\langle \nu_{i},y'\right\rangle _{k}+O(\frac{1}{k}\sum_{j=0}^{k-1}b^{-(N_{i}-N_{i-1}-j)})\nonumber \\
 & = & \left\langle \nu_{i},y'\right\rangle _{m}+O_{b}(\frac{b^{-(N_{i}-N_{i-1}-k)}}{k})\nonumber \\
 & = & \left\langle \nu_{i},y'\right\rangle _{m}+o(1)\qquad\mbox{as }i\rightarrow\infty\label{eq:18}
\end{eqnarray}
Thus as soon as $i$ is large enough that the $o(1)$ error term of
the last equation is $<\frac{1}{4}\varepsilon$, if $y$ is in the
event in \eqref{eq:22} then there is a $y'\in\mathcal{D}_{b^{N_{i}-N_{i-1}}}(y)$
such that $d(\left\langle \nu_{i},y'\right\rangle _{k},P_{i})>\frac{\varepsilon}{2}$
for some $k\in I_{i}-N_{i-1}$. Thus, 
\begin{eqnarray*}
\mu(A_{i}) & = & \sum\left\{ \mu^{(i-1)}(E)\,\left|\begin{array}{c}
E\in\mathcal{D}_{b^{N_{i}-N_{i-1}}}\mbox{ s.t. }d(\left\langle \nu_{i},y\right\rangle _{N_{i}-N_{i-1}},P_{i})>\varepsilon\\
\mbox{ for some }y\in E\;\mbox{and}\; k\in I_{i}-N_{i-1}
\end{array}\right.\right\} \\
 & \leq & \sum\left\{ \nu_{i}(E)\,\left|\begin{array}{c}
E\in\mathcal{D}_{b^{N_{i}-N_{i-1}}}\mbox{ s.t. }d(\left\langle \nu_{i},y'\right\rangle _{N_{i}-N_{i-1}},P_{i})>\frac{\varepsilon}{2}\\
\mbox{ for some }y\in E\;\mbox{and}\; k\in I_{i}-N_{i-1}
\end{array}\right.\right\} \\
 & = & \nu_{i}(y'\,:\, d(\left\langle \nu_{i},y\right\rangle _{k},P_{i})>\frac{\varepsilon}{4}\;\mbox{for some }\; k\in I_{i}-N_{i-1})\\
 & = & q_{i}(\frac{\varepsilon}{4})
\end{eqnarray*}
Since $\sum_{i}q_{i}(\frac{\varepsilon}{4})<\infty$, the proposition
is proved.\end{proof}
\begin{prop}
\label{prop:USM-construction}Let $\nu_{n},\delta_{n},N_{n}$ and
$\mu$ be as in the previous proposition. Suppose that for every $\varepsilon>0$
we have
\[
\limsup_{i\rightarrow\infty}\int\theta(U_{\varepsilon})\, dP_{i}=o(1)\qquad\qquad\mbox{as }\varepsilon\rightarrow0
\]
Then for $\mu$-a.e. $x$, if $\left\langle \mu,n\right\rangle _{k(j)}\rightarrow Q$
for some sequence $k(j)\rightarrow\infty$, then $\left\langle \mu\right\rangle _{x,k(j)\log b}^{\sqr}\rightarrow\cnt^{\sqr}Q$
and 
\begin{equation}
\cnt^{\sqr}Q=\lim_{j\rightarrow\infty}\sum_{i=1}^{n(k(j))-1}\frac{N_{i}-N_{i-1}}{k(j)}\cnt^{\sqr}P_{i}+\frac{k(j)-N_{n(k(j))-1}}{k(j)}\cnt^{\sqr}P_{n(k(j))}\label{26}
\end{equation}

\end{prop}
It remains to control the sceneries $\mu_{x,t}$ for $\mu$ as above.
We do this as follows. Given $\varepsilon>0$ let $U_{\varepsilon}$
denote the $\varepsilon$-neighborhood of the boundary of cells from
$\mathcal{D}_{b}$. 
\begin{proof}
By \eqref{eq:dsc-iteratived-distance} and approximating $1_{U_{\varepsilon}}$
by continuous functions in the usual way, our hypothesis, the assumption
(3) of the previous proposition and Borel-Cantelli imply that for
$\mu$-a.e. $x$, 
\[
\limsup_{i\rightarrow\infty}\sup_{\delta_{i}N_{i-1}\leq k\leq N_{i}-N_{i-1}}\left\langle \nu_{i},x\right\rangle _{k}(U_{\varepsilon})=o_{x}(1)\qquad\qquad\mbox{as }\varepsilon\rightarrow0
\]
Hence by the conclusion of the previous proposition, for $\mu$-a.e.
$x$, 
\begin{equation}
\limsup_{k\rightarrow\infty}\left\langle \mu,x\right\rangle _{k}(U_{\varepsilon})=o_{x}(1)\qquad\mbox{as }\varepsilon\rightarrow0\label{eq:vanishing-b-adic-boundaries}
\end{equation}
Following the argument in the proof of Proposition \ref{pro:USM-to-FD-and-CP}
(2) and (3), we see that for $\mu$-a.e. $x$, any accumulation point
$Q$ of $\left\langle \mu,x\right\rangle _{n}^{\sqr}$ is a CP-distribution,
and any accumulation point of $(\mu_{x,t}^{\sqr})$ is the centering
of such a $Q$ (the only difference between the proof of the proposition
and our situation now is that the rate of decay of the $o_{x}(1)$
term in \ref{eq:vanishing-b-adic-boundaries} is not explicit). Equation
\eqref{26} follows from the continuity of the centering operation
(again, see proof of Proposition \ref{pro:USM-to-FD-and-CP})
\end{proof}
We will also need the following variant of Proposition \ref{prop:splicing-measures}.
The conclusion of that proposition can be restated as follows: for
any $f\in C(\mathcal{M}^{\sqr}\times B_{1})$, 
\[
\frac{1}{k}\sum_{j=1}^{k}\int f(M^{j}(\mu,x))=\sum_{i=1}^{n(k)-1}\frac{N_{i}-N_{i-1}}{k}\int f\, dP_{i}+\frac{k-N_{k-1}}{k}\int f\, dP_{k}+o(1)\qquad\mbox{ as }k\rightarrow\infty
\]
We will sometimes want to vary $f$ depending on the term in the sum.
More precisely,
\begin{prop}
\label{prop:splicing-measures-2}Let $\nu_{n},P_{n},N_{n},\delta_{n}$
be as in Proposition \ref{prop:splicing-measures}, satisfying the
hypothesis (1)--(3). Let $f\in C(\mathcal{M}^{\sqr}\times B_{1})$,
$n=1,2,\ldots$, and suppose that \end{prop}
\begin{enumerate}
\item [(4)] For all $\varepsilon>0$, setting
\[
q'_{n}(\varepsilon)=\nu_{n}\left(y\,:\,\exists k\geq(1+\delta_{n})N_{n-1}\quad|\int f_{n}\, d\left\langle \nu_{n},y\right\rangle _{k}-\int f_{n}\, dP_{n}|>\varepsilon\right)
\]
we have $\sum q'_{n}(\varepsilon)<\infty$.
\end{enumerate}
Then for the measure $\mu$ constructed as above from $\nu_{n},N_{n}$,
for $\mu$-a.e. $x$ we have, as $k\rightarrow\infty$, 
\begin{equation}
\frac{1}{k}\sum_{j=0}^{k-1}f_{n(j)}(M^{j}(\mu,x))=\sum_{i=1}^{k-1}\frac{N_{i}-N_{i-1}}{k}\int f_{i}\, dP_{i}+\frac{k-N_{k-1}}{k}\int f_{k}\, dP_{k}+o(1)\label{eq:24}
\end{equation}

The proof is identical to that of the previous proposition, and we
omit it.

\subsection{\label{sub:USM-generating-non-ergodic-FD}A USM generating a non-ergodic
FD}

Let $d=1$ and $b=2$. We shall construct our examples using Proposition
\ref{prop:splicing-measures}. 

Let $N_{n}=n^{2}$ and $\delta_{n}=1-(1+\frac{1}{n})^{2}$. Then (1),(2)
of Proposition \eqref{prop:splicing-measures} are satisfied.

Let $P$ be the (ergodic) CP-distribution generated by $\lambda^{\sqr}$.
Let $x_{0}\in[-1,1]$ by a point chosen randomly according to $\lambda^{\sqr}$so
with probability one, $\delta_{x_{0}}$ generates the CP-distribution
$Q$ consisting of pairs $(\delta_{x},x)$, with $x\sim\lambda^{\sqr}$.
We fix such $x_{0}$. Define $\nu_{n}$ by
\[
\nu_{n}=\left\{ \begin{array}{cc}
\lambda^{\sqr} & n\mbox{ evem}\\
\delta_{x_{0}} & n\mbox{ odd}
\end{array}\right.
\]
so that in the notation of the proposition, $P_{2n}=P$ and $P_{2n+1}=Q$.
We must verify (3). The fact that it holds for odd $n$ is immediate:
as soon as $d(\left\langle \delta_{x_{0}},x_{0}\right\rangle _{n},Q)<\varepsilon$
we have $q_{n}(\varepsilon)=0$. As for even $n$, note that $M_{2}^{\sqr}(\lambda^{\sqr},x)=(\lambda^{\sqr},x')$
and so the speed of convergence of $\left\langle \lambda^{\sqr},x\right\rangle _{k}$
to $P$ is controlled by the speed of convergence of the second component.
Now, the fact that $2^{n}x\bmod1$, $n=1,2,\ldots$, equidistributes
exponentially fast for $\lambda^{[0,1]}$ as $n\rightarrow\infty$
follows from classical large deviations bounds, and, up to re-scaling
(from the intervals $[0,1]$ to $[-1,1]$) this is exactly what we
need. Thus $q_{n}(\varepsilon)$ decays exponentially for even $n$
and all $\varepsilon$, and (3) holds. 

Condition (4) of Proposition \ref{prop:splicing-measures-2} is also
clear for $P,Q$. Thus, the measure $\mu$ that we have constructed
is a USM, generating the centering of any accumulation point of the
expression in \eqref{eq:21}. From our choice of $N_{i}$ this expression
converges to $\frac{1}{2}P+\frac{1}{2}Q$, whose centering is $\frac{1}{2}\delta_{\lambda^{\sqr}}+\frac{1}{2}\delta_{\delta_{0}}$,
a non-ergodic FD.

We note that the construction is easily generalized to other dimensions.
We omit the details.

\subsection{\label{sub:examples}USMs whose projections misbehave}

In this section we exhibit a USM $\mu$ generating an EFD $P$, such
that $\ldim\pi\mu>E_{P}(\pi)$. The implication is that the inequality
in Theorem \ref{thm:projection-of-SMs} cannot in general be replaced
by an equality. We also sketch a construction showing that $\pi\mu$
need not be exact dimensional. 

Let $\nu$ be a self-similar measure which also is a product measure,
$\nu=\nu'\times\nu'$, and with $\dim\nu=2\dim\nu'<1$. For example
take the measure from Section \eqref{sub:non-conservation-example}.
Then $\nu$ generates an EFD $P$ with $\dim P=\dim\nu$, and, as
we saw in Section \eqref{sub:Homogeneous-measures}, $E_{\pi}(P)=\dim\pi\nu$
for every $\pi\in\Pi_{2,1}$. Specifically, for $\pi(x,y)=x$ we have
$\pi\nu=\nu'$ so $E_{\pi}(P)=\dim\pi\nu<\dim\nu$. 

Our aim is to construct a USM $\mu$ generating $P$, but such that
$\dim\pi\mu=\dim\nu>E_{P}(\pi)$. 

Let $R_{\theta}:\mathbb{R}^{2}\rightarrow\mathbb{R}^{2}$ denote rotation
by angle $\theta$. Choose a sequence $\theta_{n}\rightarrow0$ such
that $\dim\pi(R_{\theta_{n}}\nu)=\dim(P)$. We can do this because
the set of $\theta$'s for which this holds has full Lebesgue measure
by Theorem \ref{thm:Hunt-Kaloshin} (Marstrand's theorem for measures).
Let $P_{n}=R_{\theta_{n}}^{*}P$, which is an EFD by Proposition \eqref{prop:linear-images-of-FDs},
and observe that $P_{n}^{\sqr}=R_{\theta_{n}}^{\sqr}P\rightarrow P^{\sqr}$
(Basically this is because $R_{\theta_{n}}\rightarrow\id$. There
is an issue with the discontinuity of the $\sqr$ operation which
can be dealt in a similar manner as in Section \eqref{sub:Linear-images-of-FDs}.
We omit the details). Let $Q_{n}$ be a restricted base-$10$ CP-distribution
whose centering is $P_{n}$ and satisfies $\int\tau\, dQ_{n}(\tau)=\lambda^{\sqr}$
(Proposition \ref{pro:CP-distribution-in-arbitrary-coords}).

Let $\nu_{n}$ be a $Q_{n}$-typical measure, and construct $\mu$
from the sequence $\nu_{n}$ as in Section \eqref{sub:recipe-for-constructions}
using a rapidly growing sequence $N_{1},N_{2},\ldots$ satisfying
the assumptions of Proposition \eqref{prop:splicing-measures}. We
can do this since, fixing $\delta_{n}=\frac{1}{n}$, given $n$, condition
(3) of the proposition will hold if $\delta_{n}N_{n-1}$ is large
enough, which we can arrange at the $n-1$-th stage. Thus, $\mu$
is a USM. Furthermore the hypothesis of Proposition \ref{prop:USM-construction}
is satisfied, due to the fact that $\int\tau dQ_{n}(\tau)=\lambda^{\sqr}$.
Thus, for $\mu$-a.e. $x$, every accumulation point of $\left\langle \mu,x\right\rangle _{k}$
as $k\rightarrow\infty$ has a centering given by \eqref{26}. As
we observed already, $P_{n}^{\sqr}\rightarrow P^{\sqr}$, the centering
is $P$. Since this is true of any accumulation point, $\mu$ generates
$P$.

We next estimate $\ldim\pi\mu$ for $\pi(x,y)=x$. In Section \eqref{sub:Dimension-of-projected-measures}
we saw that if 
\[
\liminf_{N\rightarrow\infty}\frac{1}{N}\sum_{n=1}^{N}\frac{1}{\log b}H(\mu_{\mathcal{D}_{b^{n}}(x)},\mathcal{D}_{b^{n+1}})\geq\alpha
\]
then $\dim\pi\mu\geq\alpha$. A trivial modification of that proof
shows that if an integer sequence $m_{1}<m_{2}<\ldots$ satisfies
$\frac{m_{k+1}}{m_{k}}\rightarrow1$ and 
\begin{equation}
\liminf_{N\rightarrow\infty}\frac{1}{N}\sum_{k=1}^{N}\frac{1}{(m_{k+1}-m_{k})\log b}H(\mu_{\mathcal{D}_{b^{m_{k}}}(x)},\mathcal{D}_{b^{m_{k+1}}})\geq\alpha\label{eq:25}
\end{equation}
then $\dim\pi\mu\geq\alpha$. IT is this limit that we shall estimate. 

By \cite[Theorem 8.2]{HochmanShmerkin09}, for each $n$ we can choose
integers $u_{n}$ so that 
\[
E_{\pi}(Q_{n})>\frac{1}{u_{n}\log b}\int h_{b^{u_{n}}}(\pi\theta)\, dQ_{n}(\theta)-\frac{1}{n}
\]
where $h_{a}$ is a continuous approximation of $H(\cdot,\mathcal{D}_{a}^{1})$
as in Lemma \eqref{lem:continuous-approximation-of-entropy}. Let
\[
f_{n}(\theta,y)=\frac{1}{u_{n}\log b}h_{b^{u_{n}}}(\theta)
\]
Then by our assumption about $\nu_{n}$,
\[
E_{\pi}(Q_{n})>\lim_{k\rightarrow\infty}\int f_{n}d\left\langle \nu_{n},y\right\rangle _{k}-\frac{1}{n}\qquad\mbox{for }\nu_{n}\mbox{-a.e. }y
\]
If $\delta_{n}N_{n-1}$ is large enough then $\sup_{k>\delta_{n}N_{n-1}}|\int f_{n}d\left\langle \nu_{n},y\right\rangle _{k}-\int f_{b^{u_{n}}}dQ_{n}|<\frac{1}{n}$
for $y$ in a set of $\nu_{n}$-measure $<2^{-n}$. Clearly we can
choose $N_{n}$ which satisfy this. Then the assumption (4) of Proposition
\ref{prop:splicing-measures-2} is satisfied, hence the conclusion
\eqref{eq:24} of that proposition holds. Setting $m_{i}=u_{n(i)}$,
this implies that \eqref{eq:25} holds for $\mu$-a.e. $x$, with
$\alpha=\dim\nu$. Finally, although we do not have control over how
quickly the $u_{i}$ grow, we can make $n(i)$ grow arbitrarily slowly
by making $N_{i}$ grow rapidly, in particular we can ensure that
$n(i)/n(i-1)\rightarrow1$. Thus, for a suitable choice of $N_{i},\delta_{i}$,
we will have $\dim\pi\mu=\dim\nu>E_{\pi}(P)$, as desired. 

A variant of this constructions gives a USM $\mu$ whose $\pi$-projection
is not exact dimensional. The only change is that for even $n$ we
define the angle $\theta_{n}=0$ (and take $\theta_{n}$ as before
for odd $n$). Again, if $N_{i}$ grow rapidly enough, one can show
that the upper pointwise dimension of $\pi\mu$ is a.e. $\dim\nu$.
To control the lower pointwise dimension one ensures that 
\[
\liminf_{n\rightarrow\infty}\frac{1}{n\log b}H(\pi\mu,\mathcal{D}_{b^{n}})=\dim\pi\nu
\]
(since $\ldim\pi\mu$ is not more than this limit by Lemma \ref{lem:lower-dim-lower-boundes-entropy}).
We can do so because 
\[
H(\pi\mu,\mathcal{D}_{b^{n}})=H(\mu^{(i)},\mathcal{D}_{n-i})+O_{i}(1)
\]
and $H(\mu^{(i)},\mathcal{D}_{n-i})=H(\nu_{i},\mathcal{D}_{n-i})+O(1)$
if $n<N_{i+1}=N_{i}$. Thus if $N_{n-1}$ is chosen large enough at
each stage, $\pi\mu$ will not have exact dimension.

\bibliographystyle{plain}
\bibliography{bib}

\begin{thebibliography}{10}

\bibitem{AspbergEkstromPerssonSchmeling2013}
Magnus Aspenberg, Fredrik Ekstr{\"{o}}m, Tomas Persson, and J{\"{o}}rg
  Schmeling.
\newblock On the asymptotics of the scenery flow.
\newblock 2013.
\newblock preprint.

\bibitem{Bandt1997}
C.~Bandt.
\newblock Note on the axiomatic approach to self-similar random sets and
  measures.
\newblock {\em Mathematica Gottingensis}, (5):45--49, 1997.

\bibitem{Bandt1992}
Christoph Bandt.
\newblock The tangent distribution for self-similar measures.
\newblock In {\em Lecture at the Fifth Conference on Real Analysis}. 1992.

\bibitem{BedfordFisher96}
Tim Bedford and Albert~M. Fisher.
\newblock On the magnification of {C}antor sets and their limit models.
\newblock {\em Monatsh. Math.}, 121(1-2):11--40, 1996.

\bibitem{BedfordFisher97}
Tim Bedford and Albert~M. Fisher.
\newblock Ratio geometry, rigidity and the scenery process for hyperbolic
  {C}antor sets.
\newblock {\em Ergodic Theory Dynam. Systems}, 17(3):531--564, 1997.

\bibitem{BedfordFisherUrbanski02}
Tim Bedford, Albert~M. Fisher, and Mariusz Urba{\'n}ski.
\newblock The scenery flow for hyperbolic {J}ulia sets.
\newblock {\em Proc. London Math. Soc. (3)}, 85(2):467--492, 2002.

\bibitem{CoverThomas06}
Thomas~M. Cover and Joy~A. Thomas.
\newblock {\em Elements of information theory}.
\newblock Wiley-Interscience [John Wiley \& Sons], Hoboken, NJ, second edition,
  2006.

\bibitem{Falconer90}
Kenneth Falconer.
\newblock {\em Fractal geometry}.
\newblock John Wiley \& Sons Ltd., Chichester, 1990.
\newblock Mathematical foundations and applications.

\bibitem{Falconer97}
Kenneth Falconer.
\newblock {\em Techniques in fractal geometry}.
\newblock John Wiley \& Sons Ltd., Chichester, 1997.

\bibitem{FanLauRao02}
Ai-Hua Fan, Ka-Sing Lau, and Hui Rao.
\newblock Relationships between different dimensions of a measure.
\newblock {\em Monatsh. Math.}, 135(3):191--201, 2002.

\bibitem{Feller71}
William Feller.
\newblock {\em An introduction to probability theory and its applications.
  {V}ol. {II}.}
\newblock Second edition. John Wiley \& Sons Inc., New York, 1971.

\bibitem{Fisher04}
Albert~M. Fisher.
\newblock Small-scale structure via flows.
\newblock In {\em Fractal geometry and stochastics {III}}, volume~57 of {\em
  Progr. Probab.}, pages 59--78. Birkh\"auser, Basel, 2004.

\bibitem{Furstenberg70}
Harry Furstenberg.
\newblock Intersections of {C}antor sets and transversality of semigroups.
\newblock In {\em Problems in analysis ({S}ympos. {S}alomon {B}ochner,
  {P}rinceton {U}niv., {P}rinceton, {N}.{J}., 1969)}, pages 41--59. Princeton
  Univ. Press, Princeton, N.J., 1970.

\bibitem{Furstenberg22001}
Hillel Furstenberg.
\newblock {\em Private communication}, 2001.

\bibitem{Furstenberg08}
Hillel Furstenberg.
\newblock Ergodic fractal measures and dimension conservation.
\newblock {\em Ergodic Theory Dynam. Systems}, 28(2):405--422, 2008.

\bibitem{Gavish09}
Matan Gavish.
\newblock Measures with uniform scaling scenery.
\newblock {\em preprint}, 2009.

\bibitem{Glasner2003}
Eli Glasner.
\newblock {\em Ergodic theory via joinings}, volume 101 of {\em Mathematical
  Surveys and Monographs}.
\newblock American Mathematical Society, Providence, RI, 2003.

\bibitem{Graf95}
Siegfried Graf.
\newblock On {B}andt's tangential distribution for self-similar measures.
\newblock {\em Monatsh. Math.}, 120(3-4):223--246, 1995.

\bibitem{Hochman2010}
Michael Hochman.
\newblock Geometric rigidity of {$\times m$}-invariant measures.
\newblock {\em preprint}, 2010.

\bibitem{HochmanShmerkin09}
Michael Hochman and Pablo Shmerkin.
\newblock local entropy and dimension of projections.
\newblock {\em preprint}, 2009.

\bibitem{HuntKaloshin97}
Brian~R. Hunt and Vadim~Yu. Kaloshin.
\newblock How projections affect the dimension spectrum of fractal measures.
\newblock {\em Nonlinearity}, 10(5):1031--1046, 1997.

\bibitem{Kallenberg83}
Olav Kallenberg.
\newblock {\em Random measures}.
\newblock Akademie-Verlag, Berlin, third edition, 1983.

\bibitem{KriegMorters98}
Daniela Krieg and Peter M{\"o}rters.
\newblock Tangent measure distributions of hyperbolic {C}antor sets.
\newblock {\em Monatsh. Math.}, 126(4):313--328, 1998.

\bibitem{Mattila95}
Pertti Mattila.
\newblock {\em Geometry of sets and measures in {E}uclidean spaces}, volume~44
  of {\em Cambridge Studies in Advanced Mathematics}.
\newblock Cambridge University Press, Cambridge, 1995.
\newblock Fractals and rectifiability.

\bibitem{Morters98}
Peter M{\"o}rters.
\newblock Symmetry properties of average densities and tangent measure
  distributions of measures on the line.
\newblock {\em Adv. in Appl. Math.}, 21(1):146--179, 1998.

\bibitem{MortersPreiss98}
Peter M{\"o}rters and David Preiss.
\newblock Tangent measure distributions of fractal measures.
\newblock {\em Math. Ann.}, 312(1):53--93, 1998.

\bibitem{PatzschkeZahle90}
N.~Patzschke and U.~Z{\"a}hle.
\newblock Self-similar random measures. {IV}. {T}he recursive construction
  model of {F}alconer, {G}raf, and {M}auldin and {W}illiams.
\newblock {\em Math. Nachr.}, 149:285--302, 1990.

\bibitem{Patzschke04}
Norbert Patzschke.
\newblock The tangent measure distribution of self-conformal fractals.
\newblock {\em Monatsh. Math.}, 142(3):243--266, 2004.

\bibitem{PeresShmerkin08}
Yuval Peres and Pablo Shmerkin.
\newblock Resonance between {C}antor sets.
\newblock {\em Ergodic Theory Dynam. Systems}, 29(1):201--221, 2009.

\bibitem{Pesin97}
Yakov~B. Pesin.
\newblock {\em Dimension theory in dynamical systems}.
\newblock Chicago Lectures in Mathematics. University of Chicago Press,
  Chicago, IL, 1997.
\newblock Contemporary views and applications.

\bibitem{Varadarajan1963}
V.~S. Varadarajan.
\newblock Groups of automorphisms of {B}orel spaces.
\newblock {\em Trans. Amer. Math. Soc.}, 109:191--220, 1963.

\bibitem{Walters82}
Peter Walters.
\newblock {\em An introduction to ergodic theory}, volume~79 of {\em Graduate
  Texts in Mathematics}.
\newblock Springer-Verlag, New York, 1982.

\bibitem{Zahle88}
U.~Z{\"a}hle.
\newblock Self-similar random measures. {I}. {N}otion, carrying {H}ausdorff
  dimension, and hyperbolic distribution.
\newblock {\em Probab. Theory Related Fields}, 80(1):79--100, 1988.

\end{thebibliography}

\end{document}